\documentclass[a4paper, 12pt]{article}
\usepackage[utf8]{inputenc}
\usepackage{amssymb}\usepackage{multirow}
\usepackage{bbold}
\usepackage{xcolor}
\usepackage[T1]{fontenc}
\usepackage[utf8]{luainputenc}
\usepackage{amsmath}
\usepackage{bbm}
\usepackage{import}
\usepackage[english]{babel}
\usepackage{blindtext}
\usepackage[export]{adjustbox}
\usepackage{subcaption}
\usepackage{amssymb}
\usepackage{amsthm}
\usepackage{float}
\usepackage[ruled,vlined]{algorithm2e}
\usepackage[margin=2cm]{geometry}
\usepackage[authoryear]{natbib}
\usepackage{xcolor}
\usepackage{amsthm}
\usepackage{caption}
\usepackage{scalerel,stackengine}
\stackMath
\newcommand{\cov}{\text{Cov}}
\newcommand{\var}{\text{Var}}
\newcommand\reallywidehat[1]{%
\savestack{\tmpbox}{\stretchto{%
  \scaleto{%
    \scalerel*[\widthof{\ensuremath{#1}}]{\kern-.6pt\bigwedge\kern-.6pt}%
    {\rule[-\textheight/2]{1ex}{\textheight}}%WIDTH-LIMITED BIG WEDGE
  }{\textheight}% 
}{0.5ex}}%
\stackon[1pt]{#1}{\tmpbox}%
}
\theoremstyle{definition}
\newtheorem{definition}{Definition}[]
\newtheorem{proposition}{Proposition}[]

\usepackage{url}
\usepackage{enumitem}
\setlist[enumerate]{itemsep=0mm}
\usepackage{outlines}
\usepackage{yhmath}
\usepackage{amsthm,amsmath,amsfonts,amssymb}
\RequirePackage{mathtools}

\usepackage{xcolor}
\usepackage{bm,accents}
\usepackage{graphicx}
\usepackage{nicefrac}
\usepackage[super]{nth}

\usepackage{color}
\bibpunct{(}{)}{;}{a}{,}{,}
\def\spacingset#1{\renewcommand{\baselinestretch}%
	{#1}\small\normalsize} \spacingset{1.6}
\renewcommand{\baselinestretch}{1.4} 
\newtheorem{theorem}{Theorem}[]
\newtheorem{corollary}{Corollary}[]

\newtheorem{lemma}{Lemma}
\theoremstyle{remark}
\newtheorem{remark}{Remark}[]
\usepackage{authblk}

\title{Improved Estimators for Semi-supervised   High-dimensional Regression Model }
\author{Ilan Livne, David Azriel,  Yair Goldberg}
\affil{Technion - Israel Institute of Technology}
\date{\today}

\begin{document}
\maketitle

\begin{abstract}
  We study a linear high-dimensional regression model in a semi-supervised setting, where for many observations only the vector of covariates $X$ is given with no responses $Y$. We consider a linear regression model but do not make any sparsity assumptions on the vector of coefficients, and aim at estimating  $\var(Y|X)$. We propose an estimator, which is unbiased, consistent, and asymptotically normal. This estimator can be improved by adding zero-estimators arising from the unlabeled data. Adding zero-estimators does not affect the bias and potentially can reduce the variance.
We further illustrate our approach for other estimators, and present an algorithm that improves estimation for any given variance estimator. Our theoretical results are demonstrated in a simulation study.

\vspace{9pt}
  
\end{abstract}

\section{Introduction}

%\subsection{Motivation}
High-dimensional data analysis, where the number of predictors is larger than the sample size, is a topic of current interest. In such settings, an important goal is to estimate the signal level~$\tau^2$ and the  noise level $\sigma^2$, i.e., to quantify how much variation in the response variable can be explained by the predictors, versus how much of the variation is left unexplained. For example, in disease classification using DNA microarray data, where the number of potential predictors, say the genotypes, is enormous per each individual, one may wish to understand how disease risk is associated with genotype versus  environmental factors.\par
Estimating the signal and noise levels is important even in a low-dimensional setting.  In particular, a statistical model partitions the total variability of the response variable into two components: the variance of the fitted model $\tau^2$, and the variance of the residuals $\sigma^2$. 
This partition is at the heart of  techniques such as ANOVA and  linear regression, where $\tau^2$  and $\sigma^2$  might also be commonly referred to  as explained  versus unexplained variation, or between treatments  versus within treatments variation.
 % In ANOVA they are typically referred as between treatments variation versus within treatments variation.
 % technique of analysis of variance (ANOVA) consists of partitioning the total variability of a response variable into components that are associated with the different sources of that variability.
 Moreover, 
 %the  explained variance $\tau^2$ and the residual variance $\sigma^2$  are  important quantities in  other statistical problems:
 in model selection problems,  $\tau^2$ and $\sigma^2$  may be required for computing popular statistics, such as Cp, AIC, BIC and $R^2$. 
%In linear regression models, inference on individual coefficients require knowing $\sigma^2$.
%Estimating $\sigma^2$ is also important on its own when we wish to understand the variance decomposition of $Y$.
  Both $\tau^2$ and $\sigma^2$ are also closely  related to other important statistical problems, such as genetic heritability and signal detection. 
%In typical linear models, the residual variance estimator usually performs well and plays an important role in inference. However, it does not work in high dimensional data sets.
Hence, developing good estimators for these quantities  is a desirable goal. \par
When the number of covariates $p$ is much smaller than the number of observations $n$, and a linear model is assumed,  the ordinary least squares (henceforth, OLS) method provides us  straightforward estimators for $\tau^2$ and $\sigma^2$. 
However, when $p>n$, it becomes more challenging to perform inference on $\tau^2$ and $\sigma^2$   without further assumptions, such as  sparsity of the coefficients. In practice, the sparsity assumption may be unrealistic for some areas of interest. In this case, considering only a small number of significant coefficient can lead to biases and inaccuracies. One relevant example is the problem of missing heritability, i.e., the gap between heritability estimates from genome-wide-association-studies (GWAS) and the corresponding estimates from twin studies. For example, by 2010, GWAS studies had  identified a relatively  small number of covariates  that collectively explained around $5\%$ of the total variations in the trait \textit{height}, which is a small fraction  compared to  $80\%$ of the total variations that were explained by twin studies \citep{eichler2010missing}.   Identifying all the GWAS covariates affecting a  trait,  and  measuring  how much variation  they  capture, is believed to  bridge a significant fraction of the heritability  gap. With that in mind, methods that heavily rely  on the sparsity assumption may underestimate $\tau^2$ by their nature.
We show in this work that in the semi-supervised setting, in which for many observations only the covariates $X$ are given with no responses $Y$, one may consistently estimate the heritability without sparsity assumptions.
We use the term \textit{semi-supervised setting} to describe a setting in which the distribution of $X$ is known.  The setting where the distribution of $X$ is only partially known is not part of this work.

Estimating $\tau^2$ and $\sigma^2$ in a high-dimensional regression setting is generally a  challenging problem. 
 As mentioned above, the sparsity assumption, which means that only a relatively small number of predictors are relevant, plays an important role in this context. 
\cite{fan2012variance}
 introduced a refitted cross validation method for estimating $\sigma^2$.  Their method includes a two-staged procedure where a  variable-selection technique is performed in the first stage, and OLS is used to estimate $\sigma^2$ in the second stage. \cite{sun2012scaled}  introduced the scaled lasso algorithm that jointly estimates  the noise level and the regression coefficients by an iterative lasso procedure. Both  works  provide asymptotic distributional results  for their estimators and prove consistency under several assumptions including sparsity.
 In the context of heritability estimation, \cite{gorfine2017heritability} presented  the HERRA estimator, which is based on the above methods and is also applicable to time-to-event outcomes, in addition to  continuous or dichotomous outcomes.
 Another recent related work is   \cite{tony2020semisupervised} that considers, as we do here, a semi-supervised learning setting. In their work, Cai and Gue  proposed the CHIVE estimator of $\tau^2$, which integrates both labelled and unlabelled data and works well when the model is sparse.
 They characterize its  limiting distribution  and calculate confidence intervals for $\tau^2$. For more related works, see the literature review of \cite{tony2020semisupervised}.
 
Rather than assuming sparsity, or other structural assumptions on the coefficient vector $\beta$, a different approach for high-dimensional inference  is to assume some knowledge  about the covariates distribution. 
\cite{Dicker} uses the method-of-moments to develop several asymptotically-normal  estimators of $\tau^2$ and   $\sigma^2$, when  the covariates are assumed to be Gaussian.
\cite{schwartzman2019simple} proposed the  GWASH  estimator for estimating heritability, which   is essentially a modification of one of Dicker's estimators where  the columns of $X$ are standardized. Unlike Dicker, the GWASH estimator can also be computed from typical summary statistics, without accessing the original data.
 \cite{janson2017eigenprism}
 proposed the EigenPrism procedure to estimate $\tau^2$ and $\sigma^2$. Their method, which is based on singular value decomposition  and convex optimization techniques,  provides estimates and  confidence intervals for normal covariates.
%  In this paper we introduce a U-statistic naive estimator of $\tau^2$ and show that it is asymptotically normal and also that it is asymptotically equivalent to one of Dicker's estimators. 

In this paper we introduce a naive estimator of $\tau^2$ and show that it is asymptotically equivalent to Dicker's estimators when the covariates are normal, an assumption which is relaxed in this work. The naive  estimator is also a U-statistic and asymptotically normal.  U-statistics can  be typically used to obtain uniformly minimum variance unbiased estimators (UMVUE). However, when moments restrictions exist, U-statistics are no longer UMVUE, as shown by \cite{hoeffding1977some}. Under the assumed {semi-supervised setting}, the distribution of $X$ is known (and hence, moments of $X$ are known). Thus, the naive estimator is not UMVUE and it potentially can be improved.
We demonstrate how its variance  can be reduced  by using zero-estimators that incorporate the additional information from the unlabelled data.

%\subsection{ Main Contributions}
The contribution of this paper is threefold. First, we propose a novel approach for improving initial estimators of the signal level $\tau^2$ in the semi-supervised setting without assuming sparsity or normality of the covariates.
The key idea of this approach is to use zero-estimators that are  correlated with the initial estimator of $\tau^2$ in order to reduce variance without introducing extra bias.
Second, we define a new notion of optimality with respect to a linear family of zero-estimators. 
This allows us to suggest a necessary and sufficient condition for identifying optimal oracle-estimators. We use the term \textit{oracle} to point out that the specific coefficients that compose the optimal linear combination of zero-estimators  are dependent on the unknown parameters. 
Third, we suggest two estimators that  successfully improve initial estimators of $\tau^2$. 
We discuss in detail the improvement of the naive estimator and also apply our approach to other estimators. Thus, in fact, we provide an algorithm that has the potential to improve any given estimator of $\tau^2$.

%\subsection{Organization of the Paper}
The rest of this work in organized as follows. In Section \ref{Naive_section} we describe our setting  and introduce the naive estimator. In Section \ref{improv_of_Naive_sec} we introduce the zero-estimator approach and suggest a new notion of optimality  with respect to linear families of zero-estimators.  An optimal oracle estimator of $\tau^2$ is also presented.
In Section \ref{proposed_est_sec} we apply the zero-estimator approach to improve the naive estimator. We  then study some theoretical properties of the improved estimators. Simulation results are given in Section \ref{sim_res}. Section \ref{gener_es} demonstrates  how the zero-estimator approach can be  generalized to other estimators. A discussion is given in Section \ref{discuss}, while the proofs are provided in the Appendix.

\section{The Naive Estimator}\label{Naive_section}
\subsection{Preliminaries}
%{\bf Contents of the second section.}
We begin with describing our setting and assumptions.
Let $(X_1,Y_1),...,(X_n,Y_n)$ be i.i.d.\ observations drawn from some unknown distribution  where $X_i\in\mathbb{R}^p$ and $Y_i\in\mathbb{R} $. 
We consider a semi-supervised setting, where we have access to  infinite i.i.d.  observations of the covariates. Thus,   we essentially assume we know the covariate distribution.
%We view this assumption as a \textit{starting point}  for a more general framework, where one may have access to additional $N-n$  i.i.d\ observations  $X_{n+1},...,X_N$, in which $N$ is assumed to be finite. 
%Although this is beyond the scope of this work, our simulations suggest that our proposed method can be still useful even when $N$ is finite, i.e., when the joint distribution of the covariates is not  known exactly. \textcolor{red}{Ilan: did you run simulations with finite $N$? If not, then I think we need to remove the last sentence}
Notice that the assumption of known covariate distribution has already been presented and discussed in the context of high-dimension regression (e.g. \citealt{candes2017panning} and \citealt{janson2017eigenprism}) without using the term  ``semi-supervised learning''.
 
 For $i=1,\ldots,n$ we consider the the linear model 
\begin{equation}\label{linear_model}
Y_i=\beta^TX_i+\epsilon_i,    
\end{equation}
  where $E(\epsilon_i|X_i)=0$ and $E(\epsilon_i^2|X_i)=\sigma^2.$
We also assume that the intercept term is zero, which can be achieved in practice by centering the $Y$'s. 
Let $(X,Y)$ denote a generic observation and let
 $\sigma_Y^2$ denote the variance of $Y$. Notice that it can be decomposed into signal and noise components, 
 \begin{equation}\label{var_y_decompose}
	\sigma_{Y}^2 =\text{Var}(X^T\beta+\epsilon)
	 =\beta^T\cov(X)\beta+\text{Var}(\epsilon)
	 ={\beta ^T}{\bf{\Sigma}} \beta +\sigma^2,
	 \end{equation}
where $\text{Var}(\epsilon)=E(\epsilon^2)=\sigma^2$ and $\cov(X)=\bf{\Sigma}.$ 

The \textit{signal} component $\tau^2\equiv{\beta ^T}{\bf{\Sigma}} \beta$ can be thought of as the total variance explained by the best linear function of the covariates, while the \textit{noise} component $\sigma^2$ can be thought of as the variance left unexplained. 
We assume that $E(X)\equiv \mu$  are known and also that $\bf{\Sigma}$ is invertible. Therefore, we can apply the linear transformation
$X\mapsto{\bf\Sigma}^{-1/2}(X-~\mu)$ and assume w.l.o.g. that
$\mu=\textbf{0}$ and $\bf{\Sigma}=\textbf{I}.$
It follows by (\ref{var_y_decompose}) that  $\sigma_{Y}^2=\|\beta\|^{2}+\sigma^2$,
which implies that in order to evaluate $\sigma^2,$ it is enough to estimate both $\sigma_{Y}^2$ and $\|\beta\|^{2}$. The former can be easily evaluated from the sample, and the main challenge is to derive an estimator for $\|\beta\|^{2}$ in the high-dimensional setting.

\subsection{A Naive Estimator }
In order to find an unbiased estimator for $\|\beta\|^{2}=\sum_{j=1}^p \beta_j^2$ we first consider the estimation of $\beta_j^2$ for each $j$.  A straightforward approach   is given as follows:
Let $W_{ij}\equiv X_{ij}Y_i$ for $i=1,...,n$, and $j=1,...,p$. Notice that
  $$
  E\left( {{W_{ij}}} \right) = E\left( {{X_{ij}}{Y_i}} \right) = E\left[ {{X_{ij}}\left( {{\beta ^T}{X_i} + {\varepsilon _i}} \right)} \right] = {\beta _j},$$
   Now, since
   $\{E(W_{ij})\}^2=E(W_{ij}^2)-\text{Var}(W_{ij})$, a natural unbiased estimator for $\beta_j^2$ is
 \begin{equation}\label{beta_j_hat} 
 {\hat\beta_j^2}\equiv\frac{1}{n}\sum\limits_{i=1}^{n}W_{ij}^2-\frac{1}{n-1}\sum\limits_{i=1}^{n}(W_{ij}-\overline{W}_j)^2 = \frac{1}{n(n-1)}\sum_{i_1\ne i_2}^n W_{i_1j}W_{i_2j}, 
 \end{equation}
where $\overline{W}_j=\frac{1}{n}\sum_{i=1}^{n}W_{ij}$. Thus, unbiased estimates of $\tau^2\equiv\|\beta\|^{2}$ and $\sigma^2$ are given by
\begin{equation} \label{estimates}
 {\hat \tau ^2} = \sum\limits_{j = 1}^p {\hat \beta _j^2 = } \frac{1}{{n\left( {n - 1} \right)}}\sum\limits_{j = 1}^p {\sum\limits_{{i_1} \ne {i_2}}^n {{W_{{i_1}j}}{W_{{i_2}j}}} } , \qquad
\hat\sigma^2=\hat\sigma_{Y}^2-\hat{\tau}^2,
\end{equation}
where $\hat\sigma_{Y}^2=\frac{1}{n-1}\sum\limits_{i=1}^{n}(Y_i-\bar{Y})^2$.
We use the term \textit{ Naive} to describe $\hat\tau^2$ since its construction is relatively simple and straightforward. The  Naive estimator was also discussed  by  \cite{kong2018estimating}.   A similar estimator was proposed by \cite{Dicker}. 
Specifically, let
$$\hat \tau _{Dicker}^2 = \frac{{{{\left\| {{{\textbf{\emph{X}}}^T}{\bf{Y}}} \right\|}^2} - p{{\left\| {\bf{Y}} \right\|}^2}}}{{n\left( {n + 1} \right)}}$$
 where $\textbf{\emph{X}}$ is the $n \times p$ design matrix and ${\bf{Y}}=(Y_1,...,Y_n)^T.$ The following lemma shows that  $\hat\tau^2$ and $\hat \tau _{Dicker}^2$ are asymptotically equivalent under some conditions. 
 \begin{lemma}\label{asymptotic_normality_Naive}
Assume the linear model in \eqref{linear_model} and \({X_i}\mathop \sim\limits^{i.i.d} N\left( {\bf{0} ,\bf{I} } \right),\)  and that
$\epsilon_1,\dots, \epsilon_n \sim N(0,\sigma^2).$
When    $\tau^2+\sigma^2$  is bounded and $p/n$ converges to a constant, then, 
$$\sqrt{n}\left(\hat\tau^2-\hat\tau_{Dicker}^2\right)\overset{p}{\rightarrow} 0.$$
\end{lemma}
Note that in this paper we are interested in  a high-dimensional regression setting and therefore we study the limiting behaviour when $n$ and $p$ go together  to $\infty.$
Using Corollary 1 from \cite{Dicker}, which computes the asymptotic variance of $\hat \tau _{Dicker}^2$, and the above lemma, we obtain the following corollary. 
\begin{corollary}\label{normality_naive}
Under the assumptions of Lemma~\ref{asymptotic_normality_Naive}, 
$$
\sqrt n \left( {\frac{{{{\hat \tau }^2} - {\tau ^2}}}{\psi }} \right)\overset{D}{\rightarrow} N(0,1)\,,
$$
where 
\(\psi  = 2\left\{ {\left( {1 + \frac{p}{n}} \right){{\left( {\sigma ^2 + {\tau ^2}} \right)}^2} - \sigma^4 + 3{\tau ^4}} \right\}.\)
\end{corollary}
The variance of the naive estimator $\hat\tau^2$ under model \eqref{linear_model} (without assuming normality)  is given by the following proposition.
\begin{proposition}\label{var_naive1}
Assume model \eqref{linear_model} and additionally  that
$\beta^T\bf{A}\beta$ and $\|\textbf{A}\|_F^{2}$ are finite.   Then,
\begin{equation}\label{var_naive}
{\var} \left( {{{\hat \tau }^2}} \right) = \frac{{4\left( {n - 2} \right)}}{{n\left( {n - 1} \right)}}\left[ {{\beta ^T}{\bf{A}}\beta  - {{\left\| \beta  \right\|}^4}} \right] + \frac{2}{{n\left( {n - 1} \right)}}\left[ {\left\| {\bf{A}} \right\|_F^2 - {{\left\| \beta  \right\|}^4}} \right],
\end{equation}
where \({\bf{A}} = E\left( {{{\bf{W}}_i}{\bf{W}}_i^T} \right)\)
and $\|\textbf{A}\|_F^{2}$ denoted the Frobenius norm of $\textbf{A}.$
\end{proposition}

The following proposition shows that the naive estimator is consistent under some minimal assumptions.  
\begin{proposition}\label{consistency_naive}
Assume model \eqref{linear_model} and additionally that $\tau^2+\sigma^2=O(1)$  and 
 $\frac{\|A\|_F^{2}}{n^2}\rightarrow 0 $.
 Then, $\hat\tau^2$ is consistent.
Moreover, when  the columns of $\bf{X}$ are independent and both   $p/n$ and $E(X_{ij}^4)$ are bounded,  then $\frac{\|A\|_F^{2}}{n^2}\rightarrow 0$ holds and $\hat\tau^2$ is $\sqrt{n}$-consistent.
\end{proposition}

%%%%%%%%%%%%%%%%%%%%%%%%%%%%%%%%%%%%%%%%%%%%%%%%%%%%%%%%%%%%%%%%%%%%%%%%%%%%%%%%%%%%%%%%%%%%%%%%%%%%%%%%%%%%%%%%%%%%%%%%%%%%

\section{Oracle Estimator}\label{improv_of_Naive_sec}
%\textcolor{red}{ADD HERE AN INRO PARAGRAPH}
In this  section we introduce the zero-estimator approach and study how it can be used to improve the naive estimator. 
In Section \ref{zero_est_section} we present the zero-estimator approach and an illustration of this approach is given in Section \ref{illustration_Zero}.  Section \ref{opt_oracle_est_sec} introduces a new notion of optimality with respect to linear families for zero-estimators. We then find an optimal oracle estimator of $\tau^2$ and calculate its improvement over the naive estimator.
 \subsection{The Zero-Estimator Approach}\label{zero_est_section}
 %We are now ready to suggest our novel approach that incorporates the additional data that exist in the semi-supervised setting in order to improve the estimation of  $\tau^2$.
We describe the approach in general terms.
Consider a random variable $V \sim P$, where $P$ belongs to a family of distributions ${\cal P}$. Let $g(V)$ be a zero-estimator, i.e., $E_P[g(V)]=0$ for all $P \in {\cal P}$. Let $T(V)$ be an unbiased estimator of a certain quantity of interest $\theta$. Then, the statistic $U_c(V)$, defined by $U_c(V)=T(V) - cg(V)$ for a fixed  constant $c$, is also an unbiased estimator of $\theta$. 
%The correlation between the  Naive estimator $\hat{\tau}^2$ and unbiased estimators of zero is the key in characterizing the situations in which lower variance can be achieved. 
The variance of $U_c(V)$ is
\begin{equation} \label{eq2}
\text{Var}[U_c(V)]=\text{Var}[T(V)]+c^2\text{Var}[g(V)] - 2c\cdot\text{Cov}[T(V),g(V)].
\end{equation}
Minimizing $\text{Var}[U_c(V)]$ with  respect to $c$ yields the minimizer 
\begin{equation}\label{general_c}
c^*= \frac{\text{\cov}[T(V),g(V)]}{\var[g(V)]}.
\end{equation}
Notice that $\text{Cov}[T(V),g(V)]\neq 0$ implies $\text{Var}[U_{c^*}(V)]<\text{Var}(T(V))$. In other words, by combining a correlated unbiased estimator of zero with the initial unbiased estimator of $\theta$, one can lower the variance. Note that plugging $c^*$ in (\ref{eq2}) reveals  how much variance can be  potentially reduced, 
\begin{align}\label{variance_change}
\text{Var}[U_{c^*}(V)]=&
\var[T(V)]-[c^*]^2\var[g(V)]
\nonumber\\=&
\text{Var}[T(V)]-\frac{\{\text{Cov}[T(V),g(V)]\}^2}{\text{Var}[g(V])}
=(1-\rho^2)\text{Var}[T(V)],
\end{align}
where $\rho$ is the correlation coefficient between $T(V)$ and $g(V)$. Therefore, it is best to find an unbiased zero-estimator $g(V)$ which is highly correlated with $T(V)$, the initial unbiased estimator of $\theta$ . It is important to notice that $c^*$ is an unknown quantity and, therefore, $U_{c^*}$ is not a statistic. However,  in practice, one can estimate $c^*$ by some $\hat{c}^*$ and use the approximation $U_{\hat{c}^*}$ instead.

%We argue that this approach can be generalized to improve  different estimators for the signal and noise levels. We begin by introducing the  Naive estimator and then move to show how zero-estimator approach can be applied to improve it.

%\newpage

\subsection{Illustration of the Zero-Estimator Approach}\label{illustration_Zero}
The following example illustrates how the \textit{zero-estimator approach} can be applied to improve the  naive estimator $\hat\tau^2$ in the simple linear model setting.

\theoremstyle{definition}
\newtheorem{exmp}{Example}
\begin{exmp}[$p=1$]\label{exmp1}
Assume model \eqref{linear_model} with $X\sim N(0,1)$.  
By (\ref{variance_change}), we wish to find a zero-estimator $g(X)$ which is  correlated with $\hat\tau^2$.  
Consider the estimator $U_{c}=\hat\tau^2+cg(X)$, where $g(X)\equiv\frac{1}{n}\sum\limits_{i=1}^{n}(X_i^2-1)$ and $c$ is a fixed constant. The variance of $U_{c}$ is minimized by $c^*=-2\beta^2$ and one can verify that
 $\text{Var}(U_{c^*})=\text{Var}(\Hat{\tau}^2)-\frac{8}{n}\beta^4.$ For more details see Remark \ref{example_1} in the Appendix
 \end{exmp}

 The above example illustrates the potential of using  additional information that exists in the semi-supervised setting to lower the variance of the initial Naive estimator $\hat{\tau}^2$. However, it also raises the question: \textit{Can we achieve a lower variance by adding different zero-estimators? }
 One might  attempt to  reduce the variance by adding zero-estimators such as $g_k(X)\equiv\frac{1}{n}\sum\limits_{i=1}^{n}[X_i^k-E(X_i^k)]$, for $k>2$.
Surprisingly this attempt will fail. 
%Moreover, this result holds regardless of how $X$ is distributed.
Hence, the  unbiased oracle estimator of $\tau^2$,  $R\equiv\hat{\tau}^2-2\beta^2g(X)$, is  optimal with respect to 
%finite linear combinations of 
 zero-estimators of the form $g_k(X).$ %$\set*{\sum_{k=1}^{m}c_kg_k(X): c_k \in\mathbb{R} }_{m=1}^{\infty}$.
This unanticipated result motivated us to extend the idea of the optimal zero-estimator to a general regression setting of $p$ covariates.

\subsection{Optimal Oracle  Estimator }\label{opt_oracle_est_sec}
We now define a new oracle unbiased estimator of $\tau^2$ and prove that under some regularity assumptions this estimator is optimal with respect to a family of zero-estimators. Here, optimality means that the variance cannot be further reduced by including additional zero-estimators of \textit{that}  given family. We now specifically define our notion of optimality in a general setting.  
\begin{definition}
Let $T$ be an unbiased estimator of $\theta$ and let $g_1,g_2,...$ be a sequence of zero-estimators, i.e., $E_{\theta}(g_i)=0$ for $i\in \mathbb{N}$ and for all $\theta$. Let  $\mathcal{G} = \left\{ {\sum\limits_{k = 1}^m {{c_k}{g_k}:{c_k} \in \mathbb{R},m \in \mathbb{N}} } \right\}$ be a family of zero-estimators. 
For a zero-estimator $g^* \in \cal G$, we say that $R^* \equiv T + g^*$  is an \textit{optimal oracle  estimator (OOE) } of  $\theta$ with respect to $\cal G$, if
     $\var_{\theta}[R^*]=\text{Var}_{\theta}[T + g^*]\leqslant \text{Var}_{\theta}[T+g]$ for all $g\in\cal G$ and for all $\theta.$ 
\end{definition}
We use the term oracle since \({g^*} \equiv \sum\limits_{k = 1}^{m} {c_k^*{g_k}} \) for  some optimal coefficients \( {c_1^*,...,c_m^*} \),  which are a function   of the unknown parameter $\theta$.
The following theorem suggests a necessary and sufficient condition for obtaining an~OOE. 
\begin{theorem}\label{theorem1}
Let $\mathbf{g}_m=(g_1,...,g_m)^T$ be a vector of zero-estimators and assume the covariance matrix $M\equiv\var[\mathbf{g}_m]$ is positive definite for every $m$. Then, 
$R^*$ is an \textit{optimal oracle estimator}  (OOE) with respect to the family of zero-estimator $\cal G$ iff 
 $R^*$ is uncorrelated with  every zero-estimator $g\in\cal G$, i.e.,  $\cov_{\theta}[R^*,g]=0$ for all $g \in\cal G$ and for all $\theta$. 
\end{theorem}
Returning to our setting, define the following oracle estimator
\begin{equation}\label{OOE_est}
 {T_{oracle}} = {\hat \tau ^2} - 2\sum\limits_{j = 1}^p {\sum\limits_{j' = 1}^p {{\psi _{jj'}}} }, \end{equation}
where \({\psi _{jj'}} = {\beta _j}{\beta _{j'}}{h_{jj'}}\) and 
\(h_{jj'}^{} = \frac{1}{n}\sum\limits_{i = 1}^n {\left[ {{X_{ij}}{X_{ij'}} - E\left( {{X_{ij}}{X_{ij'}}} \right)} \right]}\),
and let the $\cal G$ be the 
 family of zero-estimators of the form
 $
 g_{k_1\ldots k_p}=\frac{1}{n}\sum_{i=1}^{n}[X_{i1}^{k_1}\cdot...\cdot X_{ip}^{k_p}-E(X_{i1}^{k_1}\cdot\ldots\cdot X_{ip}^{k_p})],
 $
 where $\left( {{k_1},...,{k_p}} \right) \in~{\left\{ {0,1,2,3,...} \right\}^p} \equiv \mathbb{N}_0^p.$
The following proposition shows that ${T_{oracle}}$ is an OOE with respect to ${\cal G}$. 
\begin{theorem}[General $p$]\label{oracle_p}
Assume model \eqref{linear_model} and additionally  that  $X$ has moments of all orders. Then, the oracle estimator $T_{oracle}$ defined in~\eqref{OOE_est} is an OOE of  $\tau^2$ with respect to~$\cal G.$
\end{theorem}
% \begin{remark}
% For the proof of Theorem \ref{oracle_p}, homoscedasticity  of $\epsilon$ is not required. 
% \end{remark}
\begin{remark}
For the proof of Theorem \ref{oracle_p} and Proposition \ref{var_naive1}, homoscedasticity  of $\epsilon$ is not required. 
\end{remark}
We now compute the variance reduction of  ${T_{oracle}}$ with respect to the naive estimator. 
The following statement is a corollary of  Proposition \ref{var_naive1}.
\begin{corollary}\label{V_T_orac}
Assume model  \eqref{linear_model}
and additionally that the columns of \({\textbf{\emph{X}}}\) are independent. Then,
\begin{equation}\label{var_T_oracle}
\var\left( {{T_{oracle}}} \right) = \var\left( {{{\hat \tau }^2}} \right) - \frac{4}{n}\left\{ {\sum\limits_{j = 1}^p {\beta _j^4\left[ {E\left( {X_{1j}^4 - 1} \right)} \right] + 2\sum\limits_{j \ne j'}^{} {\beta _j^2\beta _{j'}^2} } } \right\}.
\end{equation}
\end{corollary}
%\begin{remark}\label{relative_imp}
Moreover, in the special case where   \({X_i}\mathop \sim\limits^{i.i.d} N\left( {\bf{0} ,\bf{I} } \right)\). Then, Rewriting \eqref{var_T_oracle} yields 
\begin{equation}\label{var_T_oracle_normal}
\var\left( {{T_{oracle}}} \right) = \var\left( {{{\hat \tau }^2}} \right) - \frac{4}{n}\left\{ {2\sum\limits_{j = 1}^p {\beta _j^4 + 2\sum\limits_{j \ne j'}^{} {\beta _j^2\beta _{j'}^2} } } \right\} = \var\left( {{{\hat \tau }^2}} \right) - \frac{8}{n}{\tau ^4}.
\end{equation}
Notice that by Cauchy–Schwarz inequality, since  $E(X^2)=1$ then $E(X^4) \geq 1$, and therefore $\var(T_{oracle})<\var(\hat\tau^2).$ 
The following example provides intuition about the improvement of $\var(T_{oracle})$ over $\var(\hat\tau^2).$ 
\begin{exmp}[]\label{exp_OOE}
Consider a  setting where $n=p$; $\tau^2=\sigma^2=1$  and  \({X_i}\mathop \sim\limits^{i.i.d} N\left( {\bf{0} ,\bf{I} } \right)\).  
In this case, one can verify by \eqref{var_naive1} that $\var(\hat\tau^2)=\frac{20}{n}+O(n^{-2})$ and therefore $\var(T_{oracle})=\frac{12}{n} +O(n^{-2}).$ 
In other words: the optimal oracle estimator  $T_{oracle}$ reduces (asymptotically) the variance of the  naive estimator by $40\%$.
Moreover, when $p/n$ converges to zero, the reduction is $66\%$. See Remark \ref{rem:improve} in the  Appendix for more details about the relative improvement of the optimal oracle estimator.
\end{exmp}

%\end{remark}
%\begin{remark}
% The assumption that $X$ is Gaussian, given in Corollary \eqref{oracle_p}, is not required in order to prove that $\var(T_{oracle})<\var(\hat\tau^2).$ We use this assumption merely so simplify the calculations.
% \end{remark}

\section{ Proposed Estimators }\label{proposed_est_sec}
In this section we show how to use the zero-estimator approach to derive improved estimators over $\hat{\tau}^2$. 
In Section \ref{costs} we show that estimating all  $p^2$ optimal coefficients given in (\ref{OOE_est}) may introduce too much variance.  
Therefore,
Sections \ref{imp_single_zero_estimator} and
 \ref{selecting_covariates} introduce alternative methods to reduce the number of zero-estimators used in estimation.

\subsection{The cost of estimation}\label{costs}
The optimal oracle estimator defined in (\ref{OOE_est}) is based on adding $p^2$ zero-estimators.
Therefore, it is reasonable to suggest and study  the following estimator instead of the oracle one:
 $$
 T = {\hat \tau ^2} - 2\sum\limits_{j=1}^{p}\sum\limits_{j'=1}^{p}\hat\psi_{jj'},
 $$
 where \[\hat\psi_{jj'}= \frac{1}{{n\left( {n - 1} \right)\left( {n - 2} \right)}}\sum\limits_{{i_1} \ne {i_2} \ne {i_3}}^{} {{W_{{i_1}j}}{W_{{i_2}j'}}\left[ {{X_{{i_3}j}}{X_{{i_3}j'}} - E\left( {{X_{{i_3}j}}{X_{{i_3}j'}}} \right)} \right]}, \]
 is a U-statistics estimator of 
 \({\psi _{jj}} \equiv {\beta _j}{\beta _{j'}}{h_{jj'}}.\) Notice that   $E\left( \hat\psi_{jj'} \right) = 0$ and  that for $i_1 \neq i_2$ we have  $E(W_{i_1j}W_{i_2j'})=\beta_j\beta_{j'}$; thus, $T$ is an unbiased estimator of $\tau^2$ and  we wish to check it reduces the variance of naive estimator $\hat\tau^2$.
 This is described in the following proposition.
 \begin{proposition}\label{var_oracle}
 Assume model \eqref{linear_model} and additionally that $\tau^2+\sigma^2=O(1)$; $E(X_{ij}^4)\leq C$ for some positive constant $C,$  
%  \({X_i}\mathop \sim\limits^{i.i.d} N\left( {\bf{0} ,\bf{I} } \right)\)
and $p/n = O(1).$  Then, 
\begin{align}
\var\left( T \right)& = \var\left( {{T_{oracle}}} \right) + \frac{{8{p^2}{\sigma_Y^4}}}{{{n^3}}}+ O(n^{-2}) \nonumber\\
&= \var\left( {{{\hat \tau }^2}} \right) - \frac{4}{n}\left\{ {\sum\limits_{j = 1}^p {\beta _j^4\left[ {E\left( {X_{1j}^4 - 1} \right)} \right] + 2\sum\limits_{j \ne j'}^{} {\beta _j^2\beta _{j'}^2} } } \right\}
 + \frac{{8{p^2}{\sigma_Y^4}}}{{{n^3}}} + O(n^{-2}),\label{cost_of_estimation2}
\end{align}
where $\sigma_Y^2\equiv \tau^2+\sigma^2.$
\end{proposition}
Note that the second equation in \eqref{cost_of_estimation2} follows from  \eqref{var_T_oracle}.
To build some intuition, consider the case when \({X_i}\mathop \sim\limits^{i.i.d} N\left( {\bf{0} ,\bf{I} } \right)\) and $p=n.$ Then, the last equation can be rewritten~as  
\begin{equation}\label{cost_of_estimation}
\var\left( T \right) = \var\left( {{{\hat \tau }^2}} \right) + \frac{8}{n}\left( {2{\tau ^2}{\sigma ^2} + {\sigma ^4}} \right) +O(n^{-2}).
\end{equation}
 Notice that the  term $\frac{8}{n}\left( {2{\tau ^2}{\sigma ^2} + {\sigma ^4}} \right)$ in~(\ref{cost_of_estimation}) reflects the additional variability  that comes with the attempt at estimating all~$p^2$ optimal coefficients.
Therefore, the estimator $T$ fails to improve the naive estimator $\hat\tau^2$
and a similar result holds for $p/n\rightarrow c$ for some positive constant $c$.
Thus, alternative ways that improve the naive estimator are warranted, which are discussed next. 
%  Since estimating $p^2$ parameters is expensive in terms of additional variability, we need to find alternative ways to improve the naive estimator.
 
 \subsection{Improvement with a single zero-estimator}\label{imp_single_zero_estimator}
 A simple way to improve the naive estimator is by adding only a single zero-estimator. More specifically,  let $U_{c^*}=\hat\tau^2-c^*g_n$ where 
$c^*=\frac{\text{\cov}[\hat\tau^2,g_n]}{\var[g_n]}$ and $g_n$ is some zero-estimator. 
By  \eqref{variance_change} we have
\begin{equation}\label{var_single_coeff}
\text{Var}[U_{c^*}]=\text{Var}(\hat\tau^2)-\frac{\{\text{Cov}[\hat\tau^2,g_n]\}^2}{\text{Var}[g_n]}.
\end{equation}
Notice that $U_{c^*}$ is an oracle estimator and thus $c^*$ needs to be estimated in order to eventually construct a non-oracle estimator.
Let ${g_n} = \frac{1}{n}\sum\limits_{i = 1}^n {{g_i}}$ be the sample mean of some zero estimators $g_1,...,g_n$.
By  \eqref{general_c}, it can be shown that 
\begin{equation}\label{eq:c_star}
    {c^*} = \frac{{2\sum\limits_{j = 1}^p {{\beta _j}\theta_j} }}{{\var\left( {{g_i}} \right)}},
\end{equation} where
$\theta_j\equiv E(S_{ij})$ and
$S_{ij}=W_{ij}g_i$. Notice that  $\var\left( {{g_i}} \right)$  does not depend on $i$. 
Derivation of \eqref{eq:c_star} can be found in  Remark \ref{c_star_single} in the  Appendix. 
Here, we specifically chose  ${g_i} = \sum\limits_{j < j'}^{} {{X_{ij}}{X_{ij'}}}$
as it  worked well in the simulations but we do not argue that this is the best choice.
Let  \({T_{{c^*}}} = {\hat \tau ^2} - {c^*}{g_n}\) denote the oracle estimator for the specific choice of $g_n$, and where $c^*$ is given in \eqref{eq:c_star}. Notice that by \eqref{var_single_coeff} we have
\begin{equation}\label{var_single}
\var \left( {{T_{{c^*}}}} \right) = \var\left( {{{\hat \tau }^2}} \right) - \frac{{{{\left[ {2\sum\limits_{j = 1}^p {{\beta _j}}\theta_j } \right]}^2}}}{{n\var(g_i)}}.
\end{equation}
The following example demonstrates  the improvement of $\var(T_{c^*})$ over $\var(\hat\tau^2).$

\begin{exmp}[Example 2 -  continued]\label{exp_singel}
Consider a  setting where $n=p$; $\tau^2=\sigma^2=1$ ; \({X_i}\mathop \sim\limits^{i.i.d} N\left( {\bf{0} ,\bf{I} } \right)\) and  $\beta_j=\frac{1}{\sqrt{p}}$ for $j=1,...,p.$ Notice that this is an extreme non-sparse settings since the signal level $\tau^2$ is uniformly distributed across all $p$ covariates.
 In this case one can verify  that 
 %  \(\frac{{{{\left[ {2\sum\limits_{j = 1}^p {{\beta _j}E\left( {{S_{ij}}} \right)} } \right]}^2}}}{{n \cdot \var\left( {{g_i}} \right)}} = \frac{{{{\left[ {2{\left( {p - 1} \right)}} \right]}^2}}}{{n \cdot \left[ {p\left( {p - 1} \right)/2} \right]}} = \frac{8}{n} - O(n^{-2}).\)
% Therefore, by \eqref{var_single} we have
$\var(T_{c^*})= \frac{12}{n}+O(n^{-2}),$    
 which is approximately $40\%$ improvement over the naive estimator variance (asymptotically). For more details see Remark \ref{rem:improve_singel} in the  Appendix. 
\end{exmp}
In the view of \eqref{eq:c_star}, a straightforward U-statistic estimator for ${c^*}$ is
\begin{equation}\label{c_hat_star}
{\hat c^*} = \frac{{\frac{2}{{n\left( {n - 1} \right)}}\sum\limits_{{i_1} \ne {i_2}}^{} {\sum\limits_{j = 1}^p {{W_{{i_1}j}}{S_{{i_2}j}}} } }}{{\var\left( {{g_i}} \right)}},    
\end{equation}
 where $\var(g_i)$ is assumed known as it depends only on the marginal distribution of $X$. 
 Thus, we suggest the following estimator \begin{equation}\label{single_coeff_estimator}
{ T_{{\hat c^*}}} = {\hat \tau ^2} - {\hat c^*}{g_n}, 
\end{equation}
and prove that
 $T_{{c^*}}$ and $T_{{\hat c^*}}$ are asymptotically equivalent under some conditions.  \begin{proposition}\label{singel_asymptotic}
Assume model \eqref{linear_model} and additionally that $\tau^2+\sigma^2$ and  $p/n$ are  $O(1).$ Also, for every $j_1,j_2,j_3,j_4$ assume that  $E\left( {X_{1{j_1}}^2X_{1{j_2}}^2X_{1{j_3}}^2 X_{1{j_4}}^2} \right)$ is bounded and that the columns of the design matrix $\textbf{X}$ are independent. 
Then, \(\sqrt n \left[ {{T_{{c^*}}} - {{ T}_{{ \hat c^*}}}} \right]\overset{p}{\rightarrow}~0.\) 
\end{proposition}
We note that the requirement that the  columns of $\textbf{X}$  be independent can be relaxed to some form of weak dependence. 
% The second and third conditions are not trivial. The following corollary states a general result for which these conditions holds.  
% \begin{corollary}
% (MORE OR LESS): Assume model \eqref{linear_model} and additionally that $\tau^2+\sigma^2=O(1)$; $E(Y_i^8)=O(1)$ and 
%   $p/n = O(1).$ 
% Then, \(\frac{{{{\bf{e}}^T}{\bf{Ce}}}}{{{n^4}}} \) and  \(\frac{{\left\| {\bf{A}} \right\|_F^2}}{{{n^2}}} \) converge to $0.$
% \end{corollary}

 \subsection{Improvement by selecting small number of covariates}\label{selecting_covariates}
  Rather than using a single zero-estimator to improve the naive estimator, we   now consider  estimating a small number of coefficients of $T_{oracle}$. 
  Recall that $T_{oracle}$ is based on adding  $p^2$ zero estimators to the naive estimator. This estimation comes with high cost in terms of additional variability as shown is \eqref{cost_of_estimation}. Therefore, it is reasonable to use only a small number of zero estimators.
   Specifically, let  \({\bf{B}} \subset \left\{ {1,...,p} \right\}\) 
   be a fixed set of some indices such that
    \(\left| {\bf{B}} \right| \ll p\) and consider the  estimator
    \begin{equation}\label{only_small_coeff_est}
    {T_{\bf{B}}} = {\hat \tau ^2} - 2\sum\limits_{j,j' \in {\bf{B}}}^{} {{\hat\psi _{jj'}}}.     
    \end{equation}

% $$
% {T_{\bf{B}}} = {\hat \tau ^2} - \frac{2}{{n\left( {n - 1} \right)\left( {n - 2} \right)}}\sum\limits_{{i_1} \ne {i_2} \ne {i_3}}^{} {\sum\limits_{j,j' \in {\bf{B}}} {{W_{{i_1}j}}{W_{{i_2}j'}}\left[ {{X_{{i_3}j}}X_{{i_3}j'}^{} - E\left( {{X_{{i_3}j}}X_{{i_3}j'}^{}} \right)} \right]} }. 
% $$    
By the same argument as in Proposition   \ref{var_oracle} we now have
\begin{equation}\label{var_T_B}
{\var} \left( {{T_{\bf{B}}}} \right) = {\var} \left( {{{\hat \tau }^2}} \right) - \frac{4}{n}\left\{ {\sum\limits_{j \in {\bf{B}}}^{} {\beta _j^4\left[ {E\left( {X_{ij}^4} \right) - 1} \right] + 2\sum\limits_{j \ne j' \in {\bf{B}}}^{} {\beta _j^2\beta _{j'}^2} } } \right\}  + O\left( {{n^{ - 2}}} \right).
% \var\left( {{T_{\bf{B}}}} \right) =  {\var\left( {{{\hat \tau }^2}} \right) - \frac{8}{n}\tau _{\bf{B}}^4}  +O(n^{-2}),
\end{equation}
Also notice that  when \({X_i}\mathop \sim\limits^{i.i.d} N\left( {\bf{0} ,\bf{I} } \right)\),  \eqref{var_T_B} can be rewritten as
\begin{equation}\label{var_T_B_norm}
\var\left( {{T_{\bf{B}}}} \right) =  {\var\left( {{{\hat \tau }^2}} \right) - \frac{8}{n}\tau _{\bf{B}}^4}  +O(n^{-2}).
\end{equation}
 where $\tau _{\bf{B}}^2 = \sum\limits_{j \in {\bf{B}}}^{} {\beta _j^2}.$
Thus, if $\tau _{\bf{B}}^2$ is sufficiently large,  one can expect a significant improvement over the naive estimator by  using a small number of zero-estimators.
For example,  when $\tau _{\bf{B}}^2=0.5$; $p=n$;  $\tau^2=\sigma^2=1$, then $T_{\bf{B}}$  reduces the $\var(\hat\tau^2)$  by $10\%$. For more details see  Remark \ref{rem:improve_T_B} in the Appendix.

%  Applying a covariates selection method that chooses wisely the covariates that will be included in $\bf{B}.$
%  Although it
%  and and we wish to apply a method that chooses wisely the covariates that will be in included in $\bf{B}.$ 
%  The problem of is not a primary focus of this work
% Different covariate selection methods have different  pros and cons and there is no simple answer for the question of which  method to use in order to construct a desirable set $\bf{B}.$ 

% This problem of high-dimensional covariate selection has already been already extensively  studied. 
% % \cite{bunea2006consistent}  studied the relationship between multiple hypotheses  testing procedures and consistent covariate selection in high-dimensional linear regression setting.
% \cite{bunea2006consistent}  studied how  multiple hypotheses  testing procedures can be used for consistent covariate selection in high-dimensional  setting.
% \cite{wasserman2009high} suggested a consistent multi-stage method  for covariate selection that successfully controls the false positive error while giving a reasonable power. 
% \cite{fan2010selective} discussed the development of theory,  methods and their applications,  for high-dimensional covariate selection.   
% More recent  works about this subject include \cite{zambom2018consistent}  and \cite{oda2020fast}.
% % Note that the problem of selecting  covariates that are associated with significantly large $\beta_j$ has already been extensively studied (See 1,2,3,4....).
% %  specific covariate selection method is not our main  interest in this work,

 Notice that we do not assume  sparsity of the coefficients. 
 The sparsity assumption essentially ignores covariates that do not belong to the set $\bf{B}$. 
When $\beta_j$'s for $j \notin {\bf B}$ contribute much to the signal level $\tau^2\equiv\|\beta\|^{2}$, the sparse approach leads to disregarding a significant portion of the signal, while our estimators do account for this as {all} $p$ covariates are used in $\hat{\tau}^2$. 
%  Specifically,  the estimator  $T_{\textbf{B}}$ in \eqref{only_small_coeff_est} contains two terms: the naive estimator term, which is based on  \textit{all} $p$ indices,  and a zero-estimator term which is based on  the covariates in $ {\bf{B}}$. Moreover, notice that the estimator $T_{\hat c^*}$ in \eqref{single_coeff_estimator} also uses all $p$ covariates.

%See Appendix for  more details. 
%  \textcolor{red}{Where in the Appendix? I couldn't find it. You can give a reference like: See Remark 5 in the appendix (don't write See Appendix)}
The following example illustrates some key aspects of our proposed estimators. 
% The following example gives some intuition about the performance of $T_{c^*}$ and $T_{\textbf{B}}.$
\begin{exmp}[Example 3 -  continued]\label{extreme_example}
Let
 $n=p$; $\tau^2=\sigma^2=1$ and \({X_i}\mathop \sim\limits^{i.i.d} N\left( {\bf{0} ,\bf{I} } \right)\). 
Consider the following two extreme scenarios:
\begin{enumerate}
    \item \textit{non-sparse setting}: The signal level $\tau^2$ is uniformly distributed over all $p$ covariates, i.e., $\beta_j^2=\frac{1}{p}$ for all $j=1,...,p.$  
    \item \textit{Sparse setting}: the signal level $\tau^2$ is "point mass" distributed over the set \textbf{B}, i.e., $\tau^2_{\textbf{B}}=\tau^2.$\
\end{enumerate}
Two interesting key points:
\begin{enumerate}
    \item In the first scenario the estimator $T_{\textbf{B}}$ has the same asymptotic variance as $\hat{\tau}^2$, while the estimator $T_{c^*}$ reduces the variance by approximately $40\%$.
 \item In the second scenario the variance reduction of  $T_{\textbf{B}}$ is approximately $40\%$, while $T_{c^*}$ has the same asymptotic variance as $\hat{\tau}^2$.
 \end{enumerate}
Interestingly, in this example, the OOE estimator $T_{oracle}$  asymptotically improves  the naive by $40\%$ regardless of the scenario choice,  as shown by \eqref{var_T_oracle_normal}. For more details see Remark \ref{summary_eqample}  in the Appendix.
\end{exmp}

 %and therefore leads to disregarding a significant portion of the signal  $\|\beta\|^{2}$.
 %   This is different from our setting where we wish to identify and then use both the "important" and "unimportant" covariates in our suggested estimators. 
  %Moreover, indices that do not belong to the set $\bf{B}$ are still important to maintain unbiasedness.  
 % Indeed,  

 A desirable set of indices $\bf B$ contains  relatively small amount of covariates  that capture  a significant part of the signal level $\tau^2.$
 There are different methods to choose the covariates that will be included in $\bf{B},$ but these are not a primary focus of this work. 
 For more information about covariate selection methods see  \cite{zambom2018consistent} and \cite{oda2020fast} and references therein. In Section \ref{sim_res} below we work with a certain selection algorithm defined there.
%In practice, one may need to select the covariates to include in the set $\textbf{B}$. This is discussed next.
We call $\delta$ a covariate \textit{selection algorithm} if for every dataset \(\left( {{{\textbf{\emph{X}}}_{n \times p}},{{\bf{Y}}_{n \times 1}}} \right)\) it chooses a subset of indices \({{\bf{B}}_\delta }\) from \(\left\{ {1,...,p} \right\}\).
% We say that $\delta$ is  \textit{consistent}  if it satisfies
% $P\left( {{{\bf{B}}_\delta} = {\bf{B}}} \right)\xrightarrow[n\rightarrow\infty]{}~1,$
% where \({{\bf{B}}_\delta }\) is the set of selected indices by $\delta.$
%   Notice that the selection algorithm $\delta$  may depend on both $n$ and $p$ but this is suppressed in the notation.
 Our proposed estimator  for $\tau^2$, which is based on selecting small number of covariates, is given in Algorithm~1.\\

\begin{algorithm}[H]
 \caption{Proposed Estimator based on covariate selection }\label{selection_estimator}
\SetAlgoLined

\vspace{0.4 cm}

\textbf{Input:}
 A dataset \(\left( {{{\bf{X}}_{n \times p}},{{\bf{Y}}_{n \times 1}}} \right)\) and a selection algorithm $\gamma$.
 \begin{enumerate}
 \item Calculate the naive estimator \({\hat \tau ^2} = \frac{1}{{n\left( {n - 1} \right)}}\sum\limits_{j = 1}^p {\sum\limits_{{i_1} \ne {i_2}}^n {{W_{{i_1}j}}{W_{{i_2}j}}} } \), where $W_{ij}=X_{ij}Y_i.$
     
     \item Apply algorithm $\gamma$ to  \(\left( {{{\bf{X}}},{{\bf{Y}}}} \right)\) to construct \({{\bf{B}}_{\gamma}}.\) 
     \item  Calculate the zero-estimator terms:  
     \[\hat\psi_{jj'}\equiv\frac{2}{n(n-1)(n-2)} {\sum\limits_{{i_1} \ne {i_2} \ne {i_3}}^{} {{W_{{i_1}j}}{W_{{i_2}j'}}\left[ {{X_{{i_3}j}}{X_{{i_3}j'}} - E\left( {{X_{{i_3}j}}{X_{{i_3}j'}}} \right)} \right]} }, \] for all \(j,j' \in {{\bf{B}}_{{\gamma}}}.\)
 \end{enumerate}
\KwResult{Return \(T_{\gamma} = {\hat \tau ^2} - \sum\limits_{jj' \in {{\bf{B}}_{\gamma}}}^{} {{\hat\psi _{jj'}}}. \)
}
\end{algorithm}

Some asymptotic properties   of $T_{\gamma}$ are given by the following proposition.
% By the same arguments as in Proposition \ref{var_EOOE},
\begin{proposition}\label{limit_of_proposed}
Assume there is a set 
\({\bf{B}} \equiv~\left\{ {j: {{\beta _j^2}}  > b} \right\}\) where $b$ is a positive constant, such that
  \(\left| {\bf{B}}  \right|=p_0 \) where $p_0$ is a fixed constant. 
  Also assume that $\mathop {\lim }\limits_{n\to \infty } n\left[ { P\left( \left\{ {{{\bf{B}}_\gamma } \neq {\bf{B}}} \right\} \right)} \right]^{1/2} = 0,$ and that $E\left( {T_\gamma ^4} \right)$ and $E(T_{\bf{B}}^4)$ are bounded.  Then,
%  \begin{enumerate}
%   \item $\mathop {\lim }\limits_{n \to \infty } n\left[Bias\left( {{T_\gamma }},\tau^2 \right) \right]^2 = 0,$
%   and
%   \item \(\mathop {\lim }\limits_{n \to \infty } n\left[ {\var} \left( {{T_\gamma }} \right)- \var\left( {{T_{\bf{B}}}} \right)   \right] = 0,\)
% \end{enumerate}
% where $Bias\left( {{T_\gamma }},\tau^2 \right)\equiv E(T_\gamma)-\tau^2$ and
% $ \var\left( {{T_{\bf{B}}}} \right)$ is given by \eqref{var_T_B}. 
$$\sqrt{n}(T_{{\gamma}}-T_{\bf{B}})\overset{p}{\rightarrow} 0.$$
\end{proposition} 
Notice that the requirement $\mathop {\lim }\limits_{n\to \infty } n\left[ { P\left( \left\{ {{{\bf{B}}_\gamma } \neq {\bf{B}}} \right\} \right)} \right]^{1/2} = 0$ is stronger than just consistency. 
% The following claim is a direct corollary of Proposition \ref{limit_of_proposed}. 
% \begin{corollary}\label{cor_MSE}
% Under the assumptions of Proposition~\ref{limit_of_proposed},

% $$
% \mathop {\lim }\limits_{n \to \infty }n\left[ {MSE\left( {{T_{{\gamma}}}} \right) - MSE\left( {{T_{\bf{B}}}} \right)} \right] = 0 \text{ and }\sqrt{n}(T_{{\gamma}}-T_{\bf{B}})\overset{p}{\rightarrow} 0. $$
%\end{corollary}

\begin{remark}[Practical considirations]
%\subsection{Practical Considerations }
Some cautions regarding  the estimator $T_{\gamma}$ need to be considered in practice. When $n$ is insufficiently large, then \({{\bf{B}}_{{\gamma}}}\) might be different than \({\bf{B}}\) and  Proposition \ref{limit_of_proposed} no longer holds.
 Specifically, let \({\bf{S}} \cap {{\bf{B}}_{\gamma} }\) and \({\bf{B}} \cap {{\bf{S}}_{\gamma} }\)   be the set of \textit{false positive} and \textit{false negative} errors, respectively, where \({\bf{S}} = \left\{ {1,...,p} \right\}\backslash {\bf{B}}\) and    \({{\bf{S}}_{\gamma} } = \left\{ {1,...,p} \right\}\backslash {{\bf{B}}_{\gamma }}.\) 
 While false negatives merely result in not including some potential zero-estimator terms in our proposed estimator, false positives can lead to a  substantial bias. This is true since the expected value of a post-selected zero-estimator is not necessarily  zero anymore. 
A common approach to overcome this problem is
to randomly \textit{split}  the data into two parts where the first part is used for covariate selection and the second part is used for  evaluation of the zero-estimator terms.
\end{remark}

\subsection{Estimating the variance of the proposed estimators}
We now suggest  estimators for   $\var(\hat\tau^2$), $\var(T_{\gamma})$  and $\var(T_{\hat c^*}).$ 
Let
\begin{equation*}
\widehat {\var\left( {{{\hat \tau }^2}} \right)} = \frac{4}{n}\left[ {\frac{{\left( {n - 2} \right)}}{{\left( {n - 1} \right)}}\left[ {\hat \sigma_Y^2{{\hat \tau }^2} + {{\hat \tau }^4}} \right] + \frac{1}{{2\left( {n - 1} \right)}}\left( {p{{\hat \sigma_Y^4}} + 4\hat \sigma_Y^2{{\hat \tau }^2} + 3{{\hat \tau }^4}} \right)} \right],   
\end{equation*}
  where \(\hat \sigma_Y^2 = \frac{1}{{n - 1}}\sum\limits_{i = 1}^n {{{\left( {{Y_i} - \bar Y} \right)}^2}} \), and \(\hat \sigma _Y^4 = {\left( {\hat \sigma _Y^2} \right)^2}.\) 
The following proposition shows that \(\widehat {\var\left( {{{\hat \tau }^2}} \right)}\) is consistent under some conditions.
\begin{proposition}\label{var_naive_est}
Assume model \eqref{linear_model} and additionally that $\tau^2+\sigma^2=O(1)$,  \({X_i}\mathop \sim\limits^{i.i.d} N\left( {\bf{0} ,\bf{I} } \right)\) and $p/n = O(1).$  Then,
 $$n\left[ {\widehat {\var\left( {{{\hat \tau }^2}} \right)} - \var\left( {{\hat\tau ^2}} \right)} \right] \overset{p}{\rightarrow} 0.$$ 
 \end{proposition}

Consider now $\var(T_{\gamma})$ and let $\widehat {\var\left( {{T_{\gamma}}} \right)} =  {\widehat {\var\left( {{{\hat \tau }^2}} \right)} - \frac{8}{n}\hat \tau _{{{\bf{B}}_{\gamma}}}^4},$     
where  \(\hat \tau _{{{\bf{B}}_{\gamma}}}^2 = \sum\limits_{j \in {{\bf{B}}_{\gamma}}}^{} {\hat \beta _j^2} \) and \(\hat \tau _{{{\bf{B}}_\gamma }}^4 = {\left( {\hat \tau _{{{\bf{B}}_\gamma }}^2} \right)^2}.\) 
  The following propositions shows that 
 \(\widehat {\var\left( {{T_\gamma }} \right)}\) is consistent.
\begin{proposition}\label{consist_var_}
Under the assumptions of Propositions \ref{limit_of_proposed} and \ref{var_naive_est}, $$
 n\left[ {\widehat {\var\left( {T_{\gamma}} \right)} - \var\left( {T_{\gamma}} \right)} \right] \overset{p}{\rightarrow}~0.
$$   
 \end{proposition}

When normality of the covariates is not assumed, we suggest the following estimators: 
\begin{equation}\label{var_naive_estimate}
 \widetilde {{\var} \left( {{{\hat \tau }^2}} \right)} = \frac{{4\left( {n - 2} \right)}}{{n\left( {n - 1} \right)}}\left[ {\widehat {{\beta ^T}{\bf{A}}\beta } - \widehat {{{\left\| \beta  \right\|}^4}}} \right] + \frac{2}{{n\left( {n - 1} \right)}}\left[ {\widehat {\left\| {\bf{A}} \right\|_F^2} - \widehat {{{\left\| \beta  \right\|}^4}}} \right];  
\end{equation}
   \begin{equation}\label{var_T_gamm}
\widetilde {\var\left( {{T_\gamma }} \right)} = \widetilde {\var\left( {\hat \tau } \right)} - \frac{4}{n}\left\{ {\sum\limits_{j \in {{\bf{B}}_\gamma }}^p {\hat \beta _j^4\left[ {E\left( {X_{1j}^4 } \right)- 1} \right] + 2\sum\limits_{j \ne j' \in {{\bf{B}}_\gamma }} {\hat \beta _j^2\hat \beta _{j'}^2} } } \right\};
    \end{equation}
and
$$\widetilde {\var\left( {{T_{{{\hat c}^*}}}} \right)} = \widetilde {\var\left( {{{\hat \tau }^2}} \right)} - \frac{{{{\left[ {\frac{2}{{n\left( {n - 1} \right)}}\sum\limits_{{i_1} \ne {i_2}}^{} {\sum\limits_{j = 1}^p {{W_{{i_1}j}}{S_{{i_2}j}}} } } \right]}^2}}}{{\var\left( {{g_i}} \right)}},$$
where $\widehat {{\beta ^T}{\bf{A}}\beta } = \frac{1}{{n\left( {n - 1} \right)\left( {n - 2} \right)}}\sum\limits_{{i_1} = {i_2} \ne {i_3}}^{} {} {{\bf{W}}_{{i_1}}}\left( {{{\bf{W}}_{{i_2}}}{\bf{W}}_{{i_2}}^T} \right){{\bf{W}}_{{i_3}}}; $
$\widehat {\left\| {\bf{A}} \right\|_F^2} = \frac{1}{{n\left( {n - 1} \right)}}\sum\limits_{{i_1} \ne {i_2}}^{} {{{\left( {{\bf{W}}_{{i_1}}^T{{\bf{W}}_{{i_2}}}} \right)}^2}};$   \newline
$\widehat {{{\left\| \beta  \right\|}^4}} = {( {\frac{1}{{n\left( {n - 1} \right)}}\sum\limits_{{i_1} \ne {i_2}}^{} {{\bf{W}}_{{i_1}}^T{{\bf{W}}_{{i_2}}}} } )^2}$ are all U-statistics estimators, and  $\hat\beta_j^2$   is given by  \eqref{beta_j_hat}.  
Although we do not provide here  formal proofs, our simulations support that these estimators are consistent under the same assumptions of Proposition \ref{var_oracle}.

\section{Simulations Results} \label{sim_res}
We now provide a simulation study to  illustrate our  estimators performance. We compare the different estimators that were discussed earlier in this work:
\begin{itemize}
    \item The naive estimator $\hat\tau^2$ which is given in \eqref{estimates}.
    \item The optimal oracle estimator $T_{oracle}$ which is given in
\eqref{OOE_est}. 
    \item The estimator $T_{{\hat c^*}}$ which is based on adding a single zero-estimator  and is given in \eqref{single_coeff_estimator}.
    \item The estimator $T_{\gamma}$ which is based on selecting a small number of covariates  and is given by Algorithm \ref{selection_estimator}. Details about the specific selection algorithm we used can be found in Remark \ref{selection_algorithm} in  Appendix.
\end{itemize}
An additional estimator we include in the simulation study is the PSI (Post Selective Inference), which was calculated using the {\fontfamily{qcr}\selectfont
estimateSigma} function from the {\fontfamily{qcr}\selectfont
selectiveInference} R package. The PSI estimator is based on the LASSO method which assumes sparsity of the coefficients and therefore ignores small coefficients.  

\begin{figure}[H]
  \centering
 \includegraphics[width=0.9\textwidth]{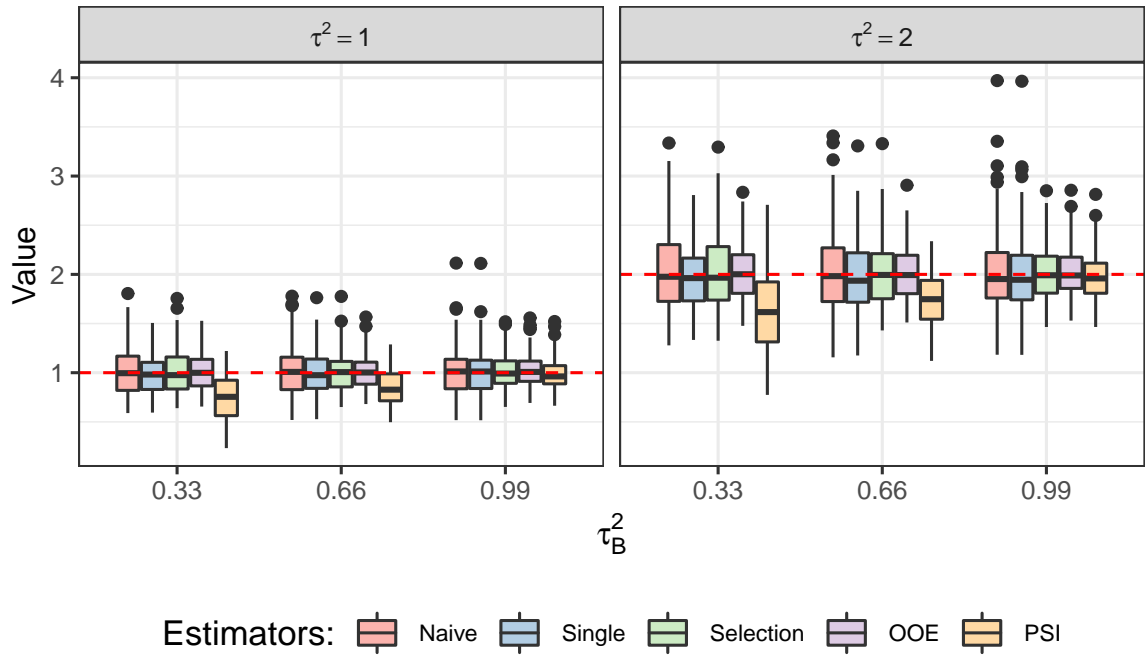}
\caption{
Boxplots representing the estimators distribution . The x-axis stands for $\tau^2_{\bf{B}}$. The red dashed is the true value of $\tau^2.$ 
% Boxplots for table \ref{table1}. The value of $\tau^2$ is shown by the red dashed line. The  naive, modified, proposed, oracle and the PSI estimators  are compared for sample sizes $n=100, 300, 500.$
}
\label{figure1}
\end{figure}

We fix  
$\beta_j^2=\frac{{\tau _{\bf{B}}^2}}{5}$ for $j=1,\dots,5$, and \(\beta _j^2 = \frac{{\tau^2 - \tau _{\bf{B}}^2}}{{p - 5}}\) for $j=6,\dots,p$, where $\tau^2$ and \(\tau _{\bf{B}}^2\) vary among different scenarios. The number of observations and covariates is $n=p=400$, and the residual variance is $\sigma^2=1$. 
For each scenario, we generated  100  independent datasets and estimated $\tau^2$ by using the different estimators. 
Boxplots of the estimates are plotted in Figure~\ref{figure1} and results of the RMSE are given in Table \ref{table1}. Code for reproducing the results is available at \url{https://git.io/Jt6bC}. 
% in the Appendix. 

Figure \ref{figure1} demonstrates that:
%\textit{Table \ref{table1}} and \textit{Table \ref{table2}} 
% one can observe that:
\begin{itemize}
 \item Both of the proposed estimators demonstrate an improvement over the naive estimator in terms of RMSE. For example, when $\tau^2=1$ and $\tau^2_{\bf{B}}=1/3,$ the Single  estimator $T_{c^*}$ improve the naive estimator by $17\%$ and when $\tau^2_{\bf{B}}=2/3$, the Selection  estimator $T_{\gamma}$  improves the naive by  $15\%$. When $\tau^2=2$ these improvements are even more substantial.
 
 \item As already been suggested in Example \ref{extreme_example}, the Selection estimator $T_{\gamma}$ works well when $\tau^2_{\textbf{B}}$ is large while the Single estimator $T_{\hat c^*}$ works well when $\tau^2_{\textbf{B}}$ is small.

\item  The PSI estimator is biased in a non-sparse setting. For example, when $\tau^2_{\bf{B}}=1/3$  the PSI has larger RMSE than the proposed estimators. When $\tau^2_{\bf{B}}=0.99$ the PSI has low bias  therefore and low RMSE. This is not surprising since the PSI estimator is based on the LASSO method which is known to work well when the true model that generates the data is sparse. 
    
    % This results in a  larger RMSE compared to the other estimators when $\tau^2_{\bf{B}}=0.25$ or $\tau^2_{\bf{B}}=0.75$, and a smaller RMSE when $\tau^2_{\bf{B}}=0.99$. This is not surprising since the PSI estimator is based on the LASSO method which
    % which is known to work well when the true model that generates the data is sparse. 
   
%   \item As $\tau^2$ increases, the  improvements of all estimators over the naive are more substantial.  For example, for $n=p=500,$ the relative improvement in RMSE   of the oracle estimator  over the naive increases from $25\%$ to $35\%$ as $\tau^2$ changes from $\tau^2=1$  to $\tau^2=2$. More generally, from other simulation results that are not shown here, we found that the larger the signal to noise ratio $\tau^2/\sigma^2$, the larger the improvement.   
   \end{itemize}
 
%\newpage
%%%%%%%%%%%%%%%%%%%%%
% Oracle plot
%%%%%%%%%%%%%%%%%%%%%%%%%%
% \begin{figure}[H]
%   \centering
%  \includegraphics[width=1\textwidth]{g_boxplot_comb_nnn.eps}
% \caption{
% Boxplots representing the estimators distribution . The x-axis stands for $\tau^2_{\bf{B}}$. The red dashed is the true value of $\tau^2.$ 
% % Boxplots for table \ref{table1}. The value of $\tau^2$ is shown by the red dashed line. The  naive, modified, proposed, oracle and the PSI estimators  are compared for sample sizes $n=100, 300, 500.$
% }
% \label{figure1}
% \end{figure}

% Oracle tables
% \begin{table}[H] %***
% \caption{Summary statistics.   The standard deviation of RMSE ($\sigma_{RMSE}$) was calculated using the delta method. 
%  The estimator with the lowest RMSE (excluding the oracle) is in  bold.}
%  \label{table1} \par
% \resizebox{\linewidth}{!}{
%  \renewcommand{\arraystretch}{0.7}   
% \scalebox{0.9}{
% \begin{tabular}{|ccccccccc|} \hline %***5truept
% $\tau^2_{\bf{B}}$ & $\tau^2$ &  $n$ &   Estimator & Mean & Bias & SE & RMSE & $1000\cdot\sigma_{RMSE}$   \\ \hline
% % \\[3pt] \hline
% 30\% & 1 & 500 & Naive & 0.99 & 0.01 & 0.209 & 0.209 & 8 \\ 
% 30\% & 1 & 500 & Selection\_O & 0.99 & 0.01 & 0.201 & 0.201 & 8 \\ 
% 30\% & 1 & 500 & Single\_O & 0.99 & 0.01 & 0.171 & 0.171 & 7 \\ 
% 30\% & 1 & 500 & OOE & 1 & 0 & 0.15 & 0.15 & 6 \\ 
% 30\% & 1 & 500 & PSI & 0.71 & 0.29 & 0.237 & 0.373 & 12 \\ 
% \hline 
% 70\% & 1 & 500 & Naive & 1.01 & -0.01 & 0.22 & 0.22 & 10 \\ 
% 70\% & 1 & 500 & Selection\_O & 1.01 & -0.01 & 0.182 & 0.182 & 8 \\ 
% 70\% & 1 & 500 & Single\_O & 1.01 & -0.01 & 0.215 & 0.215 & 10 \\ 
% 70\% & 1 & 500 & OOE & 1 & 0 & 0.16 & 0.16 & 7 \\ 
% 70\% & 1 & 500 & PSI & 0.83 & 0.17 & 0.166 & 0.24 & 7 \\ 
% \hline
% 99\% & 1 & 500 & Naive & 1.02 & -0.02 & 0.221 & 0.222 & 7 \\ 
% 99\% & 1 & 500 & Selection\_O & 1 & 0 & 0.157 & 0.157 & 6 \\ 
% 99\% & 1 & 500 & Single\_O & 1.02 & -0.02 & 0.221 & 0.222 & 7 \\ 
% 99\% & 1 & 500 & OOE & 1 & 0 & 0.149 & 0.149 & 4 \\ 
% 99\% & 1 & 500 & PSI & 0.99 & 0.01 & 0.13 & 0.13 & 4 \\ \hline
% \end{tabular}}
% }
% \end{table}

\begin{table}[H] %***
\caption{Summary statistics.   An estimate for the standard deviation of RMSE ($\hat\sigma_{RMSE}$) was calculated using the delta method. 
 The estimator with the lowest RMSE (excluding the oracle) is in  bold.}
 \label{table1} \par
\resizebox{\linewidth}{!}{
 \renewcommand{\arraystretch}{0.7}   
\scalebox{0.9}{
\begin{tabular}{|ccccccccc|} \hline %***5truept
$\tau^2_{\bf{B}}$ & $\tau^2$ &  $n$ &   Estimator & Mean & Bias & SE & RMSE & $1000\cdot\hat\sigma_{RMSE}$   \\ \hline
% \\[3pt] \hline
33\% & 1 & 400 & Naive & 1.02 & -0.02 & 0.258 & 0.258 & 19 \\ 
33\% & 1 & 400 & Selection & 1.02 & -0.02 & 0.245 & 0.244 & 18 \\ 
33\% & 1 & 400 & Single & 0.99 & 0.01 & 0.214 & \textbf{0.213} & 14 \\ 
33\% & 1 & 400 & OOE & 1.02 & -0.02 & 0.193 & 0.193 & 14 \\ 
33\% & 1 & 400 & PSI & 0.74 & 0.26 & 0.221 & 0.341 & 20 \\ 
\hline
66\% & 1 & 400 & Naive & 1.02 & -0.02 & 0.259 & 0.259 & 21 \\ 
66\% & 1 & 400 & Selection & 1.02 & -0.02 & 0.22 & \textbf{0.219} & 18 \\ 
66\% & 1 & 400 & Single & 1 & 0 & 0.234 & 0.233 & 18 \\ 
66\% & 1 & 400 & OOE & 1.02 & -0.02 & 0.185 & 0.185 & 15 \\ 
66\% & 1 & 400 & PSI & 0.84 & 0.16 & 0.172 & 0.231 & 13 \\
\hline
99\% & 1 & 400 & Naive & 1.02 & -0.02 & 0.261 & 0.261 & 28 \\ 
99\% & 1 & 400 & Selection & 1.01 & -0.01 & 0.172 & 0.171 & 13 \\ 
99\% & 1 & 400 & Single & 1.01 & -0.01 & 0.254 & 0.253 & 28 \\ 
99\% & 1 & 400 & OOE & 1.02 & -0.02 & 0.17 & 0.171 & 15 \\ 
99\% & 1 & 400 & PSI & 0.98 & 0.02 & 0.157 & \textbf{0.157} & 14 \\ \hline
\hline
33\% & 2 & 400 & Naive & 2.02 & -0.02 & 0.436 & 0.435 & 33 \\ 
33\% & 2 & 400 & Selection & 2.02 & -0.02 & 0.411 & 0.41 & 30 \\ 
33\% & 2 & 400 & Single & 1.96 & 0.04 &  0.342  & \textbf{0.342} & 22 \\ 
33\% & 2 & 400 & OOE & 2.02 & -0.02 & 0.286 & 0.286 & 21 \\ 
33\% & 2 & 400 & PSI & 1.65 & 0.35 & 0.395 & 0.529 & 32 \\ 
\hline
66\% & 2 & 400 & Naive & 2.03 & -0.03 & 0.443 & 0.441 & 38 \\ 
66\% & 2 & 400 & Selection & 2.02 & -0.02 & 0.362 & \textbf{0.36} & 30 \\ 
66\% & 2 & 400 & Single & 1.98 & 0.02 & 0.393 & 0.392 & 30 \\ 
66\% & 2 & 400 & OOE & 2.02 & -0.02 & 0.274 & 0.273 & 22 \\ 
66\% & 2 & 400 & PSI & 1.74 & 0.26 & 0.268 & 0.375 & 24 \\ 
\hline
99\% & 2 & 400 & Naive & 2.02 & -0.02 & 0.46 & 0.458 & 51 \\ 
99\% & 2 & 400 & Selection & 2 & 0 & 0.267 & 0.265 & 20 \\ 
99\% & 2 & 400 & Single & 2 & 0 & 0.446 & 0.443 & 50 \\ 
99\% & 2 & 400 & OOE & 2.02 & -0.02 & 0.251 & 0.25 & 22 \\ 
99\% & 2 & 400 & PSI & 1.97 & 0.03 & 0.243 & \textbf{0.243} & 20 \\ 
\hline
\end{tabular}}
}
\end{table}

% \begin{table}[H] \centering 
% \captionsetup{justification=centering}
%   \caption{\textit{\\Summary statistics equivalent to Table 1 with \(\tau _{\bf{B}}^2 = 0.99\)}
%      }
%   \label{table2} 
%     \renewcommand{\arraystretch}{0.7}   
%   \scalebox{0.9}{
% \begin{tabular}{@{\extracolsep{5pt}} ccccccccc} 
% \\[-1.8ex]\hline 
% \hline \\[-1.8ex] 
% $\tau^2$ &  n & p & Estimator & Mean & Bias & SE & RMSE & $1000\cdot\sigma_{RMSE}$  \\ 
% \hline\hline \\[-1.8ex] 
% 1 & 100 & 100 & Naive & 0.99 & 0.01 & 0.44 & 0.44 & 11  \\ 
% 1 & 100 & 100 & Modified & 0.99 & 0.01 & 0.42 & 0.42 & 10  \\ 
% 1 & 100 & 100 & Proposed & 0.84 & 0.16 & 0.39 & 0.42 & 9  \\ 
% 1 & 100 & 100 & Oracle & 1 & 0 & 0.34 & 0.34 & 8  \\ 
% 1 & 100 & 100 & \bf{PSI} & 1 & 0 & 0.29 & 0.29 & 7  \\ 
% \hline \\[-1.8ex] 
% 1 & 300 & 300 & Naive & 1 & 0 & 0.26 & 0.26 & 6  \\ 
% 1 & 300 & 300 & Modified & 1 & 0 & 0.23 & 0.23 & 5  \\ 
% 1 & 300 & 300 & Proposed & 0.96 & 0.04 & 0.22 & 0.23 & 5  \\ 
% 1 & 300 & 300 & Oracle & 1 & 0 & 0.2 & 0.2 & 5  \\ 
% 1 & 300 & 300 & \bf{PSI} & 0.99 & 0.01 & 0.16 & 0.16 & 4  \\ 
% \hline \\[-1.8ex] 
% 1 & 500 & 500 & Naive & 1 & 0 & 0.2 & 0.2 & 5  \\ 
% 1 & 500 & 500 & Modified & 1 & 0 & 0.17 & 0.17 & 5  \\ 
% 1 & 500 & 500 & Proposed & 0.99 & 0.01 & 0.16 & 0.16 & 4  \\ 
% 1 & 500 & 500 & Oracle & 1 & 0 & 0.15 & 0.15 & 4  \\ 
% 1 & 500 & 500 & \bf{PSI} & 0.99 & 0.01 & 0.11 & 0.11 & 3  \\ 
 
% \end{tabular}
% }
% \end{table} 
 \section{ Generalization to Other Estimators
 }\label{gener_es}
The suggested methodology in this paper is not limited  to improving only the naive estimator, but can also be generalized to other estimators. The key is to add zero-estimators that are highly correlated with our initial estimator of $\tau^2$; see Equation (\ref{variance_change}). Unlike the naive estimator, which is represented by a closed-form expression, other common estimators, such as the EigenPrism estimator \citep{janson2017eigenprism}, are computed numerically and do not have a closed-form representation. That makes the task of finding  optimal zero-estimators somewhat more challenging  since the zero-estimators'  coefficients   also need to be computed numerically. A comprehensive theory that generalizes the zero-estimate approach to other estimators, other than the naive,  is beyond the scope of this work. 
However, here we present a general algorithm that achieves improvement without claiming optimality. The algorithm is based on adding a single zero-estimator as in Section \ref{imp_single_zero_estimator}.
The algorithm below  approximates the  optimal-oracle coefficient $c^*$ given in \eqref{general_c} from bootstrap samples and then, returns a new estimator that is composed of both the initial estimator of $\tau^2$ and a single zero-estimator.   

\begin{algorithm}[H]\label{alg_emp}
\SetAlgoLined

\vspace{0.4 cm}
\textbf{Input:} 
 A dataset \(\left( {{\bf{X}},{\bf{Y}}} \right)\), an initial estimator $\tilde{\tau}^2,$ and a selection algorithm~$\gamma$.
\begin{enumerate}
 \item Calculate an initial estimator $\tilde{\tau}^2$ of $\tau^2$.
    
     \item \textbf{Bootstrap step:}  
      \begin{itemize}
        \item  Resample  with replacement $n$ observations from \(\left( {{{\bf{X}}},{{\bf{Y}}}} \right)\).
        \item Calculate the initial  estimator $\tilde{\tau}^2$ of $\tau^2.$ 
        \item Calculate the zero-estimator
       ${g_n} = \frac{1}{n}\sum\limits_{i = 1}^n {{g_i}}$ where ${g_i} = \sum\limits_{j < j'}^{} {{X_{ij}}{X_{ij'}}}$.
    \end{itemize}
    This procedure is repeated $B$ times  in order to produce
    $(\tilde{\tau}^2)^{*1},...,(\tilde{\tau}^2)^{*B}$
    and \(g_n^{*1},...,g_n^{*B}\). 
    \item   Approximate the  coefficient $c^*$ by  	\[\tilde{c}^* =  \frac{{\widehat {\cov\left( {\tilde{\tau}^2,{g_n}} \right)}}}{{\var\left( {{g_n}} \right)}}\,,\] where \(\widehat {{\cov} \left(  \cdot  \right)}\) denotes the empirical covariance from the bootstrap samples, and $\var(g_n)$ is known by the semi-supervised setting.
  \end{enumerate}
\KwResult{Return the empirical estimator \({T_{emp}} = {\tilde \tau ^2} - {\tilde c^*}{g_n}\)	.
}
  \caption{Empirical Estimator}
\end{algorithm}

We now  demonstrate the performance of the empirical estimator given by Algorithm \ref{alg_emp}  together with two initial estimators mentioned earlier: The  EigenPrism  \citep{janson2017eigenprism} and the  PSI which is described in \citet{taylor2018post} and was used in Section \ref{sim_res}.
We consider the same  setting as in Section \ref{sim_res}. 
Results are given in Tables \ref{table3}-\ref{table4} and the code for reproducing the results is available at \url{https://git.io/Jt6bC}. 

Tables \ref{table3}-\ref{table4} demonstrate that the standard error of the empirical estimators is equal to or lower than the standard error of the initial estimators, and  as $\tau^2$ increases, the  improvement  over the initial estimators is more substantial. As in Section \ref{sim_res}, the single zero-estimator approach works especially well when $\tau^2_B$ is small; otherwise, there is a small or no improvement, but also no additional variance or bias is introduced.
     This highlights the fact that the zero-estimator approach is not limited to improving only the naive estimator but rather has the potential to improve other estimators as well.
     
%  \begin{figure}[H]
%   \centering
%  \includegraphics[width=0.9\textwidth]{g_boxplot_Eigen.eps}
% \caption{
% Boxplots representing the estimators distribution . The x-axis stands for $\tau^2_{\bf{B}}$. The red dashed is the true value of $\tau^2.$ 
% % Boxplots for table \ref{table1}. The value of $\tau^2$ is shown by the red dashed line. The  naive, modified, proposed, oracle and the PSI estimators  are compared for sample sizes $n=100, 300, 500.$
% }
% \label{emp_eigen}
% \end{figure}

\begin{table}[H] %***
\caption{Summary statistics equivalent to Table 1 for the EigenPrism estimator.}
  \label{table3}  \par
\resizebox{\linewidth}{!}{
 \renewcommand{\arraystretch}{0.7}   
\scalebox{0.7}{
\begin{tabular}{|ccccccccc|} \hline %***5truept
$\tau^2_{\bf{B}}$ & $\tau^2$ &  $n$ &   Estimator & Mean & Bias & SE & RMSE & $1000\cdot\hat\sigma_{RMSE}$   \\ \hline
% \\[3pt] \hline
33\% & 1 & 400 & Eigenprism & 1.01 & -0.01 & 0.167 & 0.166 & 13 \\ 
33\% & 1 & 400 & Empirical Eigen & 1 & 0 & 0.154 & 0.154 & 13 \\ 
\hline
66\% & 1 & 400 & Eigenprism & 1.01 & -0.01 & 0.17 & 0.17 & 15 \\ 
66\% & 1 & 400 & Empirical Eigen & 1.01 & -0.01 & 0.164 & 0.163 & 15 \\
\hline 
99\% & 1 & 400 & Eigenprism & 1.01 & -0.01 & 0.175 & 0.174 & 15 \\ 
99\% & 1 & 400 & Empirical Eigen & 1 & 0 & 0.175 & 0.174 & 15 \\ 
\hline\hline 
33\% & 2 & 400 & Eigenprism & 2.01 & -0.01 & 0.245 & 0.243 & 18 \\ 
33\% & 2 & 400 & Empirical Eigen & 2 & 0 & 0.208 & 0.207 & 19 \\ 
\hline  
66\% & 2 & 400 & Eigenprism & 2 & 0 & 0.247 & 0.246 & 23 \\ 
66\% & 2 & 400 & Empirical Eigen & 2 & 0 & 0.231 & 0.23 & 24 \\ 
\hline 
99\% & 2 & 400 & Eigenprism & 2 & 0 & 0.259 & 0.257 & 23 \\ 
99\% & 2 & 400 & Empirical Eigen & 2 & 0 & 0.259 & 0.257 & 23 \\
\hline  
\end{tabular}}
} 
\end{table}

% \begin{figure}[H]
%   \centering
%  \includegraphics[width=0.9\textwidth]{g_boxplot_PSI.eps}
% \caption{
% Boxplots representing the estimators distribution . The x-axis stands for $\tau^2_{\bf{B}}$. The red dashed is the true value of $\tau^2.$ 
% % Boxplots for table \ref{table1}. The value of $\tau^2$ is shown by the red dashed line. The  naive, modified, proposed, oracle and the PSI estimators  are compared for sample sizes $n=100, 300, 500.$
% }
% \label{emp_PSI}
% \end{figure}

\begin{table}[H] %***
\caption{Summary statistics equivalent to Table 1 for the PSI estimator.}
  \label{table4} \par
\resizebox{\linewidth}{!}{
\renewcommand{\arraystretch}{0.7}   
\scalebox{0.9}{
\begin{tabular}{|ccccccccc|} \hline %***5truept
$\tau^2_{\bf{B}}$ & $\tau^2$ &  $n$ &   Estimator & Mean & Bias & SE & RMSE & $1000\cdot\hat\sigma_{RMSE}$   \\ \hline
% \\[3pt] \hline
33\% & 1 & 400 & PSI & 0.74 & 0.26 & 0.221 & 0.341 & 20 \\ 
33\% & 1 & 400 & Empirical PSI & 0.73 & 0.27 & 0.21 & 0.339 & 19 \\ 
\hline 
66\% & 1 & 400 & PSI & 0.84 & 0.16 & 0.172 & 0.231 & 13 \\ 
66\% & 1 & 400 & Empirical PSI & 0.84 & 0.16 & 0.163 & 0.227 & 12 \\ 
\hline 
99\% & 1 & 400 & PSI & 0.98 & 0.02 & 0.157 & 0.157 & 14 \\ 
99\% & 1 & 400 & Empirical PSI & 0.98 & 0.02 & 0.155 & 0.155 & 13 \\ 
\hline\hline
33\% & 2 & 400 & PSI & 1.65 & 0.35 & 0.395 & 0.529 & 32 \\ 
33\% & 2 & 400 & Empirical PSI & 1.63 & 0.37 & 0.355 & 0.51 & 31 \\ 
\hline 
66\% & 2 & 400 & PSI & 1.74 & 0.26 & 0.268 & 0.375 & 24 \\ 
66\% & 2 & 400 & Empirical PSI & 1.73 & 0.27 & 0.251 & 0.371 & 24 \\ 
\hline 
99\% & 2 & 400 & PSI & 1.97 & 0.03 & 0.243 & 0.243 & 20 \\ 
99\% & 2 & 400 & Empirical PSI & 1.97 & 0.03 & 0.237 & 0.238 & 19 \\ 
\hline  
\end{tabular}}
}
\end{table}

\section{Discussion }\label{discuss}
This paper  presents a new approach for  improving estimation of the explained variance $\tau^2$ of a high-dimensional regression model in a semi-supervised setting without assuming sparsity. The key idea is to use zero-estimator that  is  correlated with the initial unbiased estimator of $\tau^2$ in order to lower its variance without introducing additional bias. The semi-supervised setting, where the number of observations is much greater than  the number of responses, allows us to construct such zero-estimators. We introduced a new notion of optimality with respect to  zero-estimators and presented an oracle-estimator that achieves this type of optimality. 
We proposed two different (non-oracle) estimators that showed 
 a significant reduction, but not optimal, in the asymptotic variance of the  naive estimator.
Our simulations showed that our approach can be generalized to other types of initial estimators other than the naive estimator.

Many open questions remain for future research.
While our proposed estimators improved the naive estimator, it did not achieve the optimal improvement of the oracle estimator. 
%The optimal oracle-estimator remained an unattainable gold standard in term of optimality that can be achieved by using the zero-estimator approach. 
Thus, it remains  unclear if and how one can  achieve optimal improvement. 
Moreover, in this work,  strong assumption was made about the unsupervised data size, i.e., $N = \infty$.  
Thus, generalizing the suggested approach by relaxing  this assumption to allow for a more general setting with \textit{finite} $N \gg n$ is a natural direction for future work.
A more ambitious future goal would be to extend the suggested approach to generalized linear models (GLM), and specifically to logistic regression. In this case, the concepts of signal and noise levels are less clear and are more challenging to define.

%\subsection{Illustration }\label{sub_sec:Illustration}

%%%%%%%%%%%%%%%%%%%%%%%%%%%%%%%%%%%%%%%%%%%%%%%%%%%%%%%%%%%%%%%%%%%%%%%%%%%%%%%%%%%%%%%%%%%%%%%%%%%%%%%%%%%%%%%%%%%%%%%%%%%%
% \section*{Supplementary Materials}

% Contain
% the brief description of the online supplementary materials.
% \par
%%%%%%%%%%%%%%%%%%%%%%%%%%%%%%%%%%%%%%%%%%%%%%%%%%%%%%%%%%%%%%%%%%%%%%%%%%%%%%%%%%%%%%%%%%%%%%%%%%%%%%%%%%%%%%%%%%%%%%%%%%%%
%\section*{Acknowledgements}

%Write the acknowledgements here.
%\par

%%%%%%%%%%%%%%%%%%%%%%%%%%%%%%%%%%%%%%%%%%%%%%%%%%%%%%%%%%%%%%%%%%%%%%%%%%%%%%%%%%%%%%%%%%

\bibhang=1.7pc
\bibsep=2pt
\fontsize{9}{14pt plus.8pt minus .6pt}\selectfont
\renewcommand\bibname{\large \bf References}
%\begin{thebibliography}{11}
\expandafter\ifx\csname
natexlab\endcsname\relax\def\natexlab#1{#1}\fi
\expandafter\ifx\csname url\endcsname\relax
  \def\url#1{\texttt{#1}}\fi
\expandafter\ifx\csname urlprefix\endcsname\relax\def\urlprefix{URL}\fi

%% use bibfile 
  \bibliographystyle{chicago}      % Chicago style, author-year citations
  \bibliography{bib_ref}   % name your BibTeX data base

\begin{thebibliography}{}

\bibitem[\protect\citeauthoryear{Cai and Guo}{Cai and
  Guo}{2020}]{tony2020semisupervised}
Cai, T. and Z.~Guo (2020).
\newblock Semisupervised inference for explained variance in high dimensional
  linear regression and its applications.
\newblock {\em Journal of the Royal Statistical Society: Series B (Statistical
  Methodology)\/}.

\bibitem[\protect\citeauthoryear{Candes, Fan, Janson, and Lv}{Candes
  et~al.}{2017}]{candes2017panning}
Candes, E., Y.~Fan, L.~Janson, and J.~Lv (2017).
\newblock Panning for gold: Model-x knockoffs for high-dimensional controlled
  variable selection.

\bibitem[\protect\citeauthoryear{Dicker}{Dicker}{2014}]{Dicker}
Dicker, L.~H. (2014).
\newblock {Variance estimation in high-dimensional linear models}.
\newblock {\em Biometrika\/}~{\em 101\/}(2), 269--284.

\bibitem[\protect\citeauthoryear{Eichler, Flint, Gibson, Kong, Leal, Moore, and
  Nadeau}{Eichler et~al.}{2010}]{eichler2010missing}
Eichler, E.~E., J.~Flint, G.~Gibson, A.~Kong, S.~M. Leal, J.~H. Moore, and
  J.~H. Nadeau (2010).
\newblock Missing heritability and strategies for finding the underlying causes
  of complex disease.
\newblock {\em Nature Reviews Genetics\/}~{\em 11\/}(6), 446--450.

\bibitem[\protect\citeauthoryear{Fan, Guo, and Hao}{Fan
  et~al.}{2012}]{fan2012variance}
Fan, J., S.~Guo, and N.~Hao (2012).
\newblock Variance estimation using refitted cross-validation in ultrahigh
  dimensional regression.
\newblock {\em Journal of the Royal Statistical Society: Series B (Statistical
  Methodology)\/}~{\em 74\/}(1), 37--65.

\bibitem[\protect\citeauthoryear{Gorfine, Berndt, Chang-Claude, Hoffmeister,
  Le~Marchand, Potter, Slattery, Keret, Peters, and Hsu}{Gorfine
  et~al.}{2017}]{gorfine2017heritability}
Gorfine, M., S.~I. Berndt, J.~Chang-Claude, M.~Hoffmeister, L.~Le~Marchand,
  J.~Potter, M.~L. Slattery, N.~Keret, U.~Peters, and L.~Hsu (2017).
\newblock Heritability estimation using a regularized regression approach
  {(HERRA)}: Applicable to continuous, dichotomous or age-at-onset outcome.
\newblock {\em PloS one\/}~{\em 12\/}(8), e0181269.

\bibitem[\protect\citeauthoryear{Hoeffding}{Hoeffding}{1977}]{hoeffding1977some}
Hoeffding, W. (1977).
\newblock Some incomplete and boundedly complete families of distributions.
\newblock {\em The Annals of Statistics\/}, 278--291.

\bibitem[\protect\citeauthoryear{Janson, Barber, and Candes}{Janson
  et~al.}{2017}]{janson2017eigenprism}
Janson, L., R.~F. Barber, and E.~Candes (2017).
\newblock Eigenprism: inference for high dimensional signal-to-noise ratios.
\newblock {\em Journal of the Royal Statistical Society: Series B (Statistical
  Methodology)\/}~{\em 79\/}(4), 1037--1065.

\bibitem[\protect\citeauthoryear{Kong and Valiant}{Kong and
  Valiant}{2018}]{kong2018estimating}
Kong, W. and G.~Valiant (2018).
\newblock Estimating learnability in the sublinear data regime.
\newblock {\em Advances in Neural Information Processing Systems\/}~{\em 31},
  5455--5464.

\bibitem[\protect\citeauthoryear{Oda, Yanagihara, et~al.}{Oda
  et~al.}{2020}]{oda2020fast}
Oda, R., H.~Yanagihara, et~al. (2020).
\newblock A fast and consistent variable selection method for high-dimensional
  multivariate linear regression with a large number of explanatory variables.
\newblock {\em Electronic Journal of Statistics\/}~{\em 14\/}(1), 1386--1412.

\bibitem[\protect\citeauthoryear{Schwartzman, Schork, Zablocki, Thompson,
  et~al.}{Schwartzman et~al.}{2019}]{schwartzman2019simple}
Schwartzman, A., A.~J. Schork, R.~Zablocki, W.~K. Thompson, et~al. (2019).
\newblock A simple, consistent estimator of {SNP} heritability from genome-wide
  association studies.
\newblock {\em The Annals of Applied Statistics\/}~{\em 13\/}(4), 2509--2538.

\bibitem[\protect\citeauthoryear{Sun and Zhang}{Sun and
  Zhang}{2012}]{sun2012scaled}
Sun, T. and C.-H. Zhang (2012).
\newblock Scaled sparse linear regression.
\newblock {\em Biometrika\/}~{\em 99\/}(4), 879--898.

\bibitem[\protect\citeauthoryear{Taylor and Tibshirani}{Taylor and
  Tibshirani}{2018}]{taylor2018post}
Taylor, J. and R.~Tibshirani (2018).
\newblock Post-selection inference for-penalized likelihood models.
\newblock {\em Canadian Journal of Statistics\/}~{\em 46\/}(1), 41--61.

\bibitem[\protect\citeauthoryear{van~der Vaart}{van~der
  Vaart}{2000}]{van2000asymptotic}
van~der Vaart, A.~W. (2000).
\newblock {\em Asymptotic statistics}.
\newblock Cambridge university press.

\bibitem[\protect\citeauthoryear{Zambom and Kim}{Zambom and
  Kim}{2018}]{zambom2018consistent}
Zambom, A.~Z. and J.~Kim (2018).
\newblock Consistent significance controlled variable selection in
  high-dimensional regression.
\newblock {\em Stat\/}~{\em 7\/}(1), e210.

\end{thebibliography}

%%  Another method
% \begin{thebibliography}{}

% \bibitem[Curtis(1943)]{1943}
% Curtis, M. (1943).
% {\em Documents on International Affairs, 1938}, Volume~II.
%  Oxford University Press, London.

% \bibitem[Eubank(2004)]{Eubank}
% Eubank, K. (2004).
%  {\em The origins of World War II}. 
%  3rd Edition.
% Harlan Davidson, Wheeling, Ill.

% \bibitem[Gellately(1988)]{1988}
% Gellately, R. (1988).
%  The gestapo and german society: Political denunciation in the gestapo
%   case files.
% {\em The Journal of Modern History}~{\bf 60}, 654--694.

% \bibitem[Noakes and Pridham(2001)]{Noakes}
% Noakes, J. and G.~Pridham (2001).
%  {\em Nazism, 1919-1945. Vol. 3: Foreign Policy, War and Racial
%   Extermination}.
% University of Exeter Press,  Exeter.

% \end{thebibliography}

%%%%%%%%%%%%%%%%%%%%%%%%%%%%%%%%%%%%%%%%%%%%%
%   Appendix - Ilan
%%%%%%%%%%%%%%%%%%%%%%%%%%%%%%%%%%%%%%%%%%%%%

\newpage

\section{Appendix}\label{appnd}

\noindent\textbf{\textit{Proof of Lemma \ref{asymptotic_normality_Naive}}}:\\
Notice that ${{\textbf{\emph{X}}}^T}{\bf{Y}} = {\left( {\sum\limits_{i = 1}^n {{W_{i1}}} ,...,\sum\limits_{i = 1}^n {{W_{ip}}} } \right)^T}$ where  $\textbf{X}$ is the $n\times p$ design matrix and  $\textbf{Y}=(Y_1,...,Y_n)^T$.
Thus, the naive estimator can be also written as
$$\hat\tau^2  = \frac{1}{{n\left( {n - 1} \right)}}\sum\limits_{{i_1} \ne {i_2}}^{} {\sum\limits_{j = 1}^p {W_{{i_1}j}^{}{W_{{i_2}j}}  = } } \frac{{{{\left\| {{{\textbf{\emph{X}}}^T}{\bf{Y}}} \right\|}^2} - \sum\limits_{j = 1}^p {\sum\limits_{i = 1}^n {W_{ij}^2} } }}{{n\left( {n - 1} \right)}}.$$
The Dicker estimate for $\tau^2$ ia given by
$\hat\tau^2_{Dicker}\equiv \frac{{{{\left\| {{{\textbf{\emph{X}}}^T}{\bf{Y}}} \right\|}^2} - p{{\left\| {\bf{Y}} \right\|}^2}}}{{n\left( {n + 1} \right)}}.$ 
We need to prove that root-$n$  times the difference between the estimators converges in probability to zero, i.e., $\sqrt{n}\left(\hat\tau^2_{Dicker}-\hat\tau^2    \right)\overset{p}{\rightarrow}0$.
We have,

\begin{equation}\label{eq:U_Dicker_U}
\begin{split}
   \sqrt{n}\left(\hat\tau^2_{Dicker}-\hat\tau^2    \right) &= \sqrt n \left( {\frac{{{{\left\| {{{\bf{X}}^T}{\bf{Y}}} \right\|}^2} - p{{\left\| {\bf{Y}} \right\|}^2}}}{{n\left( {n + 1} \right)}} - \frac{{{{\left\| {{{\bf{X}}^T}{\bf{Y}}} \right\|}^2} - \sum\limits_{j = 1}^p {\sum\limits_{i = 1}^n {W_{ij}^2} } }}{{n\left( {n - 1} \right)}}} \right)\\
    &=  \sqrt n \left( {\frac{{\sum\limits_{j = 1}^p {\sum\limits_{i = 1}^n {W_{ij}^2} } }}{{n\left( {n - 1} \right)}} - \frac{{p{{\left\| {\bf{Y}} \right\|}^2}}}{{n\left( {n + 1} \right)}} - \frac{{2{{\left\| {{{\bf{X}}^T}{\bf{Y}}} \right\|}^2}}}{{n\left( {n - 1} \right)\left( {n + 1} \right)}}} \right).
\end{split}
\end{equation}
It is enough to prove that:
\begin{enumerate}

    \item $n^{-1.5} \left(\sum\limits_{j = 1}^p \sum\limits_{i = 1}^n {W_{ij}^2}  - p{{\left\| {\bf{Y}} \right\|}^2} \right) \overset{p}{\rightarrow}0,$ 
    \item 
$ n^{-2.5} \left( {{{{\left\| {{{\bf{X}}^T}{\bf{Y}}} \right\|}^2}}} \right)\overset{p}{\rightarrow}0.$\\
\end{enumerate}
We start with the first term,
    \begin{multline} \label{difference}
      n^{-1.5} \left(\sum\limits_{j = 1}^p \sum\limits_{i = 1}^n {W_{ij}^2}  - p{{\left\| {\bf{Y}} \right\|}^2} \right) 
  = n^{-1.5}\left( {\sum\limits_{j = 1}^p {\sum\limits_{i = 1}^n {Y_i^2X_{ij}^2} }  - p\sum\limits_{i = 1}^n {Y_i^2} } \right) \\
 = n^{-1.5}\left( {\sum\limits_{i = 1}^n {Y_i^2\sum\limits_{j = 1}^p {X_{ij}^2} }  - p\sum\limits_{i = 1}^n {Y_i^2} } \right)
 = n^{-0.5}\sum\limits_{i = 1}^n {Y_i^2} \left[ {\frac{1}{n}\sum\limits_{j = 1}^p {\left( {X_{ij}^2 - 1} \right)} } \right]
 \equiv n^{-0.5}\sum\limits_{i = 1}^n {{\omega _i}}
 \end{multline}
where ${\omega _i} = Y_i^2\left[ {\frac{1}{n}\sum\limits_j^{} {\left\{ {X_{ij}^2 - 1} \right\}} } \right].$ Notice that $\omega_i$   depends on $n$  but this is suppressed in the notation. In order to show that $n^{-0.5}\sum\limits_{i = 1}^n {{\omega _i}} \overset{p}{\rightarrow}0,$ it is enough to show that $E\left( {n^{-0.5}\sum\limits_{i = 1}^n {{\omega _i}} } \right) \to 0$ and $\var\left( {n^{-0.5}\sum\limits_{i = 1}^n {{\omega _i}} } \right) \to 0.$
Moreover, since $E\left( {n^{-0.5}\sum\limits_{i = 1}^n {{\omega _i}} } \right) = \sqrt n E\left( {{\omega _i}} \right)$ and $\var\left( {n^{-0.5}\sum\limits_{i = 1}^n {{\omega _i}} } \right) = \var\left( {{\omega _i}} \right) = E\left( {\omega _i^2} \right) - {\left[ {E\left( {{\omega _i}} \right)} \right]^2},$ it is enough to show that  $\sqrt n E\left( {{\omega _i}} \right)$ and  $E\left( {\omega _i^2} \right)$ converge to zero.\\
  Consider now \(\sqrt n E\left( {{\omega _i}} \right)\). 
  By (\ref{difference}) we have 
  $$\sum\limits_{i = 1}^n {{\omega _i}} = \frac{1}{n}\left[ {\sum\limits_{j = 1}^p {\sum\limits_{i = 1}^n {W_{ij}^2} }  - p{{\left\| {\bf{Y}} \right\|}^2}} \right] .$$
   \\ Taking expectation of both sides,
  \[ \sum\limits_{i = 1}^n {E\left( {{\omega _i}} \right)} = \frac{1}{n}\left[ {\sum\limits_{j = 1}^p {\sum\limits_{i = 1}^n {E\left( {W_{ij}^2} \right)} }  - pE\left( {{{\left\| {\bf{Y}} \right\|}^2}} \right)} \right]  .\]
  Now,  notice that 
  \begin{equation}\label{mean_W2}
E(W_{ij}^2) =E[X_{ij}^2(\beta^T X+\epsilon)^2]= {\left\| \beta  \right\|^2} + {\sigma ^2}  + \beta _j^2[E (X_{ij}^4) - 1] ={\tau ^2}+ {\sigma ^2} +  2\beta _j^2.      
  \end{equation}
 Also notice that $Y_i^2/\left( {\sigma _\varepsilon ^2 + {\tau ^2}} \right)\sim \chi_1^2,$ and hence      $E\left( {{{\left\| {\bf{Y}} \right\|}^2}} \right) = n\left( {{\tau ^2} + {\sigma ^2}} \right)$. Therefore,
 \begin{align*}
 nE\left( {{\omega _i}} \right) &= \frac{1}{n}\left[ {\sum\limits_{j = 1}^p {\sum\limits_{i = 1}^n {\left( {{\tau ^2} + {\sigma ^2} + 2\beta _j^2} \right)} }  - pn\left( {{\tau ^2} + {\sigma ^2}} \right)} \right]\\ &= \frac{1}{n}\left[ {\sum\limits_{i = 1}^n {\left[ {p\left( {{\tau ^2} + {\sigma ^2}} \right) + 2{\tau ^2}} \right] - pn\left( {{\tau ^2} + {\sigma ^2}} \right)} } \right] = 2{\tau ^2}    
 \end{align*}
    which implies that \(\sqrt n E\left( {{\omega _i}} \right) = \frac{{2{\tau ^2}}}{{\sqrt n }}\overset{}{\rightarrow}~0.\)
    
Consider now $E(\omega_i^2)$. 
By Cauchy-Schwartz,
\begin{equation*}
 E\left( {\omega _i^2} \right) = E\left( {Y_i^4{{\left[ {n^{-1}\sum\limits_{j=1}^{p} {\left\{ {X_{ij}^2 - 1} \right\}} } \right]}^2}} \right) \le {\left\{ {E\left( {Y_i^8} \right)} \right\}^{1/2}}{\left\{ {E\left( {{{\left[ {n^{-1}\sum\limits_{j=1}^{p} {\left\{ {X_{ij}^2 - 1} \right\}} } \right]}^4}} \right)} \right\}^{1/2}}.
\end{equation*}
Notice that  ${Y_i}\sim N\left( {0,{\tau ^2} + {\sigma ^2}} \right)$   by construction and therefore $E(Y_i^8)=O(1)$ as $n$ and $p$ go to infinity.
Let $V_j=X_{ij}^2-1$ and notice that $E(V_j)=0.$ We have $$	 
E\left( {{{\left[ {{n^{ - 1}}\sum\limits_{j = 1}^p {\left\{ {X_{ij}^2 - 1} \right\}} } \right]}^4}} \right) = E\left( {{{\left[ {{n^{ - 1}}\sum\limits_{j = 1}^p {{V_j}} } \right]}^4}} \right) = {n^{ - 4}}\sum\limits_{{j_1},{j_2},{j_3},{j_4}}^{} {E\left( {{V_{{j_1}}}{V_{{j_2}}}{V_{{j_3}}}{V_{{j_4}}}} \right)} .$$
The expectation 
$\sum\limits_{{j_1},{j_2},{j_3},{j_4}}^{} {E\left( {{V_{{j_1}}}{V_{{j_2}}}{V_{{j_3}}}{V_{{j_4}}}} \right)} $
  is not 0 when $j_1=j_2$  and $j_3=j_4$  (up to permutations) or when all terms are equal. In the first case we have 
  $$\sum\limits_{j \ne j'}^{} {E\left( {V_j^2V_{j'}^2} \right) = } \sum\limits_{j \ne j'}^{} {{{\left[ {E\left( {V_j^2} \right)} \right]}^2} = } p\left( {p - 1} \right){\left[ {E\left\{ {{{\left( {X_{ij}^2 - 1} \right)}^2}} \right\}} \right]^2} \le C_1p^2,$$ for a positive constant $C_1.$
  In the second case we have
  $\sum\limits_{j=1}^{p} {E\left( {V_j^4} \right) = } pE\left[ {{{\left( {X_{ij}^2 - 1} \right)}^4}} \right] \le C_2p,$ for a positive constant $C_2$.
  Hence, as $p$  and $n$  have the same order of  magnitude, we have
  $${\left\{ {E\left[ {{{\left( {{n^{ - 1}}\sum\limits_{j = 1}^p {\left\{ {X_{ij}^2 - 1} \right\}} } \right)}^4}} \right]} \right\}^{1/2}} = {\left\{ {{n^{ - 4}}\sum\limits_{{j_1},{j_2},{j_3},{j_4}}^{} {E\left( {{V_{{j_1}}}{V_{{j_2}}}{V_{{j_3}}}{V_{{j_4}}}} \right)} } \right\}^{1/2}}\le {\left\{ {{n^{ - 4}} \cdot O\left( {{p^2}} \right)} \right\}^{1/2}} \le K/n,$$ which implies  $E(\omega_i^2)\le K_1/n\rightarrow 0,$ where   $K$ and $K_1$ are  positive constants.  
This completes the proof that 
$$n^{-1.5} \left(\sum\limits_{j = 1}^p \sum\limits_{i = 1}^n {W_{ij}^2}  - p{{\left\| {\bf{Y}} \right\|}^2} \right) \overset{p}{\rightarrow}~0.$$  

We now move to prove that
$ n^{-2.5} \left( {{{{\left\| {{{\bf{X}}^T}{\bf{Y}}} \right\|}^2}}} \right)\overset{p}{\rightarrow}0.$
By Markov's inequality, for $\epsilon>0$
\[P\left( {{n^{ - 2.5}}{{\left\| {{{\bf{X}}^T}{\bf{Y}}} \right\|}^2} > \varepsilon } \right) \le {n^{ - 2.5}}E\left( {{{\left\| {{{\bf{X}}^T}{\bf{Y}}} \right\|}^2}} \right)/\varepsilon. \]
Thus, it is enough to show that \({n^{ - 2}}E\left( {{{\left\| {{{\bf{X}}^T}{\bf{Y}}} \right\|}^2}} \right)\) is bounded. 
Notice that
\begin{equation*}
    \begin{split}
E\left( {{{\left\| {{{\bf{X}}^T}{\bf{Y}}} \right\|}^2}} \right) &= \sum\limits_{{i_1},{i_2}}^{} {\sum\limits_{j = 1}^p {E\left( {{W_{{i_1}j}}{W_{{i_2}j}}} \right)} }  = \sum\limits_{i = 1}^n {\sum\limits_{j = 1}^p {E\left( {W_{ij}^2} \right)} }  + \sum\limits_{{i_1} \ne {i_2}}^{} {\sum\limits_{j = 1}^p {E\left( {{W_{{i_1}j}}{W_{{i_2}j}}} \right)} } \\ 
&= \sum\limits_{i = 1}^n {\sum\limits_{j = 1}^p {\left( {{\tau ^2} + {\sigma ^2} + 2\beta _j^2} \right)} }  + \sum\limits_{{i_1} \ne {i_2}}^{} {\sum\limits_{j = 1}^p {\beta _j^2} }  = n\left[ {p\left( {{\tau ^2} + {\sigma ^2}} \right) + 2{\tau ^2}} \right] + n\left( {n - 1} \right){\tau ^2}\\
&= n\left[ {p\left( {{\tau ^2} + {\sigma ^2}} \right) + \left( {n + 1} \right){\tau ^2}} \right],
    \end{split}
\end{equation*}
where we used \eqref{mean_W2} in the third equality.
Therefore, 
\({n^{ - 2}}E\left( {{{\left\| {{{\bf{X}}^T}{\bf{Y}}} \right\|}^2}} \right) = {n^{ - 1}}\left[ {p\left( {{\tau ^2} + {\sigma ^2}} \right) + \left( {n + 1} \right){\tau ^2}} \right]\). 
Since   $p$  and $n$  have the same order of  magnitude and $\tau^2$+$\sigma^2$ is bounded by assumption, then
\({n^{ - 2}}E\left( {{{\left\| {{{\bf{X}}^T}{\bf{Y}}} \right\|}^2}} \right)\) is also bounded. This completes the proof of 
$ n^{-2.5} \left( {{{{\left\| {{{\bf{X}}^T}{\bf{Y}}} \right\|}^2}}} \right)\overset{p}{\rightarrow}0$
and hence \( \sqrt n \left( {\hat \tau _{Dicker}^2 - {{\hat \tau }^2}} \right)\overset{p}{\rightarrow}~0.\)
\qed\\

\noindent\textbf{\textit{Proof of Corollary \ref{normality_naive}}}:\\
According to Corollary 1 in \cite{Dicker}, we have
$$\frac{{\sqrt n \left( {{\hat \tau_{Dicker} } - {\tau ^2}} \right)}}{\psi } \overset{D}{\rightarrow} N\left( {0,1} \right),$$
where 
\(\psi  = 2\left\{ {\left( {1 + \frac{p}{n}} \right){{\left( {\sigma ^2 + {\tau ^2}} \right)}^2} - \sigma^4 + 3{\tau ^4}} \right\},\) given that $p/n$ converges to a constant. 
Therefore we can write
\[\frac{{\sqrt n \left( {{{\hat \tau }^2} - {\tau ^2}} \right)}}{\psi } = \frac{1}{\psi }\left[ {\sqrt n \left( {{{\hat \tau }^2} - {{\hat \tau }_{Dicker }}} \right) + \sqrt n \left( {  {{\hat \tau }_{Dicker }} - {\tau ^2}} \right)} \right],\]
and obtain $
\sqrt n \left( {\frac{{{{\hat \tau }^2} - {\tau ^2}}}{\psi }} \right)\overset{D}{\rightarrow} N(0,1)\,
$ by Slutsky's theorem. \qed\\

\noindent\textbf {\textit{Proof of Proposition~\ref{var_naive1}}}:\\ 
 Let ${{\bf{W}}_i} = {\left( {{W_{i1}},...,{W_{ip}}} \right)^T}$ and notice that  $\hat\tau^2= \frac{1}{{n\left( {n - 1} \right)}} \sum\limits_{{i_1} \ne {i_2}}^n\sum\limits_{j = 1}^p {{W_{{i_1}j}}{W_{{i_2}j}}} $ is a U-statistic of order~2 with the kernel \(h\left( {{{\bf{w}}_1},{{\bf{w}}_2}} \right) = {\bf{w}}_1^T{{\bf{w}}_2} = \sum\limits_{j = 1}^p {{w_{1j}}{w_{2j}}} \), where \({{\bf{w}}_i} \in {\mathbb{R}^p}\).
 
 By Theorem 12.3 in \cite{van2000asymptotic},
 \begin{equation}\label{naive_var_form}
  {\var} \left( {{{\hat \tau }^2}} \right) = \frac{{4\left( {n - 2} \right)}}{{n\left( {n - 1} \right)}}{\zeta _1} + \frac{2}{{n\left( {n - 1} \right)}}{\zeta _2},   
  \end{equation}
  where \({\zeta _1} = {\cov} \left[ {h\left( {{{\bf{W}}_1},{{\bf{W}}_2}} \right),h\left( {{{\bf{W}}_1},{{\widetilde{\bf{ {W} }}_2}}} \right)} \right]\)
  and 
  ${\zeta _2} = {\cov} \left[ {h\left( {{{\bf{W}}_1},{{\bf{W}}_2}} \right),h\left( {{{\bf{W}}_1},{{\bf{W}}_2}} \right)} \right]$
  where $\widetilde{\bf{ {W} }}_2$ is an independent copy of $\bf{W}_2$.
  %By Chapter 12 in \cite{van2000asymptotic}, the variance of a U-statistic of order $r$ is given by
  %$$\var\left( U \right) = \sum\limits_{c = 1}^r {\frac{{r{!^2}}}{{c!\left( {r - c} \right){!^2}}}\frac{{\left( {n - r} \right)\left( {n - r - 1} \right)...\left( {n - 2r + c + 1} \right)}}{{n\left( {n - 1} \right)...\left( {n - r + 1} \right)}}} {\zeta _c}$$
  %where  $\zeta_c$ is the covariance of the kernels if $c$ variables are  in common. In our case $r=2,$ and thus  
  %\({\zeta _1} = {\cov} \left[ {h\left( {{{\bf{W}}_1},{{\bf{W}}_2}} \right),h\left( {{{\bf{W}}_1},{{\widetilde{\bf{ {W} }}_2}}} \right)} \right]\)
  %and 
  %${\zeta _2} = {\cov} \left[ {h\left( {{{\bf{W}}_1},{{\bf{W}}_2}} \right),h\left( {{{\bf{W}}_1},{{\bf{W}}_2}} \right)} \right]$
  %where $\widetilde{\bf{ {W} }}_2$ is an independent copy of $\bf{W}_2$.
  %Therefore, the variance of the naive estimator has the following form
  %\begin{equation}\label{naive_var_form}
  %{\var} \left( {{{\hat \tau }^2}} \right) = \frac{{4\left( {n - 2} \right)}}{{n\left( {n - 1} %\right)}}{\zeta _1} + \frac{2}{{n\left( {n - 1} \right)}}{\zeta _2.}    
  %\end{equation}
  Now, let  \({\bf{A}} = E\left( {{{\bf{W}}_i}{\bf{W}}_i^T} \right)\) be a $p \times p$ matrix and   notice that
  \begin{align*}
  {\zeta _1} &= \cov\left[ {h\left( {{{\bf{W}}_1},{{\bf{W}}_2}} \right),h\left( {{{\bf{W}}_1},{{\widetilde {\bf{W}}}_2}} \right)} \right]
  \\ &= \sum\limits_{j,j'}^p {\cov\left( {{W_{1j}}{W_{2j}},{W_{1j'}}{{\widetilde W}_{2j'}}} \right)}  = \sum\limits_{j,j'}^p {\left( {{\beta _j}{\beta _{j'}}E\left[ {{W_{1j}}{W_{1j'}}} \right] - \beta _j^2\beta _{j'}^2} \right)} \\
  &= {\beta ^T}{\bf{A}}\beta  - {\left\| \beta  \right\|^4}    
  \end{align*}
and 
\begin{align*}
{\zeta _2} &= {\cov} \left[ {h\left( {{{\bf{W}}_1},{{\bf{W}}_2}} \right),h\left( {{{\bf{W}}_1},{{\bf{W}}_2}} \right)} \right]\\
&= \sum\limits_{j,j'}^{} {{\cov} \left( {{W_{1j}}{W_{2j}},{W_{1j'}}{W_{2j'}}} \right)}
=\sum\limits_{j,j'}^{}\left( {{{\left( {E\left[ {{W_{1j}}{W_{1j'}}} \right]} \right)}^2} - \beta _j^2\beta _{j'}^2}\right)\\  
&= \left\| {\bf{A}} \right\|_F^2 - {\left\| \beta  \right\|^4},    
\end{align*}
where  $\|\textbf{A}\|_F^{2}$ is the Frobenius norm of $\textbf{A}.$ Thus, by rewriting (\ref{naive_var_form}) the variance of the naive estimator  is given by
\begin{equation}\label{var_naive_1}
{\var} \left( {{{\hat \tau }^2}} \right) = \frac{{4\left( {n - 2} \right)}}{{n\left( {n - 1} \right)}}\left[ {{\beta ^T}{\bf{A}}\beta  - {{\left\| \beta  \right\|}^4}} \right] + \frac{2}{{n\left( {n - 1} \right)}}\left[ {\left\| {\bf{A}} \right\|_F^2 - {{\left\| \beta  \right\|}^4}} \right].
\end{equation}

\noindent\textbf{\textit{Proof of Proposition \ref{consistency_naive}}}:\\
Notice that  $\hat\tau^2$ is consistent if $\var[\hat\tau^2]\xrightarrow{n\rightarrow\infty}0$ since $\hat\tau^2$ is unbiased.
Thus, by \eqref{var_naive_1} it is enough to require that $\frac{\beta^T\textbf{A}\beta}{n}\xrightarrow{n\rightarrow\infty}0$ and $\frac{\|\textbf{A}\|^{2}_F}{n^2}\xrightarrow{n\rightarrow\infty}0$. The latter is assumed  and we now show that the former also holds true.
Let $\lambda_1\geq...\geq\lambda_p$ be the eigenvalues of $\textbf{A}$  
and notice that $\textbf{A}$ is symmetric.
We have that $n^{-2}\lambda_1^2 \leq n^{-2}\sum_{j=1}^{p}\lambda_j^2=n^{-2} tr(\textbf{A}^2)=n^{-2}\|\textbf{A}\|_F^2$ and therefore (iii) implies that $\frac{\lambda_1}{n}\xrightarrow{n\rightarrow\infty}0$. 
Now,
$\frac{1}{n}\beta^T\textbf{A}\beta\equiv\frac{1}{n}\|\beta\|^2[(\frac{\beta}{\|\beta\|})^T\textbf{A}\frac{\beta}{\|\beta\|}]\leq\frac{1}{n}\|\beta\|^2\lambda_1\xrightarrow{n\rightarrow\infty}0,$ where the last limit follows from the assumption that $\tau^2=O(1),$ and from the fact that $\frac{\lambda_1}{n}\xrightarrow{n\rightarrow\infty}0$ . We conclude that $\var[\hat\tau^2] \xrightarrow{n\rightarrow\infty} 0$.

We now prove the moreover part, that is, independence of the columns of $\textbf{\emph{X}}$ implies that 
 $\frac{\|\textbf{A}\|^{2}_F}{n^2}\xrightarrow{n\rightarrow\infty}0$.
By definition we have
$ \|\textbf{A}\|_F^{2} =\sum_{j,j'}[E(W_{ij}W_{ij'})]^2.$
Notice that when $j=j'$ we have,
\begin{align*}
E\left( {W_{ij}^2} \right) =& E\left( {X_{ij}^2Y_i^2} \right) = E\left( {X_{ij}^2{{\left[ {{\beta ^T}{X_i} + {\varepsilon _i}} \right]}^2}} \right) = E\left( {X_{ij}^2\left[ {\sum\limits_{k,k'}^{} {{\beta _k}{\beta _{k'}}{X_{ik}}{X_{ik'}}}  + 2{\beta ^T}{X_i}{\varepsilon _i} + \varepsilon _i^2} \right]} \right)
\\ &= E\left( {X_{ij}^2\sum\limits_{k,k'}^{} {{\beta _k}{\beta _{k'}}{X_{ik}}{X_{ik'}}} } \right) + 0 + E\left( {X_{ij}^2\varepsilon _i^2} \right)
\\ &= E\left( {X_{ij}^2\sum\limits_{k = 1}^p {\beta _k^2X_{ik}^2} } \right) + \underbrace {E\left( {X_{ij}^2\sum\limits_{k \ne k'}^{} {{\beta _k}{\beta _{k'}}{X_{ik}}{X_{ik'}}} } \right)}_0 +\sigma^2 E\left( {X_{ij}^2} \right)
\\ &= \beta _j^2E\left( {X_{ij}^4} \right) + \sum\limits_{k \ne j}^p {\beta _k^2\underbrace {E\left( {X_{ik}^2X_{ij}^2} \right)}_1}  + {\sigma ^2}
\\ &= \beta _j^2E\left( {X_{ij}^4} \right) + {\left\| \beta  \right\|^2} - \beta _j^2 + {\sigma ^2} = {\left\| \beta  \right\|^2} + {\sigma ^2} + \beta _j^2\left[ {E\left( {X_{ij}^4 - 1} \right)} \right].
\end{align*}
Notice that $E\left( {X_{ij}^2\sum\limits_{k \ne k'}^{} {{\beta _k}{\beta _{k'}}{X_{ik}}{X_{ik'}}} } \right)=0$ follows from the assumptions that the columns of ${\bf X}$ are independent and  $E(X_{ij})=0$ for each $j$. Also notice that in the third row we used the assumption that $E(\epsilon_i^2|X_i)=\sigma^2.$ 

Similarly, when $j \neq j'$,
\begin{align*}
E\left( {{W_{ij}}{W_{ij}}} \right) = E\left( {{X_{ij}}{X_{ij'}}Y_i^2} \right) &= E\left[ {{X_{ij}}{X_{ij'}}{{\left( {{\beta ^T}{X_i} + {\varepsilon _i}} \right)}^2}} \right] = E\left[ {{X_{ij}}{X_{ij'}}{{\left( {{\beta ^T}{X_i} + {\varepsilon _i}} \right)}^2}} \right]
\\& =E\left[ {{X_{ij}}{X_{ij'}}\left( {\sum\limits_{k,k'}^{} {{\beta _k}{\beta _{k'}}{X_{ik}}{X_{ik'}}}  + 2{\beta ^T}{X_i}{\varepsilon _i} + \varepsilon _i^2} \right)} \right]
\\ &= E\left[ {{X_{ij}}{X_{ij'}}\sum\limits_{k,k'}^{} {{\beta _k}{\beta _{k'}}{X_{ik}}{X_{ik'}}} } \right] + 0 + E\left( {{X_{ij}}{X_{ij'}}\varepsilon _i^2} \right)
\\ &= 2{\beta _j}{\beta _{j'}}E\left( {X_{ij}^2X_{ij'}^2} \right) + 0 + \underbrace {E\left( {{X_{ij}}{X_{ij'}}} \right)}_0E\left( {\varepsilon _i^2} \right) = 2{\beta _j}{\beta _{j'}}E\left( {X_{ij}^2} \right)E\left( {X_{ij'}^2} \right) = 2{\beta _j}{\beta _{j'}}.
 \end{align*}
This can be written more  compactly as
\begin{equation}\label{expectation_WjWj}
 E(W_{ij}W_{ij'})=\begin{cases}
			2{\beta _j}{\beta _{j'}}, &  j\neq j'\\
            \sigma_Y^2   + \beta _j^2[E (X_{ij}^4) - 1]  , & j=j',
		 \end{cases}   
\end{equation}
where $\sigma_Y^2= {\left\| \beta  \right\|^2} + {\sigma ^2}.$
Therefore,
\begin{multline}\label{forbenouos}
 \|\textbf{A}\|_F^{2} = 4\sum_{j\neq j'}\beta_j^2\beta_{j'}^2+\sum_{j}\Big(\sigma_Y^2   + \beta _j^2[E (X_{ij}^4) - 1]\Big)^2 \le 4\|\beta\|^4+\sum_j\Big( \sigma_Y^4+
\beta_j^4[E (X_{ij}^4) - 1]^2
+2\sigma_Y^2\beta _j^2[E (X_{ij}^4) - 1]
\Big)\\
=p\sigma_Y^4+O(1).    
\end{multline}
where the last equality holds since
 $\sigma_Y^2\equiv \tau^2+\sigma^2=O(1)$,
 $E(X_{ij}^4)=O(1)$ and  by the Cauchy–Schwarz inequality we have  $\sum_j \beta_j^4\leq\sum_{j,j'}\beta_j^2\beta_{j'}^2=\|\beta\|^4=O(1).$  
Now  since $p/n = O(1)$ then  $\frac{\|\textbf{A}\|_F^{2}}{n^2}\rightarrow 0$ and we conclude that $\var(\hat\tau^2)=O(\frac{1}{n}),$ i.e., $\hat\tau^2$ is $\sqrt{n}$-consistent.

\begin{remark}\label{example_1}
\noindent\textbf{\textit{Calculations for Example \ref{exmp1}}}:\\
\begin{equation}\label{cov_ta2_g2}
\begin{split}
\text{Cov}[\hat\tau^2,g(X)] & \equiv \text{Cov}\left( {\frac{2}{n\left( {n - 1} \right)}\sum\limits_{{i_1} < {i_2}} {{W_{{i_1}}}{W_{{i_2}}}} ,\frac{1}{n}\sum_{i=1}^{n}[X_i^2-1]} \right)\\
& = \frac{2}{{{n^2}\left( {n - 1} \right)}}\sum\limits_{{i_1} < {i_2}} {\sum\limits_{i = 1}^n {\text{Cov}\left( {{X_{{i_1}}}{Y_{{i_1}}}{X_{{i_2}}}{Y_{{i_2}}},{X_i^2}} \right)} } \\
&= \frac{2}{{{n^2}\left( {n - 1} \right)}}\sum\limits_{{i_1} < {i_2}} {\sum\limits_{i = 1}^n [{E\left( {{X_{{i_1}}}{Y_{{i_1}}}{X_{{i_2}}}{Y_{{i_2}}}{X_i^2})-\beta^2} \right]} } \\
 &= \frac{4}{{{n^2}\left( {n - 1} \right)}}\sum\limits_{{i_1} < {i_2}}[{E\left( {X_{{i_1}}^3{Y_{{i_1}}}} \right)}\beta-\beta^2] \\
 &= \frac{{4\beta }}{{{n^2}\left( {n - 1} \right)}}\sum\limits_{{i_1} < {i_2}}[ {E\left( {X_{{i_1}}^3{Y_{{i_1}}}}\right)-\beta] }  \\
 &= \frac{{4\beta }}{{{n^2}\left( {n - 1} \right)}}\frac{{n\left( {n - 1} \right)}}{2}[E\left( {X_{{i_1}}^3{Y_{{i_1}}}} \right)-\beta]\\
 &= \frac{2\beta}{n} [E\left( {{X^3}Y} \right)-\beta],
 \end{split}
\end{equation}
where in the third equality we used $E(X^2)=1$ and $E(XY) \equiv \beta.$
In the fourth equality  the expectation is zero for all $i \neq i_1,i_2$.
Now, since $X\sim N(0,1)$ and $E(\epsilon|X)=0$, then 
$$E\left( {{X^3}Y} \right) = E\left( {{X^3}\left( { \beta X + \varepsilon } \right)} \right) =  \beta E\left( {{X^4}} \right) = 3\beta.$$
Therefore,
$\text{Cov}[ {\hat \beta{}^2,{{g}}(X)} ]=\frac{4\beta^2}{n}.$ 
Notice that
\begin{equation*}
{\mathop{\var}} \left[ {g\left( X \right)} \right] = {\mathop{\rm var}} \left[ {\frac{1}{n}\sum\limits_{i = 1}^n {\left( {X_i^2 - 1} \right)} } \right] = \frac{1}{n}\left[ {E\left( {{X^4}} \right) - E\left( {{X^2}} \right)} \right] = \frac{2}{n}.
\end{equation*}
Therefore, by \eqref{general_c} we get $c^*=-2\beta^2.$
Plugging-in $c^*$  back in \eqref{variance_change}  yields $\text{Var}(U_{c^*})=\text{Var}(\Hat{\tau}^2)-~\frac{8}{n}\beta^4.$ 
\end{remark}

 \noindent\textbf{\textit{{Proof of Theorem \ref{theorem1}}}}:\\
%The covariance matrix $M \equiv\cov(g)$ is positive definite.\\
 1. We now prove the first direction: OOE $\Rightarrow \cov[R^*,g]=0$ for all $g\in \cal G.$\\
 Let $R^* \equiv T+g^*$ be an OOE for $\theta$ with respect to the family of zero-estimators $\cal G$. By definition, $\var[R^*] \leqslant \var[T+g]$ for all $g \in \cal G$. For every $g=\sum_{k=1}^{m}c_kg_k$,  define $\Tilde{g}\equiv g-g^*=\sum_{k=1}^{m}(c_k-c_k^*)g_k=\sum_{k=1}^{m}\Tilde{c}_kg_k$ for some fixed $m,$ and note that $ \Tilde{g} \in \cal G$. Then,
 \begin{equation*}\label{f_c}
 \begin{split}
 \var[R^*]& \leqslant \var[T+g]=\var[T+g^*+\Tilde{g}]=  \var[R^*+\sum_{k=1}^{m}\Tilde{c}_kg_k]\\
  &= \var[R^*]+2\sum_{k=1}^{m}\Tilde{c}_k\cdot\cov[R^*,g_k]+\var[\sum_{k=1}^{m}\Tilde{c}_kg_k].
  \end{split}
 \end{equation*}
 Therefore, for all $(\tilde{c}_1,...,\tilde{c}_m)$,
 $$0 \leqslant 2\sum_{k=1}^{m}\Tilde{c_k}\cdot\cov[R^*,g_k]+\var[\sum_{k=1}^{m}\Tilde{c_k}g_k],$$ which can be represented compactly as
  \begin{equation}\label{f_c_def}
 0 \leqslant -2\mathbf{\Tilde{c}}^T \mathbf{b} +\var[\mathbf{\Tilde{c}}^T \mathbf{g_m} ]= -2\mathbf{\Tilde{c}}^T\mathbf{b} + \mathbf{\Tilde{c}}^T M \mathbf{\Tilde{c}} \equiv f(\mathbf{\Tilde{c}}),
  \end{equation}
  where $\mathbf{b} \equiv -\left( \cov[R^*,g_1],...,\cov[R^*,g_m]\right)^T$, $\mathbf{g_m} \equiv(g_1,...,g_m)^T$, $M=\cov[\mathbf{g_m}]$ and $\mathbf{\Tilde{c}}\equiv(\Tilde{c}_1,...,\Tilde{c}_m)^T.$
Notice that $f(\mathbf{\Tilde{c}})$ is a convex function in $\mathbf{\Tilde{c}}$ that satisfies $f(\mathbf{\Tilde{c}}) \geq 0$ for all $\mathbf{\Tilde{c}}$. 
%In particular, 
%\begin{equation}\label{f__star}
 %f(\mathbf{\Tilde{c}_{min}})= -2\mathbf{\Tilde{c}_{min}}^T\mathbf{b} \mathbf{\Tilde{c}_{min}}^T M \mathbf{\Tilde{c}_{min}} \geq 0   
%\end{equation}
%where $\mathbf{\Tilde{c}_{min}}$ is the minimizer of $f(\mathbf{\Tilde{c}}).$
Differentiate $f(\mathbf{\Tilde{c}})$ in order to find its minimum
 \begin{equation*}
     \nabla f(\mathbf{\Tilde{c}})=-2 \mathbf{b}  +2M\mathbf{\Tilde{c}}=0.
 \end{equation*} 
    %\begin{equation*}
 %\frac{\partial f(\mathbf{\Tilde{c}})}{\partial \mathbf{\Tilde{c}}}=-2 \mathbf{b}  +2M\mathbf{\Tilde{c}}=0.
 %\end{equation*}
  Assuming $M$ is positive definite and solving for $\mathbf{\Tilde{c}}$ yields the minimizer $\mathbf{\Tilde{c}_{min}}=M^{-1} \mathbf{b}$. 
Plug-in $\mathbf{\Tilde{c}_{min}}$ in the (\ref{f_c_def}) yields
\begin{equation}\label{f_c_star}
\begin{split}
f(\mathbf{\Tilde{c}_{min}}) \equiv -2(M^{-1}\mathbf{b})^T\mathbf{b}+(M^{-1}\mathbf{b})^TM(M^{-1}\mathbf{b})=
 -\mathbf{b}^TM^{-1}\mathbf{b} & \geq 0.
\end{split}    
\end{equation}   

Since, by assumption,  $M$ is positive definite,  so is $M^{-1}$, i.e., $\mathbf{b}^TM^{-1}\mathbf{b}>0$ for all non-zero $\mathbf{b} \in \mathbb{R}^m.$
Thus, (\ref{f_c_star}) is satisfied only if $\mathbf{b}\equiv \mathbf0$, i.e., $\cov[R^*,\mathbf{g_m}]=\mathbf 0$ which also implies $\cov[R^*,\sum_{k=1}^{m}c_kg_k]=0$ for any $c_1,...,c_m \in \mathbb{R}$. Therefore, $\cov[R^*,g]=0$ for all $g \in \cal G$. \\
 2. We now prove the other direction: if $R^*$ is uncorrelated with all zero-estimators of a given family $\cal G$ then it is an OOE.\\
 Let $R^*=T+g^*$ and $R\equiv T+g$ be  unbiased estimators of $\theta$, where $g^*,g \in \cal G$. 
  Define  $\Tilde{g} \equiv R^*-R=g^*-g$ and notice that $\Tilde{g}\in\cal G$ . Since by assumption $R^*$ is uncorrelated  with $\Tilde{g}$,
 \begin{equation*}
     \begin{split}
     0=\cov[R^*,\Tilde{g}] \equiv\cov[R^*,R^*-R]=\var[R^*]-\cov[R^*,R],
       \end{split}
 \end{equation*}
 and hence $\var[R^*]=\cov[R^*,R]$. By the Cauchy–Schwarz inequality, $(\cov[R^*,R])^2 \leq \var[R^*] \var[R]$,  we conclude that  $\var[R^*] \leq \var[R]=\var[T+g]$ for all $g\in\cal G$.  \qed\\

\noindent\textbf{\textit{Proof of Theorem  \ref{oracle_p}:}}\\
 We start by proving Theorem \ref{oracle_p} for the special case of $p=2$ and then generalize  for  $p>2$. 
By Theorem \ref{theorem1} we need to show that $\cov\left( {{T_{oracle}},{g_{{k_1}{k_2}}}} \right) = 0$  for all \(\left( {{k_1},{k_2}} \right) \in {\mathbb{N}_0^2}\)
where \(\)\({g_{{k_1}{k_2}}} = \frac{1}{n}\sum\limits_{i = 1}^n {\left[ {X_{i1}^{{k_1}}X_{i2}^{{k_2}} - E\left( {X_{i1}^{{k_1}}X_{i2}^{{k_2}}} \right)} \right]} \).
Write, 
 \[\cov\left( {{T_{oracle}},{g_{{k_1}{k_2}}}} \right) = \cov\left( {{{\hat \tau }^2} - 2\sum\limits_{j = 1}^2 {\sum\limits_{j' = 1}^2 {{\psi _{jj'}}} } ,{g_{{k_1}{k_2}}}} \right) = \cov\left( {{{\hat \tau }^2},{g_{{k_1}{k_2}}}} \right) - 2\sum\limits_{j = 1}^2 {\sum\limits_{j' = 1}^2 {\cov \left( {{\psi _{jj'}},{g_{{k_1}{k_2}}}} \right)} } .\] 
 Thus, we need to show that 
 \begin{equation}\label{what_to_show}
 \cov\left( {{{\hat \tau }^2},{g_{{k_1}{k_2}}}} \right) = 2\sum\limits_{j = 1}^2 {\sum\limits_{j' = 1}^2 {\cov \left( {{\psi _{jj'}},{g_{{k_1}{k_2}}}} \right)} }.    
 \end{equation}

 We start with calculating the LHS of ~\eqref{what_to_show}, namely   $\cov\left( {{{\hat \tau }^2},{g_{{k_1}{k_2}}}} \right)$. Recall that $\hat{\tau}^2\equiv\hat{\beta}_1^2+\hat{\beta}_2^2$ and therefore $\cov[\hat{\tau}^2,g_{k_1k_2}]=\cov(\hat{\beta}_1^2,g_{k_1k_2})+\cov(\hat{\beta}_2^2,g_{k_1k_2}).$
Now, for all $(k_1,k_2)\in\mathbb{N}_0^2$ we have
\begin{equation}\label{cov_bet_gk1k2}
\begin{split}
\text{Cov}[\hat\beta_1^2, g_{k_1k_2}]&\equiv\text{Cov}\left(\frac{2}{n(n-1)}\sum_{i_1<i_2}W_{i_11}W_{i_21},\frac{1}{n}\sum_{i=1}^{n}(X_{i1}^{k_1}X_{i2}^{k_2}-E[X_{i1}^{k_1}X_{i2}^{k_2}])\right)\\
&=\frac{2}{n^2(n-1)}\sum_{i_1<i_2}\sum_{i=1}^{n}\text{Cov}\left(X_{i_11}Y_{i_1}X_{i_21}Y_{i_2},X_{i1}^{k_1}X_{i2}^{k_2} \right)\\
&=\frac{2}{n^2(n-1)}\sum_{i_1<i_2}\sum_{i=1}^{n}\left (E[X_{i_11}Y_{i_1}X_{i_21}Y_{i_2}X_{i1}^{k_1}X_{i2}^{k_2}]-\beta_1^2E[X_{i1}^{k_1}X_{i2}^{k_2}]
\right)\\ 
&=\frac{4}{n^2(n-1)}\sum_{i_1<i_2}\left (E[X_{i_11}Y_{i_1}X_{i_21}Y_{i_2}X_{i_11}^{k_1}X_{i_12}^{k_2}]-\beta_1^2E[X_{i_11}^{k_1}X_{i_12}^{k_2}]
\right)\\
&=\frac{4}{n^2(n-1)}\sum_{i_1<i_2}\left (E[X_{i_11}^{k_1+1}Y_{i_1}X_{i_12}^{k_2}]E[X_{i_21}Y_{i_2}]-\beta_1^2E[X_{i_11}^{k_1}X_{i_12}^{k_2}]
\right)\\
&=\frac{4}{n^2(n-1)}\sum_{i_1<i_2}\left (E[X_{i_11}^{k_1+1}Y_{i_1}X_{i_12}^{k_2}]\beta_1-\beta_1^2E[X_{i_11}^{k_1}X_{i_12}^{k_2}]
\right)\\
&=\frac{4}{n^2(n-1)}\frac{n(n-1)}{2}\left (E[X_{11}^{k_1+1}Y_{1}X_{12}^{k_2}]\beta_1-\beta_1^2E[X_{11}^{k_1}X_{12}^{k_2}]
\right)\\
&=\frac{2}{n}\left (E[X_{11}^{k_1+1}Y_{1}X_{12}^{k_2}]\beta_1-\beta_1^2E[X_{11}^{k_1}X_{12}^{k_2}]
\right),
\end{split}
\end{equation}
where the calculations can be justified by  similar arguments to those presented in (\ref{cov_ta2_g2}). 
We shall use the following notation:\\
\begin{equation*}
    \begin{split}
        A &\equiv E\left[ {X_{11}^{{k_1} + 2}}  {X_{12}^{{k_2}}} \right]\\
        B &\equiv E\left[ {X_{11}^{{k_1} + 1}}  {X_{12}^{{k_2} + 1}} \right] \\
        C &\equiv E\left[ {X_{11}^{{k_1}}}  {X_{12}^{{k_2}}} \right] \\
        D&\equiv E\left[ {X_{11}^{{k_1}}}  {X_{12}^{{k_2} + 2}} \right].     \end{split}
\end{equation*}
Notice that $A, B, C$ and $D$ are functions of $(k_1,k_2)$ but this is suppressed in the notation.
Write,
\begin{equation*}\label{add_linear}
\begin{split}
E[X_{11}^{k_1+1}X_{12}^{k_2}Y_1]&=E[X_{11}^{k_1+1}X_{12}^{k_2}(\beta_1X_{11}+\beta_2X_{12}+\epsilon_1)]\\
&=\beta_1E[X_{11}^{k_1+2}X_{12}^{k_2}]+\beta_2E[X_{11}^{k_1+1}X_{12}^{k_2+1}] = \beta_1A+\beta_2B.
\end{split}    
\end{equation*}
Thus, rewrite (\ref{cov_bet_gk1k2}) and obtain
\begin{equation}\label{cov1}
    \text{Cov}[\hat\beta_1^2, g_{k_1k_2}]=
\frac{2}{n}\left( {\left[ {{\beta _1}A + {\beta _2}B} \right]{\beta _1} - \beta _1^2C} \right).
\end{equation}
Similarly, by symmetry, 
% \begin{equation*}
%     \begin{split}
% &\text{Cov}[\hat\beta_2^2, g_{k_1k_2}]=
% \frac{2}{n}\Big ([\beta_2E[X_{i2}^{k_2+2}]E[X_{i1}^{k_1}]+\beta_1E[X_{i2}^{k_2+1}]E[X_{i1}^{k_1+1}]]\beta_2
% -\beta_2^2E[X_{i1}^{k_1}]E[X_{i2}^{k_2}]
% \Big).
%     \end{split}
% \end{equation*}
% Therefore,
\begin{equation}\label{cov2}
\cov[\hat\beta_2^2, g_{k_1k_2}]=\frac{2}{n}\left( {\left[ {{\beta _2}D + {\beta _1}B} \right]{\beta _2} - \beta _2^2C} \right).
\end{equation}
Using (\ref{cov1}) and (\ref{cov2}) we get
\begin{equation}\label{1_part_show}
    \begin{split}
\cov[\hat{\tau}^2,g_{k_1k_2}]&=\cov(\hat{\beta}_1^2,g_{k_1k_2})+\cov(\hat{\beta}_2^2,g_{k_1k_2}) \\
&=\frac{2}{n}\left( {\left[ {{\beta _1}A + {\beta _2}B} \right]{\beta _1} - \beta _1^2C + \left[ {{\beta _2}D + {\beta _1}B} \right]{\beta _2} - \beta _2^2C} \right) \\
&=\frac{2}{n}\left[ {\overbrace {\beta _1^2A + \beta _2^2D}^{{L_1}} + \overbrace {2{\beta _1}{\beta _2}B}^{{L_2}} - \overbrace {C\left( {\beta _1^2 + \beta _2^2} \right)}^{{L_3}}} \right] = \frac{2}{n}\left( {{L_1} + {L_2} - {L_3}} \right).
    \end{split}
\end{equation}

We now move to calculate the RHS of \eqref{what_to_show}, namely  
 $\sum\limits_{j = 1}^2 {\sum\limits_{j' = 1}^2 {{\cov} \left( {{\psi _{jj'}},{g_{{k_1}{k_2}}}} \right)} }.$
First,  recall that    \({h_{jj}} \equiv \frac{1}{n}\sum\limits_{i = 1}^n {\left[ {{X_{ij}}{X_{ij'}} - E\left( {{X_{ij}}{X_{ij'}}} \right)} \right]} \)  and \({g_{{k_1}{k_2}}} \equiv \frac{1}{n}\sum\limits_{i = 1}^n {\left[ {X_{i1}^{{k_1}}X_{i2}^{{k_2}} - E\left( {X_{i1}^{{k_1}}X_{i2}^{{k_2}}} \right)} \right]} \) where \(\left( {{k_1},{k_2}} \right) \in~\mathbb{N}_0^2\). Hence,  
${h_{11}} \equiv \frac{1}{n}\sum\limits_{i = 1}^n {\left( {X_{i1}^2 - 1} \right)}$ which by definition is also equal to   ${g_{20}}.$ Similarly, we have
${h_{12}} = {h_{21}} \equiv\frac{1}{n}\sum\limits_{i = 1}^n {\left( {{X_{i1}}{X_{i2}}} \right) = {g_{11}}}$
 and ${h_{22}} \equiv \frac{1}{n}\sum\limits_{i = 1}^n {\left( {X_{i2}^2 - 1} \right) = {g_{02}}}.$
Thus,
\begin{multline}\label{snd_term}
 \sum\limits_{j = 1}^2 {\sum\limits_{j' = 1}^2 {{\cov} \left( {{\psi _{jj'}},{g_{{k_1}{k_2}}}} \right)} }  = \sum\limits_{j = 1}^2 {\sum\limits_{j' = 1}^2 {{\beta _j}{\beta _{j'}}{\cov} \left( {{h_{jj'}},{g_{{k_1}{k_2}}}} \right)} } \\
 = \beta _1^2{\cov} \left( {{h_{11}},{g_{{k_1}{k_2}}}} \right) + 2{\beta _1}{\beta _2}{\cov} \left( {{h_{12}},{g_{{k_1}{k_2}}}} \right) + \beta _2^2{\cov} \left( {{h_{22}},{g_{{k_1}{k_2}}}} \right)\\ 
 = \beta _1^2{\cov} \left( {{g_{20}},{g_{{k_1}{k_2}}}} \right) + 2{\beta _1}{\beta _2}{\cov} \left( {{g_{11}},{g_{{k_1}{k_2}}}} \right) + \beta _2^2{\cov} \left( {{g_{02}},{g_{{k_1}{k_2}}}} \right).
 \end{multline}
 Now, observe that for every $(k_1,k_2,d_1
 ,d_2)\in\mathbb{N}_0^4,$
\begin{equation}\label{M_matrix}
    \begin{split}
\text{Cov}[g_{k_1k_2},g_{d_1d_2}]
&=\text{Cov}\Big(\frac{1}{n}\sum_{i=1}^{n}[X_{i1}^{k_1}X_{i2}^{k_2}-E(X_{i1}^{k_1}X_{i2}^{k_2})],\frac{1}{n}\sum_{i=1}^{n}[X_{i1}^{d_1}X_{i2}^{d_2}-E(X_{i1}^{d_1}X_{i2}^{d_2})]\Big)\\
&=n^{-2}\sum_{i_1=1}^{n}\sum_{i_2=1}^{n}(E[X_{i_11}^{k_1}X_{i_12}^{k_2}X_{i_21}^{d_1}X_{i_22}^{d_2}]-E[X_{i_11}^{k_1}X_{i_12}^{k_2}]E[X_{i_21}^{d_1}X_{i_22}^{d_2}])\\
&=\frac{1}{n}\Big( E[X_{11}^{k_1+d_1}X_{12}^{k_2+d_2}]-E[X_{11}^{k_1}X_{12}^{k_2}]E[X_{11}^{d_1}X_{12}^{d_2}]\Big),
    \end{split}
\end{equation}
where the third equality holds since the terms with $i_1 \neq i_2$ vanish.
It follows from (\ref{M_matrix}) that
 \begin{equation*}\label{h1_h_5}
 \begin{array}{l}
\cov[g_{k_1k_2},g_{20}]=  \frac{1}{n}\left( E[{X_{11}^{{k_1} + 2} {X_{12}^{{k_2}}}  } ] - E[ {X_{11}^{{k_1}}}  {X_{12}^{{k_2}}}]  \right) =\frac{1}{n}(A-C),\\
\cov[g_{k_1k_2},g_{11}] =   \frac{1}{n}E[ {X_{11}^{{k_1} + 1}}  {X_{12}^{{k_2} + 1}} ] = \frac{B}{n},\\
\cov[g_{k_1k_2},g_{02}] =
\frac{1}{n}\left(E[ {X_{11}^{{k_1}}}  {X_{12}^{{k_2} + 2}} ] - E[ {X_{11}^{{k_1}}}  {X_{12}^{{k_2}}}]  \right) = \frac{1}{n}(D-C)\,\,\,\,\,\,\\
\end{array}   
\end{equation*}
Therefore, rewrite (\ref{snd_term}) to get
\begin{equation}\label{2_part_show}2\sum\limits_{j = 1}^2 {\sum\limits_{j' = 1}^2 {\cov\left( {{\psi _{jj'}},{g_{{k_1}{k_2}}}} \right)} }  = \frac{2}{n}\left[ {\overbrace {\beta _1^2A + \beta _2^2D}^{{L_1}} + \overbrace {2{\beta _1}{\beta _2}B}^{{L_2}} - \overbrace {C\left( {\beta _1^2 + \beta _2^2} \right)}^{{L_3}}} \right] = \frac{2}{n}\left( {{L_1} + {L_2} - {L_3}} \right)
\end{equation}
which is exactly the same expression as in \eqref{1_part_show}.
 Hence, equation \eqref{what_to_show} follows which  completes the proof of Theorm~\ref{oracle_p} for $p=2.$

We now generalize the proof for $p>2$.
Similarly to (\ref{what_to_show}) we want to show that
\begin{equation}\label{want_to_show_p}
\cov\left( {{{\hat \tau }^2},{g_{{k_1}...{k_p}}}} \right) = 2\sum\limits_{j = 1}^p {\sum\limits_{j' = 1}^p {\cov \left( {{\psi _{jj'}},{g_{{k_1}...{k_p}}}} \right)} }.
\end{equation}

We begin by calculating  the LHS of \eqref{want_to_show_p}, i.e., the covariance between $\hat{\tau}^2$ and $g_{k_1...k_p}.$ 
By the same type of calculations as in (\ref{cov_bet_gk1k2}), for all $(k_1,...,k_p)\in\mathbb{N}_0^p$ we have
\begin{multline*}
\cov \left[ {\hat \beta _j^2,{g_{{k_1},...,k_p}}} \right] = \\\frac{2}{n}\left\{ {\left[ {{\beta _j}E\left( {X_{1j}^{{k_j} + 2}\prod\limits_{m \ne j}^{} {X_{1m}^{{k_m}}} } \right) + \sum\limits_{j \ne j'}^{} {{\beta _{j'}}E\left( {X_{1j}^{{k_j} + 1}X_{1j'}^{{k_{j'}} + 1}\prod\limits_{m \ne j,j'}^{} {X_{1m}^{{k_m}}} } \right)} } \right]{\beta _j} - \beta _j^2 {E\left( \prod\limits_{m = 1}^p{X_{1m}^{{k_m}}} \right)} } \right\}
\end{multline*}
Summing the above expressions for $j=1,\ldots,p$, yields
\begin{equation}\label{cov_tau2_zero_estimates}
\begin{split}
 \cov\left[ {{{\hat \tau }^2},{g_{{k_1},....,{k_p}}}} \right]&=\sum_{j=1}^{p}\cov\left[ {\hat \beta _j^2,{g_{{k_1},....,{k_p}}}} \right]  \\
  &= \frac{2}{n}\sum\limits_{j = 1}^p {\beta _j^2} E\left( {X_{1j}^{{k_j}
 + 2}} \prod\limits_{m \ne j}^{}  {X_{1m}^{{k_m}}} \right) \\
 &+ \frac{2}{n}\sum\limits_{j \neq j'}^{} {{\beta _j}{\beta _{j'}}} E\left( {X_{1j}^{{k_j} + 1}}  {X_{1j'}^{{k_{j'}} + 1}} \prod\limits_{m \ne j,j'}^{} { {X_{1m}^{{k_m}}} }\right) \\
 &- \frac{2}{n}\sum\limits_{j = 1}^p {\beta _j^2 {E\left(\prod\limits_{m = 1}^p {X_{1m}^{{k_m}}} \right)} }\\
 &\equiv\frac{2}{n} \big({L_1} + {L_2} -{L_3}\big),
\end{split}    
\end{equation}
where $L_1,L_2$ and $L_3$ are just a generalization of the notation given in \eqref{1_part_show}. 
 Again, notice that $L_1,L_2$ and $L_3$ are functions of $k_1,...,k_p$ but this is suppressed in the notation.

We now move to calculate the RHS of \eqref{want_to_show_p}, namely $2\sum\limits_{j = 1}^p {\sum\limits_{j' = 1}^p {\cov \left( {{\psi _{jj'}},{g_{{k_1}...{k_p}}}} \right)} }
$. Since $\psi_{jj'} = \beta_j\beta _{j'}h_{jj'}$ we have,
\begin{equation}\label{whatever}
 \sum\limits_{j = 1}^p {\sum\limits_{j' = 1}^p {{\cov} \left( {{\psi _{jj'}},{g_{{k_1}...{k_p}}}} \right)} }  = \sum\limits_{j = 1}^p {\sum\limits_{j' = 1}^p {{\beta _j}{\beta _{j'}}{\cov} \left( {{h_{jj'}},{g_{{k_1}...{k_p}}}} \right)} }.   
\end{equation}
 Again, notice the relationship between $h_{jj'}$ and $g_{k_1...k_p}:$ when $j=j'$ we have \({h_{jj}} \equiv \frac{1}{n}\sum\limits_{i = 1}^n {\left( {X_{ij}^2 - 1} \right)}=g_{0...2...0}, \) (i.e., the $j$-th entry is 2 and all others are 0), and for $j \neq j'$ we have \({h_{jj'}} \equiv \frac{1}{n}\sum\limits_{i = 1}^n {{X_{ij}}{X_{ij'}}} = g_{0...1...1...0},\) (i.e., the $j$-th and $j$'-th entries are 1 and all other entries are~0). Hence,
\begin{multline}\label{whatever1}
\sum\limits_{j = 1}^p {\sum\limits_{j' = 1}^p {{\cov} \left( {{\psi _{jj'}},{g_{{k_1}...{k_p}}}} \right) = } } \sum\limits_{j = 1}^p {\sum\limits_{j' = 1}^p {{\beta _j}{\beta _{j'}}{\cov} \left( {{h_{jj'}},{g_{{k_1}...{k_p}}}} \right)} } \\ = \sum\limits_{j = 1}^p {\beta _j^2{\cov} \left( {{g_{0...2...0}},{g_{{k_1}...{k_p}}}} \right)}  + \sum\limits_{j \ne j'}^{} {{\beta _j}{\beta _{j'}}{\cov} \left( {{g_{0...1...1...0}},{g_{{k_1}...{k_p}}}} \right)}.    \end{multline}
Now, similar to (\ref{M_matrix}), for all pairs of index vectors $(k_1,...,k_p)\in\mathbb{N}_0^p$, and $(k_1',...,k_p')\in\mathbb{N}_0^{p}$
\begin{equation}
{\mathop\cov} \left( {{g_{{k_1},...,{k_p}}},{g_{k_1',...,k_p'}}} \right) = \frac{1}{n}\left\{ {E\left( {\prod\limits_{j = 1}^p {X_{1j}^{{k_j} + {k_{j}'}}} } \right) - E\left( {\prod\limits_{j = 1}^p {X_{1j}^{{k_j}}} } \right)E\left( {\prod\limits_{j = 1}^p {X_{1j}^{k_j'}} } \right)} \right\}
\end{equation}
This implies that   
\[\cov\left[ {g_{0...2...0},{g_{{k_1},...,{k_p}}}} \right] = \frac{1}{n}\left[ {E\left( {X_{1j}^{{k_j} + 2}} \prod\limits_{m \ne j} {X_{1m}^{{k_m}}} \right) -   {E\left( \prod\limits_{m = 1}^p{X_{1m}^{{k_m}}} \right)} } \right]\]
and 
\[\cov\left[ {{g_{{0...1...1...0}}},{g_{{k_1},...,{k_p}}}} \right] = \frac{1}{n} {E\left( {X_{1j}^{{k_j} + 1}}  {X_{1j'}^{{k_{j'}} + 1}} \prod\limits_{m \ne j,j'}^{} { {X_{1m}^{{k_m}}} }\right)}  .\]
Hence, rewrite (\ref{whatever1}) to see that
\begin{multline}\label{want_to_show_1p}
2\sum\limits_{j = 1}^p {\sum\limits_{j' = 1}^p {\cov \left( {{\psi _{jj'}},{g_{{k_1}...{k_p}}}} \right)  } } = 2\sum\limits_{j = 1}^p {\beta _j^2{\cov} \left( {{g_{0...2...0}},{g_{{k_1}...{k_p}}}} \right)}  + 2\sum\limits_{j \ne j'}^{} {{\beta _j}{\beta _{j'}}{\cov} \left( {{g_{0...1...1...0}},{g_{{k_1}...{k_p}}}} \right)} \\
= \frac{2}{n} \sum\limits_{j = 1}^p \beta _j^2\left[ {E\left( X_{1j}^{{k_j} + 2} \prod\limits_{m \ne j}  {X_{1m}^{{k_m}}} \right) -   {E\left( \prod\limits_{m = 1}^p{X_{1m}^{{k_m}}} \right)} } \right]  +\\
\frac{2}{n}\sum\limits_{j \ne j'} {\beta _j}{\beta _{j'}} E\left( {X_{1j}^{{k_j} + 1}}  X_{1j'}^{{k_{j'}} + 1} \prod\limits_{m \ne j,j'}  {X_{1m}^{{k_m}}} \right)  
 = \frac{2}{n}\big(L_1-L_3+L_2\big),   
\end{multline}
which is exactly the same expression as in \eqref{cov_tau2_zero_estimates}. Hence,  equation \eqref{want_to_show_p} follows  which completes the proof of Theorem~\ref{oracle_p}.\qed

 \noindent\textbf{\textit{Proof of Corollary~\ref{V_T_orac}}}:\\ 
Write,
\begin{equation}\label{var_t_oracle_expand}
{\var} \left( {{T_{oracle}}} \right) = \var\left( {{{\hat \tau }^2} - 2\sum\limits_{j,j'}^{} {{\psi _{jj'}}} } \right) = {\var} \left( {{{\hat \tau }^2}} \right) - 4\sum\limits_{j,j'}^{} {{\beta _j}{\beta _{j'}}{\cov} \left( {{{\hat \tau }^2},{h_{jj'}}} \right) + 4{\var} \left( {\sum\limits_{j,j'}^{} {{\psi _{jj'}}} } \right)}. \end{equation}
Consider 
\(\sum\limits_{j,j'}^{} {{\beta _j}{\beta _{j'}}{\cov} \left( {{{\hat \tau }^2},{h_{jj'}}} \right)}. \)
We have
\begin{align*}
    \sum\limits_{j,j'}^{} {{\beta _j}{\beta _{j'}}{\cov} \left( {{{\hat \tau }^2},{h_{jj'}}} \right)}  &= \sum\limits_{j = 1}^p {\beta _j^2{\cov} \left( {{{\hat \tau }^2},{h_{jj}}} \right) + } \sum\limits_{j \ne j'}^{} {{\beta _j}{\beta _{j'}}{\cov} \left( {{{\hat \tau }^2},{h_{jj'}}} \right)}\\  
    &= \sum\limits_{j = 1}^p {\beta _j^2{\cov} \left( {{{\hat \tau }^2},{g_{0...2...0}}} \right) + } \sum\limits_{j \ne j'}^{} {{\beta _j}{\beta _{j'}}{\cov} \left( {{{\hat \tau }^2},{g_{0...1...1...0}}} \right)}\\
    &= \sum\limits_{j = 1}^p {\beta _j^2\left[ \frac{{2\beta _j^2}}{n}\left( {E\left( {X_{1j}^4} \right) - 1} \right) \right] + } \sum\limits_{j \ne j'}^{} {{\beta _j}{\beta _{j'}}\left[ {\frac{4}{n}{\beta _j}{\beta _{j'}}} \right]} \\ 
    &= \frac{2}{n}\sum\limits_{j = 1}^p {\beta _j^4\left[ {\left( {E\left( {X_{1j}^4} \right) - 1} \right)} \right] + } \frac{4}{n}\sum\limits_{j \ne j'}^{} {\beta _j^2\beta _{j'}^2} 
       \end{align*}
    where the second and third equality are justified by \eqref{whatever1} and \eqref{cov_tau2_zero_estimates} respectively.
Consider now 
\({\var} \left( {\sum\limits_{j,j'}^{} {{\psi _{jj'}}} } \right)\).
Write,
\begin{equation*}
    \begin{split}
{\var} \left( {\sum\limits_{j,j'}^{} {{\psi _{jj'}}} } \right) & = {\cov} \left( {\sum\limits_{j,j'}^{} {{\beta _j}{\beta _{j'}}{h_{jj'}}} ,\sum\limits_{j,j'}^{} {{\beta _j}{\beta _{j'}}{h_{jj'}}} } \right)\\ 
&= \sum\limits_{{j_1},{j_2},{j_3}{j_4}}^{} {{\beta _{{j_1}}}{\beta _{{j_2}}}{\beta _{{j_3}}}{\beta _{{j_4}}}{\cov} \left( {{h_{{j_1}{j_2}}},{h_{{j_3}{j_4}}}} \right)}  \\ &=\frac{1}{{{n^2}}}\sum\limits_{{j_1},{j_2},{j_3},{j_4}}^{} {{\beta _{{j_1}}}{\beta _{{j_2}}}{\beta _{{j_3}}}{\beta _{{j_4}}}\sum\limits_{{i_1},{i_2}}^{} {} {\cov} \left( {{X_{{i_1}{j_1}}}{X_{{i_1}{j_2}}},{X_{{i_2}{j_3}}}{X_{{i_2}{j_4}}}} \right)}\\
&=\frac{1}{{{n^2}}}\sum\limits_{{j_1},{j_2},{j_3},{j_4}}^{} {{\beta _{{j_1}}}{\beta _{{j_2}}}{\beta _{{j_3}}}{\beta _{{j_4}}}\sum\limits_{{i_1},{i_2}}^{} {} \left[ {E\left( {{X_{{i_1}{j_1}}}{X_{{i_1}{j_2}}}{X_{{i_2}{j_3}}}{X_{{i_2}{j_4}}}} \right) - E\left( {{X_{{i_1}{j_1}}}{X_{{i_1}{j_2}}}} \right)E\left( {{X_{{i_2}{j_3}}}{X_{{i_2}{j_4}}}} \right)} \right]}   \\
&=n^{-2}\sum\limits_{{j_1},{j_2},{j_3},{j_4}}^{} {{\beta _{{j_1}}}{\beta _{{j_2}}}{\beta _{{j_3}}}{\beta _{{j_4}}}\sum\limits_{{i=1}}^{n}\left[ {E\left( {{X_{i{j_1}}}{X_{i{j_2}}}{X_{i{j_3}}}{X_{i{j_4}}}} \right) - E\left( {{X_{i{j_1}}}{X_{i{j_2}}}} \right)E\left( {{X_{i{j_3}}}{X_{i{j_4}}}} \right)} \right]}         
\\
&=\frac{1}{n}\sum\limits_{{j_1},{j_2},{j_3},{j_4}}^{} {{\beta _{{j_1}}}{\beta _{{j_2}}}{\beta _{{j_3}}}{\beta _{{j_4}}}\left[ {E\left( {{X_{1{j_1}}}{X_{1{j_2}}}{X_{1{j_3}}}{X_{1{j_4}}}} \right) - E\left( {{X_{1{j_1}}}{X_{1{j_2}}}} \right)E\left( {{X_{1{j_3}}}{X_{1{j_4}}}} \right)} \right]},         
    \end{split}
\end{equation*}
 where the fifth equality holds since the summand is 0 for all $i_1 \neq i_2$. 
  The summation is not zero in only three cases:\\
  1) $j_1=j_4 \neq j_2=j_3$ \\
  2) $j_1=j_3 \neq j_2=j_4$\\
  3) $j_1=j_2=j_3=j_4.$\\
For the first two cases the summation equals  
\(\frac{1}{n}\sum\limits_{j \ne j'}^{} {\beta _j^2\beta _{j'}^2}. \) For the third case the summation equals to \(\frac{1}{n}\sum\limits_{j = 1}^n {\beta _j^4\left[ {E\left( {X_{1j}^4 - 1} \right)} \right]} \). Overall we have
\[{\var} \left( {\sum\limits_{j,j'}^{} {{\psi _{jj'}}} } \right) = \overbrace {\frac{1}{n}\sum\limits_{j \ne j'}^{} {\beta _j^2\beta _{j'}^2} }^{case\,1} + \overbrace {\frac{1}{n}\sum\limits_{j \ne j'}^{} {\beta _j^2\beta _{j'}^2} }^{case\,2} + \overbrace {\frac{1}{n}\sum\limits_{j = 1}^n {\beta _j^4\left[ {E\left( {X_{1j}^4 - 1} \right)} \right]} }^{case\,3}.\]
Rewrite (\ref{var_t_oracle_expand}) to get
\begin{align*}
{\var} \left( {{T_{oracle}}} \right) &= {\var} \left( {{{\hat \tau }^2}} \right) - 4\sum\limits_{j,j'}^{} {{\beta _j}{\beta _{j'}}{\cov} \left( {{{\hat \tau }^2},{h_{jj'}}} \right)}  + 4{\var} \left( {\sum\limits_{j,j'}^{} {{\psi _{jj'}}} } \right)\\ 
&= {\var} \left( {{{\hat \tau }^2}} \right) - 4\left[ \frac{2}{n}\sum\limits_{j = 1}^p {\beta _j^4\left[ {\left( {E\left( {X_{1j}^4} \right) - 1} \right)} \right] + } \frac{4}{n}\sum\limits_{j \ne j'}^{} {\beta _j^2\beta _{j'}^2}  \right] +  \frac{4}{n}\left\{ {2\sum\limits_{j \ne j'}^{} {\beta _j^2\beta _{j'}^2}  + \sum\limits_{j = 1}^p {\beta _j^4\left[ {E\left( {X_{1j}^4 - 1} \right)} \right]} } \right\} \\
&= \var\left( {{{\hat \tau }^2}} \right) - \frac{4}{n}\left\{ {\sum\limits_{j = 1}^p {\beta _j^4\left[ {E\left( {X_{1j}^4 - 1} \right)} \right] + 2\sum\limits_{j \ne j'}^{} {\beta _j^2\beta _{j'}^2} } } \right\}.\qed
\end{align*}

\begin{remark}\label{rem:improve}
\textbf{Calculations for Example \ref{exp_OOE}}.

 Recall that by \eqref{var_naive} we have
\begin{equation*}
{\var} \left( {{{\hat \tau }^2}} \right) = \frac{{4\left( {n - 2} \right)}}{{n\left( {n - 1} \right)}}\left[ {{\beta ^T}{\bf{A}}\beta  - {{\left\| \beta  \right\|}^4}} \right] + \frac{2}{{n\left( {n - 1} \right)}}\left[ {\left\| {\bf{A}} \right\|_F^2 - {{\left\| \beta  \right\|}^4}} \right].
\end{equation*}
Now,  when we assume  standard Gaussian covariates, one can verify that
$  {\beta ^T}{\bf{A}}\beta  - {\left\| \beta  \right\|^4} = \sigma_Y^2{\tau ^2} + {\tau ^4}$
and
$ \left\| {\bf{A}} \right\|_F^2 - {\left\| \beta  \right\|^4} = p{\sigma_Y^4} + 4\sigma_Y^2{\tau ^2} + 3{\tau ^4},$ where $\sigma_Y^2=\sigma^2+\tau^2.$  Thus, in this case  we can write
\begin{equation}\label{var_naive_normal}
\var\left( {{{\hat \tau }^2}} \right) = \frac{4}{n}\left[ {\frac{{\left( {n - 2} \right)}}{{\left( {n - 1} \right)}}\left[ {\sigma_Y^2{\tau ^2} + {\tau ^4}} \right] + \frac{1}{{2\left( {n - 1} \right)}}\left( {p{\sigma_Y^4} + 4\sigma_Y^2{\tau ^2} + 3{\tau ^4}} \right)} \right].    
\end{equation}
Plug-in $\tau^2=\sigma^2=1$ to get
\begin{equation}\label{var_naive_example}
\var(\hat\tau^2)=\frac{20}{n}+O(n^{-2}),    
\end{equation}
 and $\var(T_{oracle})=\var(\hat\tau^2)-\frac{8}{n}\tau^4 = \frac{12}{n}+O(n^{-2})$ by \eqref{var_T_oracle_normal}.
More generally, the asymptotic improvement of $T_{oracle}$ over the naive estimator is:
\begin{align*}
\mathop {\lim }\limits_{n,p \to \infty } \frac{{\var\left( {{{\hat \tau }^2}} \right) - \var\left( {{T_{oracle}}} \right)}}{{\var\left( {{{\hat \tau }^2}} \right)}}
&= \mathop {\lim }\limits_{n,p \to \infty }\frac{{8{\tau ^4}/n}}{{\frac{4}{n}\left[ {\frac{{\left( {n - 2} \right)}}{{\left( {n - 1} \right)}}\left( {\sigma_Y^2{\tau ^2} + {\tau ^4}} \right) + \frac{1}{{2\left( {n - 1} \right)}}\left( {p{\sigma_Y^4} + 4\sigma_Y^2{\tau ^2} + 3{\tau ^4}} \right)} \right]}}\\ 
&= \frac{{2{\tau ^4}}}{{3{\tau ^4} + \frac{{4p{\tau ^4} + 4\sigma_Y^2{\tau ^2} + 3{\tau ^4}}}{{2n}}}} = \frac{2}{{3 + 2\frac{p}{n}}},
 \end{align*}
 where we used  the fact that $\sigma_Y^2= \tau^2+\sigma^2= 2\tau^2$ in the second equality. 
Now, notice that  when $p=n$ then the reduction is $\frac{2}{{3 + 2}} = 40\%$ and when $p/n$ converges to zero, the reduction is $66\%.$
\end{remark}

 \noindent\textbf{\textit{Proof of Proposition~\ref{var_oracle}}}:\\ 
Write,
\begin{align}\label{ntsh1}
 \var\left( T \right) 
 &= \var\left[ {{{\hat \tau }^2} - 2\sum\limits_{j = 1}^p {\sum\limits_{j' = 1}^p {{{\hat \psi }_{jj'}}} } } \right] \nonumber\\
 &= \var\left( {{{\hat \tau }^2}} \right) - 4{\cov} \left( {{{\hat \tau }^2},\sum\limits_{j = 1}^p {\sum\limits_{j' = 1}^p {{{\hat \psi }_{jj'}}} } } \right) + 4\var\left( {\sum\limits_{j = 1}^p {\sum\limits_{j' = 1}^p {{{\hat \psi }_{jj'}}} } } \right).
\end{align}
We start with calculating the middle term. Let $p_n(k)\equiv n(n-1)(n-2)\cdots (n-k).$ Write,
\begin{align}
   &{\cov} \big( {{{\hat \tau }^2},\sum\limits_{j = 1}^p {\sum\limits_{j' = 1}^p {{{\hat \psi }_{jj'}}} } } \big)\nonumber \\
&= {\cov} \left( {\frac{1}{{n\left( {n - 1} \right)}}\sum\limits_{{i_1} \ne {i_2}}^{} {\sum\limits_{j = 1}^p {{W_{{i_1}j}}{W_{{i_2}j}}} } ,\frac{1}{{n\left( {n - 1} \right)\left( {n - 2} \right)}}\sum\limits_{j,j'}^{} {\sum\limits_{{i_1} \ne {i_2} \ne {i_3}}^{} {{W_{{i_1}j}}{W_{{i_2}j'}}\left[ {{X_{{i_3}j}}{X_{{i_3}j'}} - E\left( {{X_{{i_3}j}}{X_{{i_3}j'}}} \right)} \right]} } } \right) \nonumber\\
&= C_n\sum\limits_I {\sum\limits_{J} {{\cov} \left( {{W_{{i_1}{j_1}}}{W_{{i_2}{j_1}}},{W_{{i_3}{j_2}}}{W_{{i_4}{j_3}}}\left[ {{X_{{i_5}{j_2}}}{X_{{i_5}{j_3}}} - E\left( {{X_{{i_5}{j_2}}}{X_{{i_5}{j_3}}}} \right)} \right]} \right)} }, \label{c1} \end{align}
 where $C_n\equiv \frac{1}{p_n(1)\cdot p_n(2)},$
$I$ is the set of all quintuples of indices $(i_1,i_2,i_3,i_4,i_5)$ such that   $i_1 \neq i_2$ and $i_3 \neq i_4 \neq i_5,$ and $J$ is the set of all triples of indices $(j_1,j_2,j_3)$.
For the set $I$, there are $\binom{2}{1} \cdot 3 =6$ different cases to consider when  one of $\{i_1,i_2 \}$ is equal to one of $\{i_3,i_4,i_5\}$, and an additional $ \binom{2}{2} \cdot 3! =6$ cases to consider when two of  $\{i_1,i_2\}$ are equal to two of $\{i_3,i_4,i_5\}$. 
% For the set of indices $(j_1,j_2,j_3),$ we have three different cases to consider:  (1) when all indices are the equal; (2) when all indices different; (3) when only two indices of  $(j_1,j_2,j_3)$ are equal to each other. 
Similarly, for the set~$J$ there are three cases to consider when only two indices of $\{j_1,j_2,j_3\}$ are equal to each other, (e.g., $j_1=j_2 \neq j_3$) ; one case to consider  when no pair of indices is equal to each other %( $j_1 \neq j_2 \neq j_3$ );
% ,i.e., $j_1 \neq j_2 \neq j_3;$ 
and; one case to consider  when all three indices are equal.
Thus, there are total of $(6+6)\times (3+1+1)=60$ cases to consider. Here we demonstrate only one such case. 
% cases consider let \({J_1} = \left\{ {\left( {{j_1},{j_2},{j_3}} \right):{j_2} = {j_3}} \right\}\) and \({J_2} = \left\{ {\left( {{j_1},{j_2},{j_3}} \right):{j_2} \ne {j_3}} \right\}\)
Let \({I_1} = \left\{ {\left( {{i_1},\dots,{i_5}} \right):{i_1} = {i_5} \ne {i_2} \ne {i_3}\, \ne \,{i_4}} \right\}\) and \({J_1} = \left\{ {\left( {{j_1},{j_2},{j_3}} \right):{j_1} = {j_2} = {j_3}} \right\}.\)
% $I_1$ be the set of all quintuples $(i_1,i_2,i_3,i_4,i_5)$ such that  $i_2 = i_5 \neq i_1 \neq i_3 \neq i_4,$ and $J_1$ be the set of all triples $(j_1,j_2,j_3)$ such that $j_1=j_2=j_3.$ 
Write,
% \begin{equation}\label{wwwwx}
% \begin{split}
% &C_n\sum\limits_{I_1} {\sum\limits_{J_1} {{\cov} \left( {{W_{{i_1}{j_1}}}{W_{{i_2}{j_1}}},{W_{{i_3}{j_2}}}{W_{{i_4}{j_3}}}\left[ {{X_{{i_5}{j_2}}}{X_{{i_5}{j_3}}} - E\left( {{X_{{i_5}{j_2}}}{X_{{i_5}{j_3}}}} \right)} \right]} \right)} }\\
% &=C_n\sum\limits_{I_1} {\sum\limits_{J_1} {\cov\left( {{W_{{i_1}{j_1}}}{W_{{i_2}{j_1}}},{W_{{i_3}{j_2}}}{W_{{i_4}{j_2}}}\left[ {X_{{i_1}{j_2}}^2 - 1} \right]} \right)} } \\ 
%  &=C_n\sum\limits_{I_1} {\sum\limits_{J_1} {E({{W_{{i_2}{j_1}}})E({W_{{i_3}{j_2}}})E({W_{{i_4}{j_2}}})} E\left( {{W_{{i_1}{j_1}}}\left[ {X_{{i_1}{j_2}}^2 - 1} \right]} \right)} } \\ 
%  &=C_n \sum\limits_{I_1} {\sum\limits_{J_1} {{\beta _{{j_1}}}\beta _{{j_2}}^2E\left( {{W_{{i_1}{j_1}}}\left[ {X_{{i_1}{j_2}}^2 - 1} \right]} \right)} }.   
% \end{split}
% \end{equation}
\begin{equation}\label{wwwwx}
\begin{split}
&C_n\sum\limits_{I_1} {\sum\limits_{J_1} {{\cov} \left( {{W_{{i_1}{j_1}}}{W_{{i_2}{j_1}}},{W_{{i_3}{j_2}}}{W_{{i_4}{j_3}}}\left[ {{X_{{i_5}{j_2}}}{X_{{i_5}{j_3}}} - E\left( {{X_{{i_5}{j_2}}}{X_{{i_5}{j_3}}}} \right)} \right]} \right)} }\\
&=C_n\sum\limits_{I_1} {\sum\limits_{j=1}^{p} {\cov\left( {{W_{{i_1}{j}}}{W_{{i_2}{j}}},{W_{{i_3}{j}}}{W_{{i_4}{j}}}\left[ {X_{{i_1}{j}}^2 - 1} \right]} \right)} } \\ 
 &=C_n\sum\limits_{I_1} {\sum\limits_{j=1}^{p} {E({{W_{{i_2}{j}}})E({W_{{i_3}{j}}})E({W_{{i_4}{j}}})} E\left( {{W_{{i_1}{j}}}\left[ {X_{{i_1}{j}}^2 - 1} \right]} \right)} } \\ 
 &=C_n \sum\limits_{I_1} {\sum\limits_{j=1}^{p} \beta_{{j}}^3E\left( {{W_{{i}{j}}}\left[ {X_{{i}{j}}^2 - 1} \right]} \right) }.   
\end{split}
\end{equation}
 Now, notice that 
\begin{equation}\label{from1}
\begin{split}
 E\left[ {{W_{ij}}\left( {X_{ij}^2 - 1} \right)} \right] &= E\left[ {{X_{ij}}{Y_i}\left( {X_{ij}^2 - 1} \right)} \right]\\
&= E\left[ {X_{ij}^3\left( {{\beta ^T}X + {\varepsilon _i}} \right)} \right] - {\beta _j}\\
&= {\beta _j}E\left( {X_{ij}^4} \right) - {\beta _j}\\
&= \beta_j[E(X_{ij}^4)-1].
\end{split}
   \end{equation}
%   where the third equality holds since the covariates are standard normal.
Rewrite \eqref{wwwwx} to get
\begin{align*}
 C_n\sum\limits_{I_1} {\sum\limits_{j=1}^{p} {\beta _{{j}}^3E\left[ {W_{i{j}}^{}\left( {X_{i{j}}^2 - 1} \right)} \right]} } 
&= C_n\sum\limits_{I_1} {\sum\limits_{j = 1}^p {\beta _j^3\left(    \beta_j[E(X_{ij}^4)-1]   \right)} }\\
&=\frac{p_n(3)}{p_n(1)\cdot p_n(2)} \sum\limits_{j = 1}^p {\beta _j^4}[E(X_{ij}^4)-1]\\
&=\frac{(n-3)}{n(n-1)}\sum\limits_{j = 1}^p {\beta _j^4}[E(X_{ij}^4)-1]\\
&=\frac{1}{n}\sum\limits_{j = 1}^p {\beta _j^4}[E(X_{ij}^4)-1]+ O(n^{-2}),   
\end{align*} 
where we used \eqref{from1} to justify the first equality.  
% Similarly, let  $J_2$ be the set of all triples $(j_1,j_2,j_3)$ such that $j_2 \neq j_3.$   In this case the covariance over all combination of subsets $I$ and $J_2$  is
% \begin{equation*}
%   \frac{2}{n}\sum\limits_{j \ne j'}^{} {\beta _j^2\beta _{j'}^2}  + O\left( {{n^{ - 2}}} \right)
% \end{equation*}
By the same type of calculation, one can compute the covariance in \eqref{c1} over all 60 and obtain that
\begin{equation}\label{par1}
  {\cov} \big( {{{\hat \tau }^2},\sum\limits_{j = 1}^p {\sum\limits_{j' = 1}^p {{{\hat \psi }_{jj'}}} } }\big) = 
  %\frac{4}{n}{\tau ^4}
  \frac{2}{n}\left\{ {\sum\limits_{j = 1}^p {\beta _j^4\left[ {E\left( {X_{1j}^4 - 1} \right)} \right] + 2\sum\limits_{j \ne j'}^{} {\beta _j^2\beta _{j'}^2} } } \right\}
  + O\left( {{n^{ - 2}}} \right). 
\end{equation}

We now move to calculate the last term of \eqref{ntsh1}. Recall that
\[\hat\psi_{jj'}= \frac{1}{{n\left( {n - 1} \right)\left( {n - 2} \right)}}\sum\limits_{{i_1} \ne {i_2} \ne {i_3}}^{} {{W_{{i_1}j}}{W_{{i_2}j'}}\left[ {{X_{{i_3}j}}{X_{{i_3}j'}} - E\left( {{X_{{i_3}j}}{X_{{i_3}j'}}} \right)} \right]}. \]
Therefore,
\begin{align}\label{var_to_cmp}
&\var\big( {\sum\limits_{j = 1}^p {\sum\limits_{j' = 1}^p {{{\hat \psi }_{jj'}}} } } \big) = \sum\limits_J^{} {{\cov} \left( {{{\hat \psi }_{{j_1}{j_2}}},{{\hat \psi }_{{j_3}{j_4}}}} \right)} \\
&= \frac{1}{{{{\left[ {n\left( {n - 1} \right)\left( {n - 2} \right)} \right]}^2}}}\sum\limits_J^{} {{\cov} \left( {\sum\limits_{{i_1} \ne {i_2} \ne {i_3}}^{} {{W_{{i_1}{j_1}}}{W_{{i_2}{j_2}}}{X_{{i_3}{j_1}}}{X_{{i_3}{j_2}}}} ,\sum\limits_{{i_1} \ne {i_2} \ne {i_3}}^{} {{W_{{i_1}{j_3}}}{W_{{i_2}{j_4}}}{X_{{i_3}{j_3}}}{X_{{i_3}{j_4}}}} } \right)} \nonumber\\ 
&= p_n^{-2}(2)\sum\limits_J^{} {\sum\limits_I^{} {{\cov} \left( {{W_{{i_1}{j_1}}}{W_{{i_2}{j_2}}}{X_{{i_3}{j_1}}}{X_{{i_3}{j_2}}},{W_{{i_4}{j_3}}}{W_{{i_5}{j_4}}}{X_{{i_6}{j_3}}}{X_{{i_6}{j_4}}}} \right)}\,, }\nonumber     
\end{align}
where $J$ is now defined to be the set of all quadruples $(j_1,j_2,j_3,j_4),$ and $I$ is now defined to be the set of  all sextuples $(i_1,...,i_6)$   such that  $i_1 \neq  i_2 \neq i_3$ and $ i_4 \neq i_5 \neq i_6.$
For the set $I$, there are three different cases to consider: (1) when one of $\{i_1,i_2,i_3\}$ is equal to one of $\{i_4,i_5,i_6\}$; (2) when two of $\{i_1,i_2,i_3\}$  are equal to two of $\{i_4,i_5,i_6\};$  and (3) when  $\{i_1,i_2,i_3\}$  are equal to  $\{i_4,i_5,i_6\}.$
There are $\binom{3}{1} \cdot 3=9$ options for the first case, $\binom{3}{2}\cdot 3!=18$ for the second case, and $\binom{3}{3}\cdot 3!=6$ options for the third case.
% For the set of indices $(j_1,j_2,j_3,j_4),$ we have three different cases to consider:  (1) when all indices are the equal; (2) when all indices different; (3) when only two indices of  $(j_1,j_2,j_3)$ are equal to each other. 
For the set $J$, there are five different cases to consider: (1)~when there is only \textit{one}  pair of equal indices (e.g., $j_1=j_2 \neq j_3 \neq j_4$);   
(2) when there are \textit{two}  pairs of equal indices (e.g., $j_1=j_2 \neq j_3=j_4$);
(3) when only three indices are equal  (e.g., $j_1=j_2=j_3 \neq j_4$); (4) when all four indices are equal  and; (5) all four indices are different from each other. Note that there are $\binom{4}{2}=6$ combinations for the first case, $\binom{4}{2}=6$ for the second case, $\binom{4}{3}=4$ combinations for the third case, and a single combination for each of the last two cases. Thus, there are total of $ (9+18+6)\times (6+6+4+1+1)=594.$  Again we demonstrate only one such calculation.
Let \({I_2} = \left\{ {\left( {{i_1},...,{i_6}} \right):{i_1} = {i_4},{i_2} = {i_5},{i_3} = {i_6}} \right\}\)  and \({J_2} = \left\{ {\left( {{j_1},{j_2},{j_3},{j_4}} \right):{j_1} = {j_3} \ne {j_2} = {j_4}} \right\}.\)
In the view of \eqref{var_to_cmp},
\begin{align*}
&p_n^{ - 2}(2)\sum\limits_{{J_2}}^{} {\sum\limits_{{I_2}}^{} {\cov\left( {{W_{{i_1}{j_1}}}{W_{{i_2}{j_2}}}{X_{{i_3}{j_1}}}{X_{{i_3}{j_2}}},{W_{{i_4}{j_3}}}{W_{{i_5}{j_4}}}{X_{{i_6}{j_3}}}{X_{{i_6}{j_4}}}} \right)} }  = \\
&= p_n^{ - 2}(2)\sum\limits_{{J_2}} {\sum\limits_{{I_2}} {\cov\left( {{W_{{i_1}{j_1}}}{W_{{i_2}{j_2}}}{X_{{i_3}{j_1}}}{X_{{i_3}{j_2}}},{W_{{i_1}{j_1}}}{W_{{i_2}{j_2}}}{X_{{i_3}{j_1}}}{X_{{i_3}{j_2}}}} \right)} } \\
&= p_n^{ - 2}(2)\sum\limits_{{J_2}} {\sum\limits_{{I_2}} {E\left( {W_{{i_1}{j_1}}^2} \right)E\left( {W_{{i_2}{j_2}}^2} \right)E\left( {X_{{i_3}{j_1}}^2} \right)E\left( {X_{{i_3}{j_2}}^2} \right)} } \\
&= p_n^{ - 2}(2)\sum\limits_{{J_2}} {\sum\limits_{{I_2}} {\left( {\sigma _Y^2 + \beta _{{j_1}}^2\left\{ {E\left( {X_{i{j_1}}^4} \right) - 1} \right\}} \right)\left( {\sigma _Y^2 + \beta _{{j_2}}^2\left\{ {E\left( {X_{i{j_2}}^4} \right) - 1} \right\}} \right)} } \\
&\le p_n^{ - 2}(2)\sum\limits_{{J_2}} {\sum\limits_{{I_2}} {\left( {\sigma _Y^2 + \beta _{{j_1}}^2\left( {C - 1} \right)} \right)\left( {\sigma _Y^2 + \beta _{{j_2}}^2\left( {C - 1} \right)} \right)} } \\
&= p_n^{ - 1}(2)\sum\limits_{{j_1} \ne {j_2}}^{} {\left[ {\sigma _Y^4 + \sigma _Y^2\left( {C - 1} \right)\left( {\beta _{{j_1}}^2 + \beta _{{j_2}}^2} \right) + {{\left( {C - 1} \right)}^2}\beta _{{j_1}}^2\beta _{{j_2}}^2} \right]} \\
&= p_n^{ - 1}(2)\left[ {p\left( {p - 1} \right)\sigma _Y^4 + \sigma _Y^2\left( {C - 1} \right)\sum\limits_{{j_1} \ne {j_2}}^{} {\left( {\beta _{{j_1}}^2 + \beta _{{j_2}}^2} \right)}  + {{\left( {C - 1} \right)}^2}\sum\limits_{{j_1} \ne {j_2}}^{} {\beta _{{j_1}}^2\beta _{{j_2}}^2} } \right]\\
&\le p_n^{ - 1}(2)\left[ {p\left( {p - 1} \right)\sigma _Y^4 + \sigma _Y^2\left( {C - 1} \right)\left( {2p{\tau ^2}} \right) + {{\left( {C - 1} \right)}^2}{\tau ^4}} \right],
\end{align*}
where the fourth equality we use  \(E\left( {W_{ij}^2} \right) = \sigma_Y^2 + \beta _j^2[E(X_{ij}^2)-1]\), which is given by \eqref{expectation_WjWj}, and in the fifth equality we used the assumption that $E(X_{ij}^4)\leq C$ for some positive $C.$
Since  we assume $p/n= O(1),$ 
the above expression can be further simplified to
\(\frac{{{p^2}{\sigma_Y^4}}}{{{n^3}}}~+~O\left( {{n^{ - 2}}} \right).\)

By  the same type of calculation, one can
compute the covariance in \eqref{var_to_cmp} over all 594 cases and obtain that
\begin{equation}\label{par2}
\var\big( {\sum\limits_{j = 1}^p {\sum\limits_{j' = 1}^p {{{\hat \psi }_{jj'}}} } } \big) = 
\frac{1}{n}\left\{ {\sum\limits_{j = 1}^p {\beta _j^4\left[ {E\left( {X_{1j}^4 - 1} \right)} \right] + 2\sum\limits_{j \ne j'}^{} {\beta _j^2\beta _{j'}^2} } } \right\}
+\frac{{{2p^2}{\sigma_Y^4}}}{{{n^3}}}~+~O\left( {{n^{ - 2}}} \right).    \end{equation}
Lastly, plug-in \eqref{par1} and \eqref{par2} into \eqref{ntsh1} to get
\begin{align*}
\var\left( T \right) 
&=\var\left( {{{\hat \tau }^2}} \right) - 4{\cov} \left( {{{\hat \tau }^2},\sum\limits_{j = 1}^p {\sum\limits_{j' = 1}^p {{{\hat \psi }_{jj'}}} } } \right) + 4\var\left( {\sum\limits_{j = 1}^p {\sum\limits_{j' = 1}^p {{{\hat \psi }_{jj'}}} } } \right)\\
&= \var\left( {{{\hat \tau }^2}} \right) - 4\left(  
  \frac{2}{n}\left\{ {\sum\limits_{j = 1}^p {\beta _j^4\left[ {E\left( {X_{1j}^4 - 1} \right)} \right] + 2\sum\limits_{j \ne j'}^{} {\beta _j^2\beta _{j'}^2} } } \right\}
   \right) \\
  & + 4\left( { \frac{1}{n}\left\{ {\sum\limits_{j = 1}^p {\beta _j^4\left[ {E\left( {X_{1j}^4 - 1} \right)} \right] + 2\sum\limits_{j \ne j'}^{} {\beta _j^2\beta _{j'}^2} } } \right\}
+\frac{{{2p^2}{\sigma_Y^4}}}{{{n^3}}}  } \right) + O\left( {{n^{ - 2}}} \right)\\
&= \var\left( {{{\hat \tau }^2}} \right) -\frac{4}{n} \left\{ {\sum\limits_{j = 1}^p {\beta _j^4\left[ {E\left( {X_{1j}^4 - 1} \right)} \right] + 2\sum\limits_{j \ne j'}^{} {\beta _j^2\beta _{j'}^2} } } \right\} 
+\frac{{{8p^2}{\sigma_Y^4}}}{{{n^3}}} +O( {n^{ - 2}})   \\
&= \var\left( {{T_{oracle}}} \right) + \frac{{4{p^2}{\sigma_Y^4}}}{{{n^3}}} + O\left( {{n^{ - 2}}} \right),    
\end{align*}
where the last equality holds by \eqref{var_T_oracle_normal}. \qed

\begin{remark}\label{c_star_single}
\textbf{\textit{Calculations for equation \ref{eq:c_star}}}:\\ 
Write,
\begin{align*}
\cov \left( {{{\hat \tau }^2},{g_n}} \right) &= \cov\left( {\frac{1}{{n\left( {n - 1} \right)}}\sum\limits_{{i_1} \ne {i_2}}^{} {\sum\limits_{j = 1}^p {{W_{{i_1}j}}{W_{{i_2}j}}} } ,\frac{1}{n}\sum\limits_{i = 1}^n {{g_i}} } \right)\\
&= \frac{1}{{{n^2}\left( {n - 1} \right)}}\sum\limits_{{i_1} \ne {i_2}}^{} {\sum\limits_{j = 1}^p {\sum\limits_{i = 1}^n {E\left( {{W_{{i_1}j}}{W_{{i_2}j}}{g_i}} \right)} } } \\
&= \frac{2}{{{n^2}\left( {n - 1} \right)}}\sum\limits_{{i_1} \ne {i_2}}^{} {\sum\limits_{j = 1}^p {E\left( {{W_{{i_1}j}}{g_{{i_1}}}} \right)E\left( {{W_{{i_2}j}}} \right)} } \\
&= \frac{2}{{{n^2}\left( {n - 1} \right)}}\sum\limits_{{i_1} \ne {i_2}}^{} {\sum\limits_{j = 1}^p {E\left( {{W_{{i_1}j}}{g_{{i_1}}}} \right){\beta _j}} } \\ 
&= \frac{2}{n}\sum\limits_{j = 1}^p {E\left( {{S_{ij}}} \right){\beta _j}} ,
 \end{align*}
where $S_{ij}\equiv W_{ij}g_i$. Also notice that 
\(\var\left( {{g_n}} \right) = \var \left( {\frac{1}{n}\sum\limits_{i = 1}^n {{g_i}} } \right) = \frac{{\var \left( {{g_i}} \right)}}{n}.\)
Thus, by  \eqref{general_c} we get 
$${c^*} = \frac{{{\cov} \left( {{{\hat \tau }^2},{g_n}} \right)}}{{\var \left( {{g_n}} \right)}} = \frac{{2\sum\limits_{j = 1}^p {E\left( {{S_{ij}}} \right){\beta _j}} }}{{\var \left( {{g_i}} \right)}}.$$
\end{remark} 

\begin{remark}\label{rem:improve_singel}
\textbf{\textit{Calculations for Example \ref{exp_singel}}}:\\ 
In order to calculate $\var(T_{c^*})$ we need to calculate the numerator and denominator of \eqref{var_single}. Consider first $\theta_j \equiv E(S_ij).$
Write,
\[\begin{array}{l}{\theta _j} \equiv E\left( {{S_{ij}}} \right) = E\left( {{X_{ij}}{Y_i}{g_i}} \right) = E\left( {{X_{ij}}\left( {{\beta ^T}{X_i} + \varepsilon_i } \right){g_i}} \right) = \\E\left( {{X_{ij}}\left( {\sum\limits_{m = 1}^p {{\beta _m}{X_{im}}}  + \varepsilon_i } \right)\sum\limits_{k < k'}^{} {{X_{ik}}{X_{ik'}}} } \right) = \sum\limits_{m = 1}^p {} \sum\limits_{k < k'}^{} {{\beta _m}E\left( {{X_{ij}}{X_{im}}{X_{ik}}{X_{ik'}}} \right)} \end{array}\]
where in the last equality we used the assumption that $E(\epsilon|X)=0$.
Since the columns of $\bf{X}$ are independent,  the summation is not zero (up to permutations) when \(j = k\) and \(m = k'\). In this case we have 
\begin{equation*}
{\theta _j} = \sum\limits_{m = 1}^p {\sum\limits_{k < k'}^{} {{\beta _m}E\left( {{X_{ij}}{X_{im}}{X_{ik}}{X_{ik'}}} \right)} }  = \sum\limits_{m \neq j}^p {{\beta _m}E\left( {X_{ij}^2X_{im}^2} \right)}  = \sum\limits_{m \neq j}^p {{\beta _m}E\left( {X_{ij}^2} \right)E\left( {X_{im}^2} \right)}  = \sum\limits_{m \neq j}^p {{\beta _m}}.    
\end{equation*}
Notice that in the forth equality we used the assumption that $E(X_{ij}^2)=1$ for all  $j=1,...,p.$ 
Thus,
\begin{equation}\label{numerator_single}
\sum\limits_{j = 1}^p {{\beta _j}E\left( {{S_{ij}}} \right) = } \sum\limits_{j = 1}^p {{\beta _j}\sum\limits_{m \ne j}^p {{\beta _m}}  = } \sum\limits_{j = 1}^p {{\beta _j}\left( {\sum\limits_{m = 1}^p {{\beta _m}}  - {\beta _j}} \right) = {{\left( {\sum\limits_{j = 1}^p {{\beta _j}} } \right)}^2} - \sum\limits_{j = 1}^p {\beta _j^2}  = } {\left( {\sum\limits_{j = 1}^p {{\beta _j}} } \right)^2} - {\tau ^2}.    
\end{equation}

plug-in $\tau^2=1$ and $\beta_j=\frac{1}{\sqrt{p}}$ to get
the numerator of \eqref{var_single}:
$${\left[ {2\sum\limits_{j = 1}^p {{\beta _j}} E\left( {{S_{ij}}} \right)} \right]^2} = 4{\left[ {{{\left( {\sum\limits_{j = 1}^p {{\beta _j}} } \right)}^2} - {\tau ^2}} \right]^2} = 4{\left[ {{{\left( {p\frac{1}{{\sqrt p }}} \right)}^2} - 1} \right]^2} = 4{\left( {{p^2} - 1} \right)^2}.$$
Consider now the denominator of \eqref{var_single}.
Write,
 \begin{equation*}
 \var\left( {{g_i}} \right) = E\left( {g_i^2} \right) = E\left[ {{{\left( {\sum\limits_{j < j'}^{} {{X_{ij}}{X_{ij'}}} } \right)}^2}} \right] = \sum\limits_{{j_1} < {j_2}}^{} {\sum\limits_{{j_3} < {j_4}}^{} {E\left( {{X_{i{j_1}}}{X_{i{j_2}}}{X_{i{j_3}}}{X_{i{j_4}}}} \right)} }.
 \end{equation*}
Since we assume that the columns of $\textbf{\emph{X}}$ are independent, the summation is not zero  when \({j_1} = {j_3}\)  and 
\({j_2} = {j_4}.\)  Thus, 
\begin{equation}\label{var_g_i}
\var\left( {{g_i}} \right) = \sum\limits_{{j_1} < {j_2}}^{} {E\left( {X_{i{j_1}}^2X_{i{j_2}}^2} \right) = \sum\limits_{{j_1} < {j_2}}^{} {E\left( {X_{i{j_1}}^2} \right)E\left( {X_{i{j_2}}^2} \right)} }  = p\left( {p - 1} \right)/2.
  \end{equation}
Notice that we used the assumption that since we assume that ${\bf{\Sigma}}={\bf{I}}$ in the last equality.
Now, recall by \eqref{var_naive_example} that
$\var\left( {{{\hat \tau }^2}} \right) = \frac{{20}}{n} + O\left( {\frac{1}{{{n^2}}}} \right).$
Therefore, we have
\begin{equation}\label{var_T_c_star_exmp1}
\var\left( {{T_{{c^*}}}} \right) = \var\left( {{{\hat \tau }^2}} \right) - \frac{{{{\left[ {2\sum\limits_{j = 1}^p {{\beta _j}E\left( {{S_{ij}}} \right)} } \right]}^2}}}{{n\var\left( {{g_i}} \right)}} = \frac{{20}}{n} + O\left( {\frac{1}{{{n^2}}}} \right) - \frac{{4{{\left( {p - 1} \right)}^2}}}{{n \cdot \left[ {p\left( {p - 1} \right)/2} \right]}} = \frac{{12}}{n} + O\left( {\frac{1}{{{n^2}}}} \right),
\end{equation}
where we used the assumption that $n=p$ in the last equality.

\end{remark}

\noindent\textbf{\textit{Proof of Proposition \ref{singel_asymptotic}}}:\\ 
We need to prove that
\(\sqrt n \left[ {{T_{{c^*}}} - {{ T}_{{\hat c^*}}}} \right]\overset{p}{\rightarrow} 0.\)
Write,
\[\sqrt n \left[ {{T_{c^*}} - T_{\hat c^*} } \right] = \sqrt n \left[ {{{\hat \tau }^2} - {c^*}{g_n} - \left( {{{\hat \tau }^2} - {{\hat c}^*}{g_n}} \right)} \right] =   {\sqrt n {g_n}}  {\left( {{{\hat c}^*} - {c^*}} \right)} .\]  
By Markov and Cauchy-Schwarz inequalities, it is enough to show that
$$p\left\{ {\left| {\sqrt n {g_n}\left( {{\hat c^*} - {{ c}^*}} \right)} \right| > \varepsilon } \right\} \le \frac{{E\left\{ {\left| {\sqrt n {g_n}\left( {{\hat c^*} - {{ c}^*}} \right)} \right|} \right\}}}{\varepsilon } \le \frac{{\sqrt {nE\left( {g_n^2} \right)E\left[ {{{\left( {{\hat c^*} - {{ c}^*}} \right)}^2}} \right]} }}{\varepsilon }\rightarrow 0. $$
Since $E(g_n^2)=\frac{\var(g_i)}{n}$ and $E[(\hat c^*- c^*)^2]=\var(\hat c^*),$ it enouth to show that $\var(g_i)\var(\hat c^*)\rightarrow 0.$
Notice that by \eqref{c_hat_star} we have
 \begin{equation}\label{U_statistic_C_star}
  \var \left( {{g_i}} \right)\var \left( {{{\hat c}^*}} \right) = \frac{{\var \left( U \right)}}{{\var \left( {{g_i}} \right)}}
 \end{equation}
 where $U \equiv \frac{2}{{n\left( {n - 1} \right)}}\sum\limits_{{i_1} \ne {i_2}}^{} {} \sum\limits_{j = 1}^p {{W_{{i_1}j}}  {S_{{i_2}j}}}$
  is a U-statistic  of order 2 with the kernel
 \(\left( {{{\bf{W}}_1},{{\bf{S}}_2}} \right) = {\bf{W}}_1^T{{\bf{S}}_2} = \sum\limits_{j = 1}^p {{W_{1j}}{S_{2j}}}. \)
 By Theorem 12.3 in \cite{van2000asymptotic},  the variance of $U$ 
 %a U-statistic of order 2 has a known formula. Specifically,
  \begin{equation}\label{var_U_statis}
 \var\left( U \right) = \frac{{4\left( {n - 2} \right)}}{{n\left( {n - 1} \right)}}{\delta _1} + \frac{2}{{n\left( {n - 1} \right)}}{\delta _2},    
 \end{equation}
  where \({\delta _1} = Cov\left[ {h\left( {{{\bf{W}}_1},{{\bf{S}}_2}} \right),h\left( {{{\bf{W}}_1},{{{\bf{\tilde S}}}_2}} \right)} \right]\) and \({\delta _2} = Cov\left[ {h\left( {{{\bf{W}}_1},{{\bf{S}}_2}} \right),h\left( {{{\bf{W}}_1},{{\bf{S}}_2}} \right)} \right]\).
  Consider now the  denominator of \eqref{U_statistic_C_star}. Write,
 \begin{equation*}
 \var\left( {{g_i}} \right) = E\left( {g_i^2} \right) = E\left[ {{{\left( {\sum\limits_{j < j'}^{} {{X_{ij}}{X_{ij'}}} } \right)}^2}} \right] = \sum\limits_{{j_1} < {j_2}}^{} {\sum\limits_{{j_3} < {j_4}}^{} {E\left( {{X_{i{j_1}}}{X_{i{j_2}}}{X_{i{j_3}}}{X_{i{j_4}}}} \right)} }.
 \end{equation*}
Since we assume that the columns of $\textbf{\emph{X}}$ are independent, the summation is not zero  when \({j_1} = {j_3}\)  and 
\({j_2} = {j_4}.\)  Thus, 
\begin{equation}
\var\left( {{g_i}} \right) = \sum\limits_{{j_1} < {j_2}}^{} {E\left( {X_{i{j_1}}^2X_{i{j_2}}^2} \right) = \sum\limits_{{j_1} < {j_2}}^{} {E\left( {X_{i{j_1}}^2} \right)E\left( {X_{i{j_2}}^2} \right)} }  = p\left( {p - 1} \right)/2.
  \end{equation}
Notice that since we assume that ${\bf{\Sigma}}={\bf{I}}$ then $E(X_{ij}^2)=1$ for all $i=1,...,n$ and $j=1,...,p$ .
 Now, since we assume that $n/p=O(1),$  
  by \eqref{U_statistic_C_star} and  \eqref{var_U_statis} it is enough to prove that $\frac{{{\delta _1}}}{{n^3 }} \to 0$ and 
$\frac{{{\delta _2}}}{{{n^4}}} \to 0.$

Consider first $\frac{\delta_1}{n^3}.$ Write,
\begin{multline*}
    {\delta _1} = Cov\left[ {h\left( {{{\bf{W}}_1},{{\bf{S}}_2}} \right),h\left( {{{\bf{W}}_1},{{{\bf{\tilde S}}}_2}} \right)} \right] = Cov\left[ {\sum\limits_{j = 1}^p {{W_{1j}}{S_{2j}}} ,\sum\limits_{j = 1}^p {{W_{1j}}{{\tilde S}_{2j}}} } \right]  \\
    = \sum\limits_{j,j'}^{} {\left\{ {E\left( {{W_{1j}}{W_{1j'}}} \right){\theta _j}{\theta _{j'}} - {\beta _j}{\beta _{j'}}{\theta _j}{\theta _{j'}}} \right\}}  =
    {\theta ^T}{\bf{A}}\theta  - {\left( {\sum\limits_{j = 1}^p {{\theta _j}{\beta _j}} } \right)^2}, 
\end{multline*}
   where \({\theta _j} \equiv E\left( {{S_{ij}}} \right)\),  \(\theta  = {\left( {{\theta _1},...,{\theta _p}} \right)^T}\) , \({\bf{A}} = E\left( {{\bf{W}}{{\bf{W}}^T}} \right)\) and \({\bf{W}} = \left( {{W_{i1}},...,{W_{ip}}} \right)\). Thus, we need to show that \(\frac{{{\theta ^T}{\bf{A}}\theta }}{{n^3}} \to 0.\)
   Let $\lambda_1 \leq\lambda_2 \leq ... \leq\lambda_p$ be the eigenvalues of {\bf{A}}. Notice that \(\lambda _1^2 \le \sum\limits_{j = 1}^p {\lambda _j^2 = } trace\left( {{{\bf{A}}^2}} \right) = \left\| {\bf{A}} \right\|_F^2\), and since by Corollary \ref{consistency_naive}  we have   \(\frac{{\left\| {\bf{A}} \right\|_F^2}}{{{n^2}}}  \to 0\) then also $\frac{\lambda_1}{n}\rightarrow 0.$
      Now, since $\max_{\theta} \left( {\frac{{{\theta ^T}{\bf{A}}\theta }}{{{{\left\| \theta  \right\|}^2}}}} \right) = {\lambda_1}$ then
$\frac{{{\theta ^T}{\bf{A}}\theta }}{{n^3}} = \frac{1}{{n^3}}{\left\| \theta  \right\|^2} \cdot \left( {\frac{{{\theta ^T}{\bf{A}}\theta }}{{{{\left\| \theta  \right\|}^2}}}} \right) \le \frac{1}{{n^3}}{\left\| \theta  \right\|^2} \cdot {\lambda _1}.$ Thus, it is enough to show that \(\frac{{{{\left\| \theta  \right\|}^2}}}{{{n^2}}}\) is bounded.
% Now, by Cauchy–Schwarzand and \eqref{var_g_i},
% \[\theta _j^2 = {\left[ {E\left( {W_{ij}^{}{g_i}} \right)} \right]^2} \le E\left( {W_{ij}^2} \right)E\left( {g_i^2} \right) = E\left( {X_{ij}^2Y_i^2} \right){\mathop{\rm var}} \left( {{g_i}} \right) \le {\left[ {E\left( {X_{ij}^4} \right)E\left( {Y_i^4} \right)} \right]^{1/2}}{p^2}\]
Write,
\[\begin{array}{l}{\theta _j} \equiv E\left( {{S_{ij}}} \right) = E\left( {{X_{ij}}{Y_i}{g_i}} \right) = E\left( {{X_{ij}}\left( {{\beta ^T}{X_i} + \varepsilon_i } \right){g_i}} \right) = \\E\left( {{X_{ij}}\left( {\sum\limits_{m = 1}^p {{\beta _m}{X_{im}}}  + \varepsilon_i } \right)\sum\limits_{k < k'}^{} {{X_{ik}}{X_{ik'}}} } \right) = \sum\limits_{m = 1}^p {} \sum\limits_{k < k'}^{} {{\beta _m}E\left( {{X_{ij}}{X_{im}}{X_{ik}}{X_{ik'}}} \right)} \end{array}\]
where in the last equality we used the assumption that $E(\epsilon|X)=0$.
Since we assume that the columns of $\bf{X}$ are independent,  the summation is not zero (up to permutations) when \(j = k\) and \(m = k'\). In this case we have by Cauchy–Schwartz 
\[{\theta _j} = \sum\limits_{m = 1}^p {\sum\limits_{k < k'}^{} {{\beta _m}E\left( {{X_{ij}}{X_{im}}{X_{ik}}{X_{ik'}}} \right)} }  = \sum\limits_{m \neq j}^p {{\beta _m}E\left( {X_{ij}^2X_{im}^2} \right)}  = \sum\limits_{m \neq j}^p {{\beta _m}E\left( {X_{ij}^2} \right)E\left( {X_{im}^2} \right)}  = \sum\limits_{m \neq j}^p {{\beta _m}}  \le \sqrt p \left\| \beta  \right\|.\]
Notice that in the forth equality we used the assumption that $E(X_{ij}^2)=1$ for all $i=1,...,n$ and $j=1,...,p.$ 
Thus \({\left\| \theta  \right\|^2} = \sum\limits_{j = 1}^p {\theta _j^2}  \le \sum\limits_{j = 1}^p {\left( {p{{\left\| \beta  \right\|}^2}} \right)}  = {p^2}{\tau ^2} = O\left( {{p^2}} \right).\)
Since we assume that $p/n=O(1),$ then \(\frac{{{{\left\| \theta  \right\|}^2}}}{{{n^2}}}\) is indeed bounded. 

Consider now $\frac{{{\delta _2}}}{{{n^4}}}.$
Write,
\begin{multline*}
    {\delta _2} = Cov\left[ {h\left( {{{\bf{W}}_1},{{\bf{S}}_2}} \right),h\left( {{{\bf{W}}_1},{{\bf{S}}_2}} \right)} \right] = Cov\left[ {\sum\limits_{j = 1}^p {{W_{1j}}{S_{2j}}} ,\sum\limits_{j = 1}^p {{W_{1j}}{S_{2j}}} } \right]  = \sum\limits_{j,j'}^{} {\left\{ {E\left( {{W_{1j}}{W_{1j'}}} \right)E\left( {{S_{2j}}{S_{2j'}}} \right) - {\beta _j}{\beta _{j'}}{\theta _j}{\theta _{j'}}} \right\}} 
\end{multline*}
 Since $p/n= O(1)$ by assumption, it is enough to show that 
 $\eta\equiv \sum\limits_{j,j'}^{} {E\left( {{W_{1j}}{W_{1j'}}} \right)E\left( {{S_{2j}}{S_{2j'}}} \right)  } = O(p^3).$
Write,
\begin{align*}
&\sum\limits_{j,j'}^{} {E\left( {{W_{1j}}{W_{1j'}}} \right)E\left( {{S_{2j}}{S_{2j'}}} \right) = } \sum\limits_{j,j'}^{} {E\left( {{W_{1j}}{W_{1j'}}} \right)E\left( {{W_{2j}}{W_{2j'}}g_2^2} \right)} 
\\& = \sum\limits_{j,j'}^{} {E\left( {{X_{1j}}{X_{1j'}}Y_1^2} \right)E\left[ {{X_{2j}}{X_{2j'}}Y_2^2{{\left( {\sum\limits_{k < k'}^{} {{X_{2k}}{X_{2k'}}} } \right)}^2}} \right]} 
\\ &= \sum\limits_{j,j'}^{} {E\left( {{X_{1j}}{X_{1j'}}Y_1^2} \right)E\left[ {{X_{2j}}{X_{2j'}}Y_2^2\sum\limits_{{k_1} < {k_2}}^{} {\sum\limits_{{k_3} < {k_4}}^{} {{X_{2{k_1}}}{X_{2{k_2}}}{X_{2{k_3}}}{X_{2{k_4}}}} } } \right]}
\\ &= \sum\limits_{{j_1},{j_2}}^{} {\sum\limits_{{k_1} < {k_2}}^{} {\sum\limits_{{k_3} < {k_4}}^{} {E\left( {{X_{1{j_1}}}{X_{1{j_2}}}Y_1^2} \right)} } E\left( {{X_{2{j_1}}}{X_{2{j_2}}}Y_2^2{X_{2{k_1}}}{X_{2{k_2}}}{X_{2{k_3}}}{X_{2{k_4}}}} \right)}. 
    \end{align*}
Plug-in $Y=\beta^TX+\epsilon$ to get
\begin{align*}
 \eta &= \sum\limits_{{j_1},{j_2}}^{} {\sum\limits_{{k_1} < {k_2}}^{} {\sum\limits_{{k_3} < {k_4}}^{} {E\left( {{X_{1{j_1}}}{X_{1{j_2}}}Y_1^2} \right)} } E\left( {{X_{2{j_1}}}{X_{2{j_2}}}Y_2^2{X_{2{k_1}}}{X_{2{k_2}}}{X_{2{k_3}}}{X_{2{k_4}}}} \right)}
\\ &= \sum\limits_{{j_1},{j_2}}^{} {\sum\limits_{{k_1} < {k_2}}^{} {\sum\limits_{{k_3} < {k_4}}^{} {E\left[ {{X_{1{j_1}}}{X_{1{j_2}}}{{\left( {{\beta ^T}{X_1} + {\varepsilon _1}} \right)}^2}} \right]} } E\left[ {{X_{2{j_1}}}{X_{2{j_2}}}{{\left( {{\beta ^T}{X_2} + {\varepsilon _2}} \right)}^2}{X_{2{k_1}}}{X_{2{k_2}}}{X_{2{k_3}}}{X_{2{k_4}}}} \right]} 
\\ &= \sum\limits_{{j_1},{j_2}}^{} {\sum\limits_{{k_1} < {k_2}}^{} {\sum\limits_{{k_3} < {k_4}}^{} {E\left[ {{X_{1{j_1}}}{X_{1{j_2}}}\left( {\sum\limits_{j,j'}^{} {{\beta _j}{\beta _{j'}}{X_{1j}}{X_{1j'}} + \varepsilon _1^2} } \right)} \right]} } E\left[ {{X_{2{j_1}}}{X_{2{j_2}}}\left( {\sum\limits_{j,j'}^{} {{\beta _j}{\beta _{j'}}{X_{2j}}{X_{2j'}}}  + \varepsilon _2^2} \right){X_{2{k_1}}}{X_{2{k_2}}}{X_{2{k_3}}}{X_{2{k_4}}}} \right]}
\\ &= \sum\limits_{{j_1},{j_2}}^{} {\sum\limits_{{j_3},{j_4}}^{} {\sum\limits_{{j_5},{j_6}}^{} {} } \sum\limits_{{k_1} < {k_2}}^{} {\sum\limits_{{k_3} < {k_4}}^{} {{\beta _{{j_3}}}{\beta _{{j_4}}}{\beta _{{j_5}}}{\beta _{{j_6}}}} } \overbrace {E\left( {{X_{1{j_1}}}{X_{1{j_2}}}{X_{1{j_3}}}{X_{1{j_4}}}} \right)}^{{L_1}}\overbrace {E\left( {{X_{2{j_1}}}{X_{2{j_2}}}{X_{2{j_5}}}{X_{2{j_6}}}{X_{2{k_1}}}{X_{2{k_2}}}{X_{2{k_3}}}{X_{2{k_4}}}} \right)}^{L{  _2}}} 
\\ &+ {\sigma ^4}\sum\limits_{{j_1},{j_2}}^{} {\sum\limits_{{k_1} < {k_2}}^{} {\sum\limits_{{k_3} < {k_4}}^{} {} } \overbrace {E\left( {{X_{1{j_1}}}{X_{1{j_2}}}} \right)}^{{L_3}}\overbrace {E\left( {{X_{2{j_1}}}{X_{2{j_2}}}{X_{2{k_1}}}{X_{2{k_2}}}{X_{2{k_3}}}{X_{2{k_4}}}} \right)}^{{L_4}}},
\end{align*}
where in the forth equality we used the assumption that $E(\epsilon^2|x)=\sigma^2$. Now, notice that:
\begin{itemize}
    \item  $L_1$ is not zero (up to permutation) when $j_1=j_2$ and $j_3=j_4.$
    \item $L_2$ is not zero (up to permutation) when $j_1=j_2; j_5=j_6; k_1=k_3; k_2=k_4.$
    \item $L_3$ is not zero when $j_1=j_2$
    \item $L_4$ is not zero (up to permutation) when $j_1=j_2;  k_1=k_3; k_2=k_4.$
\end{itemize}
 Putting it all together to get
 \[\eta =   \sum\limits_{{j_1}}^{} {\sum\limits_{{j_3}}^{} {\sum\limits_{{j_5}}^{} {\sum\limits_{{k_1} < {k_2}}^{} {\beta _{{j_3}}^2\beta _{{j_5}}^2} E\left( {X_{i{j_1}}^2X_{1{j_3}}^2} \right)E\left( {X_{2{j_1}}^2X_{2{j_5}}^2X_{2{k_1}}^2X_{2{k_2}}^2} \right)} } } \,\, + {\sigma ^4}\sum\limits_{{j_1}}^{} {\sum\limits_{{k_1} < {k_2}}^{} { {} E\left( {X_{1{j_1}}^2} \right)} E\left( {X_{2{j_1}}^2X_{2{k_1}}^2X_{2{k_2}}^2} \right)}.  \]
Since we assume that $E\left( {X_{1{j_1}}^2X_{1{j_2}}^2X_{1{j_3}}^2 X_{1{j_4}}^2} \right)=O(1)$ for every $j_1,j_2,j_3,j_4$, all the expectations above are also $O(1)$.  Also recall that $\tau^2+ \sigma^2 = O(1)$, and $\sum\limits_{{j_3}}^{} {\sum\limits_{{j_5}}^{} {\beta _{{j_3}}^2\beta _{{j_5}}^2} } \, = {\tau ^4} = O\left( 1 \right).$
Thus, we obtain that $\eta= O(p^3).$  This completes the   proof that \(\sqrt n \left[ {{T_{{c^*}}} - {{ T}_{{\hat c^*}}}} \right]\overset{p}{\rightarrow} 0.\)

\begin{remark}\label{rem:improve_T_B}
We now calculate the asymptotic improvement of $T_{\bf{B}}$ over the naive estimator. For simplicity, consider the case when $\tau^2=\sigma^2=1.$
Recall the variance of $\hat{\tau}^2$ and $T_{\bf{B}}$ given in \eqref{var_naive} and \eqref{var_T_B}, respectively. Write,
\begin{align*}
\mathop {\lim }\limits_{n,p \to \infty } \frac{{\var\left( {{{\hat \tau }^2}} \right) - \var\left( {{T_{\bf{B}}}} \right)}}{{\var\left( {{{\hat \tau }^2}} \right)}}
&= \mathop {\lim }\limits_{n,p \to \infty }\frac{{8{\tau_{\bf{B}} ^4}/n}}{{\frac{4}{n}\left[ {\frac{{\left( {n - 2} \right)}}{{\left( {n - 1} \right)}}\left( {\sigma_Y^2{\tau ^2} + {\tau ^4}} \right) + \frac{1}{{2\left( {n - 1} \right)}}\left( {p{\sigma_Y^4} + 4\sigma_Y^2{\tau ^2} + 3{\tau ^4}} \right)} \right]}}\\ 
&= \frac{{2{\tau_{\bf{B}}^4}}}{{3{\tau ^4} + \frac{{4p{\tau ^4} + 4\sigma_Y^2{\tau ^2} + 3{\tau ^4}}}{{2n}}}} = \frac{0.5}{{3 + 2\frac{p}{n}}},
 \end{align*}
 where we used \eqref{var_naive} in the first equality, and the fact that $\sigma_Y^2= 2\tau^2=2$ in the second equality. 
Now, notice that  when $p=n$  and $\tau_{\bf{B}}^2=0.5$ then the reduction is $\frac{0.5}{{3 + 2}} = 10\%$ and when $p/n$ converges to zero, the reduction is $16\%.$
\end{remark}

\begin{remark}\label{summary_eqample}
\textbf{\textit{Calculations for Example \ref{extreme_example}}}:\\ 
Consider the first scenario where $\beta_j^2=\frac{1}{p}.$
Recall that we assume that the set $\textbf{B}$ is a fixed set of indices such that \(\left| {\bf{B}} \right| \ll p\). Therefore, we have \(\tau _{\bf{B}}^2 = \sum\limits_{j \in {\bf{B}}}^{} {\beta _j^2 = O\left( {\frac{1}{p}} \right)} \).
Now, by \eqref{var_T_B_norm} we have $\var(T_{\textbf{B}})=\var(\hat\tau^2)-\frac{8}{n}\tau^2_{\textbf{B}}+O(n^{-2})$ and by Remark \ref{rem:improve_singel} we have $\var(\hat\tau^2) = \frac{20}{n}+O(n^{-2}).$
Using the assumption that $n=p$ we can conclude that \(\var\left( {{T_{\bf{B}}}} \right) = \frac{{20}}{n} + O\left( {\frac{1}{{{n^2}}}} \right).\)
Hence, in this scenario, $T_{\textbf{B}}$ and the naive estimator have the same asymptotic variance. In contrast, recall that in Example \ref{exp_singel} we showed that the asymptotic variance of $T_{c^*}$ is $40\%$ lower than the variance of the naive estimator.

Consider now the second scenario where \(\hat \tau _{\bf{B}}^2 = {\tau ^2} = 1\). By \eqref{var_T_B_norm} we have
\begin{equation*}
\var\left( {{T_{\bf{B}}}} \right) =  {\var\left( {{{\hat \tau }^2}} \right) - \frac{8}{n}\tau _{\bf{B}}^4}  +O(n^{-2})=\frac{12}{n}+ O(n^{-2}).
\end{equation*}
Hence, in this scenario the asymptotic variance of $T_{\textbf{B}}$ is $40\%$ smaller than the variance of the naive estimator.
Consider now $\var(T_{c^*}).$ By Cauchy–Schwarz inequality
\({\left( {\sum\limits_{j \in {\bf{B}}}^{} {{\beta _j}} } \right)^2} \le \sum\limits_{j \in {\bf{B}}}^{} {\beta _j^2 \cdot \left| {\bf{B}} \right|}  = \tau^2_{\textbf{B}}\left| {\bf{B}} \right| = O(1)\), where the last equality holds since we assume that \({\bf{B}} \subset \left\{ {1,...,p} \right\}\) 
   be a fixed set of some indices such that
    \(\left| {\bf{B}} \right| \ll p.\) Now, 
By \eqref{numerator_single} we have \[\sum\limits_{j = 1}^p {{\beta _j}} {\theta _j} = {\left( {\sum\limits_{j = 1}^p {{\beta _j}} } \right)^2} - {\tau ^2} = {\left( {\sum\limits_{j \in {\bf{B}}}^p {{\beta _j}}  + \overbrace {\sum\limits_{j \notin {\bf{B}}}^p {{\beta _j}} }^0} \right)^2} - \tau _{\bf{B}}^2 \le \left| {\bf{B}} \right| - 1 = O\left( 1 \right).\]
Now, recall  that
$\var\left( {{{\hat \tau }^2}} \right) = \frac{{20}}{n} + O\left( {\frac{1}{{{n^2}}}} \right)$ and $\var(g_i)=p(p-1)/2$ by \eqref{var_naive_example} and \eqref{var_g_i} respectively.
Therefore, we have
\begin{equation*}
\var\left( {{T_{{c^*}}}} \right) = \var\left( {{{\hat \tau }^2}} \right) - \frac{{{{\left[ {2\sum\limits_{j = 1}^p {{\beta _j}\theta_j} } \right]}^2}}}{{n\var\left( {{g_i}} \right)}} = \frac{{20}}{n} + O\left( {\frac{1}{{{n^2}}}} \right) -O\left(\frac{1}{np^2}\right) = \frac{{20}}{n} + O\left( {\frac{1}{{{n^2}}}} \right),
\end{equation*}
Hence, in this scenario, $T_{c^*}$ and the naive estimator have the same asymptotic variance. 

Lastly, recall that in Example \ref{exp_OOE} we already showed that, asymptotically, the variance of $T_{oracle}$ (i.e., the optimal oracle estimator) is $40\%$ lower than the naive variance (without any assumptions about the structure of the coefficient vector~$\beta$).

\end{remark}

\noindent\textbf{\textit{Proof of Proposition \ref{limit_of_proposed}}}:\\ 
In order to prove that \(\sqrt n \left( {{T_\gamma } - {T_{\bf{B}}}} \right) \overset{p}{\rightarrow} 0,\) it is enough to show that
 \begin{align}
  &E\left\{ {\sqrt n \left( {{T_\gamma } - {T_{\bf{B}}}} \right)} \right\} \to 0 , \label{b1}\\
   &\var\left\{ {\sqrt n \left( {{T_\gamma } - {T_{\bf{B}}}} \right)} \right\} \to 0.   \label{b2}
 \end{align}
%  We need to prove that
%  \begin{align}
%   &\mathop {\lim }\limits_{n \to \infty } n\left[  {{ {E\left( {{T_\gamma }} \right) - {\tau ^2}} }} \right]^2 = 0 , \label{b1}\\
%   &\mathop {\lim }\limits_{n \to \infty } n\left[ {\var\left( {{T_{\bf{B}}}} \right) - {\var} \left( {{T_\gamma }} \right)} \right] = 0.   \label{b2}
%  \end{align}
 We start with the first equation.
  Let $A$ denote the event that the selection algorithm $\gamma$ perfectly identifies the  set of large coefficients, i.e., \(A = \left\{ {{{\bf{B}}_\gamma } = {\bf{B}}} \right\}.\) Let $p_A\equiv P(A)$ denote the probability that $A$ occurs, and let $\mathbb{1}_A$ denote the indicator of   $A.$ Notice that $E(T_{\textbf{B}})=\tau^2$ and ${T_\gamma }\mathbb{1}_A = {T_{\rm B}}\mathbb{1}_A.$ Thus,
\begin{align}
  E\left\{ {\sqrt n \left( {{T_\gamma } - {T_{\bf{B}}}} \right)} \right\} =  \sqrt{n}\left[ {E\left( {{T_\gamma }} \right) - {\tau ^2}} \right] &=\sqrt{n}  \left(  E\left[ {{T_\gamma }(1-\mathbb{1}_A)} \right]  +E\left[ {{T_\gamma }\mathbb{1}_A} \right]- {\tau ^2} \right)\nonumber\\ 
&=  
 \sqrt{n}E\left[ {{T_\gamma }(1-\mathbb{1}_A)} \right]+
\sqrt{n}\left[ E\left( {{T_{\bf{B}}}\mathbb{1}_A} \right)- {\tau ^2}\right] , \label{dsg}  
     \end{align}
where the last equality holds since  ${T_\gamma }\mathbb{1}_A = {T_{\rm B}}\mathbb{1}_A.$
For the convenience of notation, let $C$ be an upper bound of the maximum over all  first four moments of $ {T_\gamma} $ and $ T_{\bf{B}},$ and consider  the first term of \eqref{dsg}. By the Cauchy–Schwarz inequality,
\begin{equation}\label{TT12}
\sqrt{n}E\left[ {{T_\gamma }\left( {1 - \mathbb{1}_A)} \right)} \right] \le \sqrt{n}{\left\{ {E\left[ {T_\gamma ^2} \right]} \right\}^{1/2}}{\left\{ {E\left[ {{{\left( {1 - \mathbb{1}_A} \right)}^2}} \right]} \right\}^{1/2}} \le \sqrt{n}{C^{1/2}}{\left\{ {1 - {p_A}} \right\}^{1/2}} \underset{n \to \infty}{\rightarrow} 0,    
\end{equation}
where the last inequality holds since  \(\mathop {\lim }\limits_{n \to \infty } n\left( {1 - {p_A}} \right)^{1/2} = 0\) by assumption. 
We now consider the second term of \eqref{dsg}. Write, \begin{equation} %\label{yre}
\sqrt{n}\left[ {E\left( {{T_{\bf{B}}}\mathbb{1}_A} \right) - {\tau ^2}} \right] =  \sqrt{n}E\left( {{T_{\bf{B}}}\mathbb{1}_A - {T_{\bf{B}}}} \right) =  
-\sqrt{n}E\left[ {{T_{\bf{B}}}\left( {1 - \mathbb{1}_A} \right)} \right],    
\end{equation}
 and notice that
by the same type of argument as in \eqref{TT12} we have
$\sqrt{n}E\left[ {{T_{\bf{B}}}\left( {1 - \mathbb{1}_A} \right)} \right] \underset{n \to \infty}{\rightarrow}~0.$
% The only subtle difference between \eqref{yre} and \eqref{TT12} is that  \(E\left( {T_{\bf{B}}^2} \right)\) is finite by construction, i.e.,
% \(E\left( {T_{\bf{B}}^2} \right) = {\mathop{\var}} \left( {{T_{\bf{B}}}} \right) + {\left[ {E\left( {{T_{\bf{B}}}} \right)} \right]^2} = O\left( {{n^{ - 1}}} \right) + {\tau ^4} \le C_2\) for some constant $C_2.$
This completes the proof of \eqref{b1}.

We now move to show that
$\var\left\{ {\sqrt n \left( {{T_\gamma } - {T_{\bf{B}}}} \right)} \right\} \to 0$. Write,
\begin{align*}
\var\left\{ {\sqrt n \left( {{T_\gamma } - {T_{\bf{B}}}} \right)} \right\}
&= n\var\left( {{T_\gamma } - {T_{\bf{B}}}} \right)\\
&= n\left[ {\var\left( {{T_\gamma }} \right) + \var\left( {{T_{\bf{B}}}} \right) - 2{\cov} \left( {{T_\gamma },{T_{\bf{B}}}} \right)} \right]\\
&= n\left\{ {E\left( {T_\gamma ^2} \right) - {{\left[ {E\left( {{T_\gamma }} \right)} \right]}^2} + E\left( {T_{\bf{B}}^2} \right) - {\tau ^4} - 2\left[ {E\left( {{T_\gamma }{T_{\bf{B}}}} \right) - E\left( {{T_\gamma }} \right){\tau ^2}} \right]} \right\}\\
&= n\left\{ {E\left( {T_\gamma ^2} \right) - E\left( {{T_\gamma }{T_{\bf{B}}}} \right) + E\left( {T_{\bf{B}}^2} \right) - E\left( {{T_\gamma }{T_{\bf{B}}}} \right) + E\left( {{T_\gamma }} \right)\left[ {{\tau ^2} - E\left( {{T_\gamma }} \right)} \right] - {\tau ^2}\left[ {{\tau ^2} - E\left( {{T_\gamma }} \right)} \right]} \right\}\\ 
&= n\left\{ {\underbrace {E\left( {T_\gamma ^2} \right) - E\left( {{T_\gamma }{T_{\bf{B}}}} \right)}_{{\theta _1}} + \underbrace {E\left( {T_{\bf{B}}^2} \right) - E\left( {{T_\gamma }{T_{\bf{B}}}} \right)}_{{\theta _2}} - \underbrace {{{\left[ {{\tau ^2} - E\left( {{T_\gamma }} \right)} \right]}^2}}_{{\theta _3}}} \right\}.
\end{align*}
Thus, it is enough to show that $n\theta_1\rightarrow 0$, $n\theta_2\rightarrow 0$ and $n\theta_3\rightarrow 0.$

We start with showing that $n\theta_1\rightarrow 0.$ Notice that $T_{\textbf{B}}^2\mathbb{1}_A=T_{\textbf{B}}T_{\gamma}\mathbb{1}_A=T_{\gamma}^2\mathbb{1}_A$
Thus,
\begin{align*}
  n{\theta _1} 
  &= n\left\{ {E\left( {T_\gamma ^2} \right) - E\left( {{T_\gamma }{T_{\bf{B}}}} \right)} \right\}\\ 
  &= n\left\{ {E\left( {T_\gamma ^2} \right) - E\left[ {{T_\gamma }{T_{\bf{B}}}\left( {1 - {\mathbb{1}_A}} \right)} \right] - E\left( {{T_\gamma }{T_{\bf{B}}}{\mathbb{1}_A}} \right)} \right\}  \\
  &= n\left\{ {E\left( {T_\gamma ^2} \right) - E\left[ {{T_\gamma }{T_{\bf{B}}}\left( {1 - {\mathbb{1}_A}} \right)} \right] - E\left( {T_\gamma ^2{\mathbb{1}_A}} \right)} \right\}  \\
  &= n\left\{ {E\left[ {T_\gamma ^2\left( {1 - {\mathbb{1}_A}} \right)} \right] - E\left[ {{T_\gamma }{T_{\bf{B}}}\left( {1 - {\mathbb{1}_A}} \right)} \right]} \right\}.
\end{align*}
Now, notice that  $n(E\left[ {T_\gamma ^2\left( {1 - {\mathbb{1}_A}} \right)} \right])\rightarrow 0$ by  similar arguments as in \eqref{TT12}, with a slight modification of using the existence of the fourth moments of $T_{\gamma}$ and $T_{\bf{B}}$, rather than the second moments. Also, by Cauchy–Schwarz inequality  we have,
\begin{align*}
nE\left[ {{T_\gamma }{T_{\bf{B}}}\left( {1 - {1_A}} \right)} \right] 
&\le n{\left\{ {E\left( {T_\gamma ^2T_{\bf{B}}^2} \right)} \right\}^{1/2}}{\left\{ {E\left[ {{{\left( {1 - {1_A}} \right)}^2}} \right]} \right\}^{1/2}}\\
&\le n{\left\{ {E\left( {T_\gamma ^4} \right)E\left( {T_{\bf{B}}^4} \right)} \right\}^{1/4}}{\left\{ {1 - {p_A}} \right\}^{1/2}}\\
&\le n{C^{1/2}}{\left\{ {1 - {p_A}} \right\}^{1/2}} \to 0,
\end{align*}
where $C$ is an upper bound of the maximum over all  first four moments of $ {T_\gamma} $ and $ T_{\bf{B}}.$
Therefore, $n\theta_1 \rightarrow 0.$

Consider now $n\theta_2.$ Write,
\begin{align*}
n{\theta _2} &= n\left\{ {E\left( {T_{\bf{B}}^2} \right) - E\left( {{T_\gamma }{T_{\bf{B}}}} \right)} \right\} \\
&= n\left\{ {E\left( {T_{\bf{B}}^2} \right) - E\left[ {{T_\gamma }{T_{\bf{B}}}\left( {1 - {\mathbb{1}_A}} \right)} \right] - E\left( {{T_\gamma }{T_{\bf{B}}}{\mathbb{1}_A}} \right)} \right\}  \\
&= n\left\{ {E\left( {T_{\bf{B}}^2} \right) - E\left[ {{T_\gamma }{T_{\bf{B}}}\left( {1 - {\mathbb{1}_A}} \right)} \right] - E\left( {T_{\bf{B}}^2{\mathbb{1}_A}} \right)} \right\} \\
&= n\left\{ {E\left[ {T_{\bf{B}}^2\left( {1 - {\mathbb{1}_A}} \right)} \right] - E\left[ {{T_\gamma }{T_{\bf{B}}}\left( {1 - {\mathbb{1}_A}} \right)} \right]} \right\}\rightarrow 0, 
\end{align*}
and notice that the last equation follows by similar arguments.  

Consider now $n\theta_3.$ Write,
\begin{align*}
n\theta_3 &= n\left[ {E\left( {{T_{\bf{B}}}} \right) - E\left( {{T_\gamma }} \right)} \right]\\ 
&= n\left[ {E\left[ {{T_{\bf{B}}}\left( {1 - {\mathbb{1}_A}} \right) + {T_{\bf{B}}}{\mathbb{1}_A}} \right] - E\left( {{T_\gamma }} \right)} \right]\\
&= n\left\{ {E\left[ {{T_{\bf{B}}}\left( {1 - {\mathbb{1}_A}} \right)} \right] + E\left( {{T_{\bf{B}}}{\mathbb{1}_A} - {T_\gamma }} \right)} \right\}\\
&= n\left\{ {E\left[ {{T_{\bf{B}}}\left( {1 - {\mathbb{1}_A}} \right)} \right] + E\left( {{T_\gamma }{\mathbb{1}_A} - {T_\gamma }} \right)} \right\}\\
&= n\left\{ {E\left[ {{T_{\bf{B}}}\left( {1 - {\mathbb{1}_A}} \right)} \right] - E\left[ {{T_\gamma }\left( {1 - {\mathbb{1}_A}} \right)} \right]} \right\}\rightarrow 0,
 \end{align*}
where the last equation follows by  similar arguments as in \eqref{TT12}. This completes the proof of \eqref{b2} and we conclude that \(\sqrt n \left( {{T_\gamma } - {T_{\bf{B}}}} \right) \overset{p}{\rightarrow} 0.\) \qed

  \noindent\textbf{\textit{Proof of Proposition~ \ref{var_naive_est}}}:\\ 
We wish to prove that
\begin{align}\label{cons_tau}
  n\left[ {\widehat {\var\left( {{{\hat \tau }^2}} \right)} - \var\left( {{\hat\tau ^2}} \right)} \right] \overset{p}{\rightarrow} 0. 
 \end{align}  
 Recall by \eqref{var_naive} that
\begin{equation*}
{\var} \left( {{{\hat \tau }^2}} \right) = \frac{{4\left( {n - 2} \right)}}{{n\left( {n - 1} \right)}}\left[ {{\beta ^T}{\bf{A}}\beta  - {{\left\| \beta  \right\|}^4}} \right] + \frac{2}{{n\left( {n - 1} \right)}}\left[ {\left\| {\bf{A}} \right\|_F^2 - {{\left\| \beta  \right\|}^4}} \right].
\end{equation*}
Now,  when we assume  standard Gaussian covariates, one can verify that
$  {\beta ^T}{\bf{A}}\beta  - {\left\| \beta  \right\|^4} = \sigma_Y^2{\tau ^2} + {\tau ^4}$
and
$ \left\| {\bf{A}} \right\|_F^2 - {\left\| \beta  \right\|^4} = p{\sigma_Y^4} + 4\sigma_Y^2{\tau ^2} + 3{\tau ^4},$ where $\sigma_Y^2=\sigma^2+\tau^2.$  Thus, in this case  we can write
\begin{equation}
\var\left( {{{\hat \tau }^2}} \right) = \frac{4}{n}\left[ {\frac{{\left( {n - 2} \right)}}{{\left( {n - 1} \right)}}\left[ {\sigma_Y^2{\tau ^2} + {\tau ^4}} \right] + \frac{1}{{2\left( {n - 1} \right)}}\left( {p{\sigma_Y^4} + 4\sigma_Y^2{\tau ^2} + 3{\tau ^4}} \right)} \right].    
\end{equation}
  In order to prove that \eqref{cons_tau} holds, it is enough to prove the consistency of $\hat\tau^2$ and $\hat\sigma_Y^2.$
  Consistency of the sample variance~$\hat\sigma_Y^2$ is a standard result,   and since  $\hat\tau^2$ is an unbiased estimator, it is enough to show that its variance converges to zero as $n \to \infty.$ 
  Since we assume $\hat \tau^2 +\sigma^2 = O(1)$ and $p/n = O(1),$  we have by~\eqref{var_naive_normal} that   $\var(\hat\tau^2)\underset{n \to \infty}{\rightarrow}0$,  
   and \eqref{cons_tau} follows.

  \noindent\textbf{\textit{Proof of Proposition~ \ref{consist_var_}}}:\\ 
We now move to prove that  
\begin{equation}\label{sh2}
n\left[ {\widehat {\var\left( {T_{\gamma}} \right)} - \var\left( {{T_{\gamma}}} \right)} \right] \overset{p}{\rightarrow}0, \end{equation}
   Recall that by Proposition \ref{limit_of_proposed} we have
   $\mathop {\lim }\limits_{n \to \infty } n\left[ {\var\left( {{T_{\bf{B}}}} \right) - {\var} \left( {{T_\gamma }} \right)} \right] = 0.$
Hence, it is enough to show that
$$n\left[ {\widehat {\var\left( {T_{\gamma}} \right)} - \var\left( T_{\bf{B}} \right)} \right] \overset{p}{\rightarrow}~0.$$
Since we assume \({X_i}\mathop \sim\limits^{i.i.d} N\left( {\bf{0} ,\bf{I} } \right)\) then by \eqref{var_T_B} we have
   $\var\left( {{T_{\bf{B}}}} \right) =  {\var\left( {{{\hat \tau }^2}} \right) - \frac{8}{n}\tau _{\bf{B}}^4}  +O(n^{-2}).$
%   ${\var} \left( {{T_{\bf{B}}}} \right) = {\var} \left( {{{\hat \tau }^2}} \right) - \frac{4}{n}\left\{ {\sum\limits_{j \in {\bf{B}}}^{} {\beta _j^4\left[ {E\left( {X_{ij}^4} \right) - 1} \right] + 2\sum\limits_{j \ne j' \in {\bf{B}}}^{} {\beta _j^2\beta _{j'}^2} } } \right\}  + O\left( {{n^{ - 2}}} \right)$
 Recall that by definition we have  $\widehat {\var\left( {{T_{\gamma}}} \right)} =  {\widehat {\var\left( {{{\hat \tau }^2}} \right)} - \frac{8}{n}\hat \tau _{{{\bf{B}}_{\gamma}}}^4}.$ Also recall that $\widehat{\var({\hat \tau }^2})$ is consistent by Proposition \ref{var_naive_est}.   Thus, it is enough to prove that
    $\hat \tau _{{{\bf{B}}_\gamma }}^2-\tau _{\bf{B}}^2\overset{p}{\rightarrow}0.$ 
 Now, since we assumed that  $n\left[ { P\left( \left\{ {{{\bf{B}}_\gamma } \neq {\bf{B}}} \right\} \right)} \right]^{1/2} \xrightarrow[n\rightarrow\infty]{}0$ then clearly $P\left( {{{\bf{B}}_\gamma} = {\bf{B}}} \right)\xrightarrow[n\rightarrow\infty]{}1.$ Thus,
   it is enough to show that    $\hat \tau _{{{\bf{B}} }}^2-\tau _{\bf{B}}^2\overset{p}{\rightarrow}0.$
     Recall that  $E(\hat\beta_j^2)=\beta_j^2$ for $j=1,...,p$ and notice that $\var(\hat\beta_j^2)\underset{n \to \infty}{\rightarrow}0$ by similar arguments that were used to derive \eqref{var_naive}. Hence, we have $\hat\beta_j^2-\beta_j^2\overset{p}{\rightarrow}0.$ Since  
   we assumed that ${\bf{B}}$ is finite, we have
   $$\hat \tau _{\bf{B}}^2- \tau _{\bf{B}}^2 = \sum\limits_{j \in {\bf{B}}}^{} {\left( {\hat \beta _j^2 - \beta _j^2} \right)}  \overset{p}{\rightarrow} 0,  $$
      and \eqref{sh2} follows.

 \begin{remark}\label{selection_algorithm}
We use the  the following simple selection algorithm $\gamma$:

\begin{algorithm}[H]\label{alg1}
\SetAlgoLined
\vspace{0.4 cm}

 \textbf{Input:}
 A dataset \(\left( {{{\bf{X}}_{n \times p}},{{\bf{Y}}_{n \times 1}}} \right)\).
\begin{enumerate}
  \item Calculate $\hat\beta_1^2,...,\hat\beta_p^2$ where  $\hat\beta_j^2$  is given in (\ref{beta_j_hat}) for $j=1,...,p.$   
  
  \item Calculate the  differences
  \({\lambda _j} = \hat \beta _{\left( j \right)}^2 - \hat \beta _{\left( {j - 1} \right)}^2\) for $j=2,\ldots,p$ where \(\hat \beta _{\left( 1 \right)}^2 < \hat \beta _{\left( 2 \right)}^2 < ... < \hat \beta _{\left( p \right)}^2\) denotes the  order statistics. 
  \item Select the covariates  \({{\bf{B}}_\gamma} = \left\{ {j:\hat \beta _{\left( j \right)}^2 > \hat \beta _{\left( {{j^*}} \right)}^2} \right\}\),  where \({j^*} = \mathop {\arg \max }\limits_j {\lambda _j}\).
  \end{enumerate}
\KwResult{Return $\bf{B_\gamma}$. }
  \caption{Covariate selection $\gamma$}
\end{algorithm}
The algorithm above finds the largest gap between the ordered estimated squared coefficients and then uses this gap as a threshold to select a set of  coefficients $\bf{B_\gamma} \subset \left\{ {1,...,p} \right\}.$ The algorithm  works well in  scenarios where a relatively large gap truly separates  between  larger coefficients  and the  smaller coefficients of the vector $\beta$.
 \end{remark}

\newpage
\vskip .65cm
\noindent
Ilan Livne (ilan.livne@campus.technion.ac.il)

\noindent
David Azriel (davidazr@technion.ac.il)

\noindent
Yair Goldberg (yairgo@technion.ac.il) 
\vskip 2pt

\noindent
The Faculty of Industrial Engineering and Management, Technion.
\end{document}